\documentclass[12pt]{article}

\usepackage{verbatim} 
\usepackage{amssymb,amsmath,psfrag}
\usepackage{pstricks,ifthen}
\usepackage{color}
\usepackage{bm}
\usepackage{graphicx}
\allowdisplaybreaks

\def\epsilon{\varepsilon}
\def\hat{\widehat}

\newcommand{\va}{\underline{a}}
\newcommand{\vb}{\underline{b}}
\newcommand{\vd}{\underline{d}}

\newcommand{\vpsi}{\underline{\psi}}
\newcommand{\vphi}{\underline{\phi}}

\newcommand{\vq}{\underline{q}}
\newcommand{\vQ}{\underline{Q}}
\newcommand{\vB}{\underline{B}}
\newcommand{\vP}{\underline{P}}

\newcommand{\vI}{\underline{I}}

\newcommand{\vv}{\underline{v}}

\newcommand{\vx}{\underline{x}}

\newcommand{\vz}{\underline{z}}
\newcommand{\vZ}{\underline{Z}}

\newcommand{\vR}{\underline{R}}
\newcommand{\vzero}{\underline{0}}
\newcommand{\vV}{\underline{V}}

\newcommand{\vVh}{\underline{V}^h}

\newcommand{\del}{\underline{\nabla}}
\newcommand{\eps}{\varepsilon}
\newcommand{\ts}{\textstyle}
\newcommand{\ds}{\displaystyle}

\newcommand{\beq}{\begin{equation}}
\newcommand{\eeq}{\end{equation}}

\renewcommand{\theequation}{\arabic{section}.\arabic{equation}}

\newcounter{ind}
\def\eqlabstart{%
\setcounter{ind}{\value{equation}}\addtocounter{ind}{1}%
\setcounter{equation}{0}%
\renewcommand{\theequation}{\arabic{section}.\arabic{ind}\alph{equation}}}
\def\eqlabend{%
\renewcommand{\theequation}{\arabic{section}.\arabic{equation}}%
\ifthenelse{\value{equation}=0}{\addtocounter{ind}{-1}}{}%
\setcounter{equation}{\value{ind}}}

\newcommand{\veta}{\underline{\eta}}
\newcommand{\vc}{\underline{c}}
\newcommand{\vS}{\underline{S}}
\newcommand{\vG}{\underline{G}}

\def\proof{\par{\it Proof}. \ignorespaces}
\def\endproof{\vbox{\hrule height0.6pt\hbox{%
\vrule height1.3ex width0.6pt\hskip0.8ex%
\vrule width0.6pt}\hrule height0.6pt}}

\newtheorem{thm}{\sc Theorem}[section]
\newtheorem{lem}{\sc Lemma}[section]

\newtheorem{rem}{\sc Remark}[section]

\newenvironment{AMS}%
{{\upshape\bfseries AMS subject classifications. }\ignorespaces}{}
\newenvironment{keywords}{{\upshape\bfseries Key words. }\ignorespaces}{}

\begin{document}
\title{
A Quasi-Variational Inequality Problem\\ Arising in 
the Modeling of Growing  Sandpiles}
\author{John W. Barrett\footnotemark[2] \quad {\small and} \quad 
        Leonid Prigozhin\footnotemark[3]}

\renewcommand{\thefootnote}{\fnsymbol{footnote}}
\footnotetext[2]{Department of Mathematics, 
Imperial College London, London, SW7 2AZ, UK.}
\footnotetext[3]{Department of Solar Energy and Environmental Physics,
Blaustein Institutes for Desert Research,
Ben-Gurion University of
the Negev, Sede Boqer Campus, 84990 Israel.}

\date{}

\maketitle

\begin{abstract}

Existence of a solution to the quasi-variational inequality problem arising in 
a model for sand surface 
evolution has been an open problem for a long time.
Another long-standing open problem concerns determining the dual variable,  
the flux of sand pouring down 
the evolving sand surface,
which is also of practical 
interest in a variety of applications of this model.
Previously, these problems were solved for the special case 
in which the inequality is simply variational. 
Here, we introduce a regularized mixed formulation involving both the primal 
(sand surface) and dual (sand flux) variables. We derive, analyse and   
compare two methods for the approximation, and numerical solution, of this mixed problem. 
We prove subsequence convergence of both approximations, as the mesh discretization
parameters tend to zero; 
and hence prove     
existence of a solution to this mixed model and the associated 
regularized quasi-variational inequality problem.
One of these numerical approximations, 
in which the flux is approximated by the divergence-conforming 
lowest order Raviart--Thomas element, leads to an efficient algorithm to compute not 
only the evolving pile surface, but also the flux of pouring sand. Results of 
our numerical experiments confirm the validity of the regularization employed.
\end{abstract}


\vspace{2mm}

\begin {keywords} Quasi-variational inequalities; critical-state problems; 
primal and mixed formulations; finite elements; existence; convergence analysis. \end{keywords}

\renewcommand{\thefootnote}{\arabic{footnote}}

\section{Introduction}
\label{intro}
\setcounter{equation}{0}

Let a cohesionless granular material (sand), characterized by its angle of repose $\alpha$, 
be poured out onto a rigid surface $y=w_0(\vx)$, where $y$ is vertical, $\vx\in\Omega\subset
 \mathbb{R}^d$, $d=1$ or 2, and $\Omega$ is a domain with boundary $\partial \Omega$. 
 The support surface
$w_0\in W^{1,\infty}_0(\Omega)$
and the nonnegative density of the distributed source  $f \in L^2(0,T;L^2(\Omega))$ are given. 
We consider the growing sandpile $y=w(\vx,t)$
and set an open boundary condition $w|_{\partial\Omega}=0$. 
Denoting by $\vq(\vx,t)$ the horizontal projection of the flux of material pouring 
down the evolving pile surface, 
we can write the mass balance equation
\beq \frac{\partial w}{\partial t} +\del \,.\, \vq=f.\label{bal}\eeq
The quasi-stationary model of sand surface evolution, see Prigozhin \cite{Fill,River,Psand},  
assumes the flow 
of sand is confined to a thin surface layer and directed towards the steepest descent of 
the pile surface. Wherever the support surface is covered by sand, the pile slope should not 
exceed the critical value; that is, $w>w_0\ \Rightarrow\ |\del w|\le k_0$, 
where $k_0=\tan\alpha$ 
is the internal friction coefficient. Of course, 
the uncovered parts of support can be steeper. 
This model does not allow for any flow on the subcritical parts of the pile surface; that is,
$|\del w|< k_0\ \Rightarrow\ \vq=\underline{0}$.
These constitutive relations can be conveniently reformulated 
for a.e.\ $(\vx,t) \in \Omega \times (0,T)$ as
\begin{align}|\del w| \leq M(w)
\qquad \mbox{and} \qquad 
M(w)\,|\vq| +\del w\,.\,\vq =0,\label{cond1}
\end{align}
where for any $\vx \in \overline{\Omega}$
\begin{align}
M(w)(\vx)&:=\left\{ \begin{array}{ll}  k_0 \qquad &w(\vx)>w_0(\vx),\\
\max(k_0,|\del w_0(\vx)|) \qquad & w(\vx)\le w_0(\vx).\end{array}\right.
\label{M}
\end{align}

Let us define, for any $\eta\in C(\overline{\Omega})$, the closed convex non-empty set
\beq K(\eta):=\left\{ \varphi\in W^{1,\infty}_0(\Omega)\ :\ |\del\varphi|\le M(\eta)\ \ 
\mbox{a.e. in}\ \Omega\right\}. \label{Kset}\eeq
Since $M(w)\,|\vq|+\del \varphi \,.\, \vq \geq 0$ for any $\varphi\in K(w)$,
we have, on noting (\ref{cond1}), that $w \in K(w)$ and $\del(\varphi-w)\,.\,\vq\geq 0$. 
A weak form of the latter inequality is:
for a.a. $t \in (0,T)$
\beq 
\int_{\Omega} \del\,.\,\vq\,(w-\varphi)\, {\rm d}\vx \geq 0 
\qquad \forall \ \varphi \in K(w).
\label{sbd1}
\eeq  
Combining (\ref{sbd1}) and (\ref{bal}) yields an evolutionary quasi-variational
inequality for the evolving pile surface: Find $w \in K(w)$ such that 
for $a.a.\ t \in (0,T)$
\beq 
\int_{\Omega} \left(\frac{\partial w}{\partial t} - f \right)(\varphi-w)\, 
{\rm d}\vx \geq 0 
\qquad \forall \ \varphi \in K(w).
\label{P1}
\eeq
Assuming there is no sand  on the support initially,  we set
\beq w(\cdot,0)=w_0(\cdot). 
\label{ic}\eeq

A solution $w$ to this quasi-variational inequality problem, (\ref{P1}) and (\ref{ic}), 
if it exists, should be 
a monotonically 
non-decreasing function in time for any
$f\geq 0$, see \S 3 in Prigozhin \cite{Psand}. 
However, existence and uniqueness of a solution has only been 
proved for  support surfaces with no steep slopes; that is, $|\del w_0|\leq k_0$, see 
Prigozhin \cite{Fill,Psand}.
In this case $K(w)\equiv K:=\left\{ \varphi\in W^{1,\infty}_0(\Omega)\ :\ 
|\del\varphi|\le k_0\ \ \mbox{a.e. in}\ \Omega\right\}$ and the quasi-variational inequality 
becomes simply a variational inequality.
Independently, the variational inequality for supports without 
steep slopes has been derived 
and studied in Aronson, Evans and Wu \cite{AEW} as the $p\rightarrow \infty$ 
limit of the evolutionary $p$-Laplacian equation.

The quasi-variational inequality (\ref{P1}) can, of course, 
be considered not only with the initial condition 
(\ref{ic}). However, if 
$w(\cdot,0) = \widetilde{w}_0(\cdot) \geq w_0(\cdot)$ and $\widetilde{w}_0$ does not 
belong to the admissible set $K(\widetilde{w}_0)$, an instantaneous solution reconstruction 
takes place. Such discontinuous solutions, interpreted as simplified descriptions of 
collapsing piles with overcritical slopes, were studied in the variational inequality case 
in Evans, Feldman and Gariepy \cite{EFG}, and 
Dumont and Igbida \cite{DI}. Since we assumed the initial condition (\ref{ic}) and, obviously, 
$w_0\in K(w_0)$, one could expect a solution continuously evolving in time. 
However, for the quasi-variational inequality with the open boundary condition 
$w|_{\partial\Omega}=0$, an uncontrollable influx of material from outside can occur 
through the parts of the boundary where $\del w_0\,.\,\underline{\nu}\geq k_0$,
where $\underline{\nu}$ is the outward unit normal to $\partial \Omega$.  
This makes the solution non-unique and, possibly, discontinuous. 
Such an influx is prevented in our model by assuming that 
\beq\del w_0\,.\,\underline{\nu}< k_0\ \ \mbox{on}\ 
\partial\Omega,\label{influx}\eeq 
which implies that $\del w\,.\,\underline{\nu}< k_0$ on  
$\partial\Omega$ for $t>0$.

For the variational inequality version of the sand model,
equivalent dual and mixed variational formulations have recently been proposed; 
see, e.g., Barrett and Prigozhin \cite{Dual} and Dumont and Igbida \cite{DIa}. 
Such formulations are often advantageous, 
because they allow one to determine not only the evolving sand surface $w$ but also the 
surface flux $\vq$, which is of interest too in various applications; see Prigozhin \cite{PChEJ,River},
and Barrett and Prigozhin \cite{BP5}.
In such formulations, and this is their additional advantage, the difficult to deal with 
gradient constraint is replaced by a simpler, although non-smooth, nonlinearity.

Here we will also use a mixed variational formulation of a regularized version of the growing 
sandpile model involving both variables. Instead of excluding the surface flux $\vq$ from 
the model formulation, as in the transition to (\ref{P1}) above, we now note that the first 
condition in 
(\ref{cond1}) holds if for a.e. $(\vx,t) \in \Omega \times (0,T)$
\beq M(w)\,|\vv|+\del w\,.\,\vv\geq 0\label{cond3a}\eeq
for any test flux $\vv$.
Hence we 
can reformulate the conditions 
(\ref{cond1}) for a.a.\ $t \in (0,T)$ as
\beq \int_{\Omega}\left[M(w)\,(|\vv|-|\vq|)-w\,\del\,.\,(\vv-\vq)\right] {\rm d}\vx\geq 0\label{cond3}\eeq
for any test flux $\vv$, and consider a mixed formulation of the sand model 
as (\ref{bal}) and (\ref{cond3}).

The quasi-variational inequality (\ref{P1}) is a difficult problem; in particular, due to 
the discontinuity of the nonlinear operator $M$, which determines the gradient constraint 
in (\ref{Kset}). Furthermore, the natural function space for the flux $\vq$ is the 
space of vector-valued bounded Radon measures having $L^2$ divergence. If $\vq$ is such a 
measure,  the discontinuity of $M(w)$ also makes it  difficult to give a sense to 
the term $\int_{\Omega}M(w)\, |\vq| \,{\rm d}\vx $ 
in the inequality (\ref{cond3}) of the mixed formulation.

In this work we consider a regularized version of the growing sandpile model 
with a  
continuous operator $M_\epsilon : C(\overline{\Omega}) \rightarrow C(\overline{\Omega})$, 
determined as follows.
For a fixed small $\varepsilon >0$, we approximate the initial data
$w_0\in W^{1,\infty}_0(\Omega)$ by $w_0^\epsilon \in W^{1,\infty}_0(\Omega) \cap
C^1(\overline{\Omega})$, and
$M(\cdot)$ by the continuous function $M_\eps(\cdot)$ such that 
for any $\vx \in \overline{\Omega}$
\begin{align}
&M_{\varepsilon}(\eta)(\vx)
\nonumber \\
&\; :=\left\{ \begin{array}{ll}  k_0 \quad& \eta(\vx)\geq w_0^\epsilon(\vx)
+ \varepsilon,\\
k_1^\epsilon(\vx) + (k_0-k_1^\epsilon(\vx)) \, \left(\ds\frac{\eta(\vx)-w_0^\epsilon(\vx)}{\varepsilon}\right)
\quad & \eta(\vx) \in [w_0^\epsilon(\vx),w_0^\epsilon(\vx)+\varepsilon],\\[4mm]
k_1^\epsilon(\vx)
:= \max(k_0,|\del w_0^\epsilon(\vx)|)
\quad & \eta(\vx)\le w_0^\epsilon(\vx).\end{array}\right.
\label{Meps}
\end{align}
We note that $M_\epsilon$ is such
that for all $\eta_1,\,\eta_2 \in C(\overline{\Omega})$
\begin{align}
|M_{\epsilon}(\eta_1)-M_{\epsilon}(\eta_2)|_{0,\infty,\Omega}
\leq \frac{k_{1,\infty}^\epsilon-k_0}{\varepsilon}\,
|\eta_1-\eta_2|_{0,\infty,\Omega}\,,
\label{Mepscont}
\end{align}
where
\begin{align}
k_{1,\infty}^\epsilon := \max_{x \in \overline{\Omega}} k_1^\epsilon(\vx)\,.
\label{k1max}
\end{align}
In addition, it follows for any $\vx \in \overline{\Omega}$ that
\begin{align}
\eta_1(\vx) \geq \eta_2(\vx) \quad \Rightarrow \quad 0 < k_0 \leq M_\epsilon(\eta_1(\vx))
\leq M_\epsilon(\eta_2(\vx)) \leq k^\epsilon_1(\vx)\,.
\label{Mepsineq}
\end{align}

We note that the analysis of the sand quasi-variational inequality problem studied in this paper
is far more involved than that of the superconductivity quasi-variational inequality problem
studied by the present authors in \cite{BP5}. In the superconductivity context, 
${M} :
{\mathbb R} \rightarrow [{M}_0,{M}_1] \subset {\mathbb R}$
with ${M}_0>0$.    
In \cite{BP5}, we exploit the fact that $|\del w(\vx)| \leq M(w(\vx))$ can be rewritten as
$|\del [F(w(\vx))]| \leq 1$ for all $\vx \in \Omega$, where $F'(\cdot)=
[{M}(\cdot)]^{-1}$ and $F(0)=0$.    
Clearly, such a reformulation is not applicable to $M(\cdot)$, (\ref{M}), 
or $M_\epsilon(\cdot)$, (\ref{Meps}).

In addition, we note that in the very recent preprint by Rodrigues and Santos \cite{RS}
an existence result can be deduced for the primal quasi-variational inequality problem 
(\ref{P1}) for a continuous and positive $M(\cdot)$, such as $M_\eps(\cdot)$, and $f \in 
W^{1,\infty}(\Omega \times (0,T))$. Assuming $w^0 \in K(w^0)\cap C_0(\overline{\Omega})$,
they show that $w \in L^\infty(0,T;W^{1,\infty}_0(\Omega)) \cap W^{1,\infty}(0,T;
[C_0(\overline{\Omega})]^\star)$. 
Their proof is based on the method of vanishing viscosity and constraint penalization.  

The outline of this paper is as follows. 
In the next section we introduce two fully practical finite element
approximations, (Q$^{h,\tau}_A$) and (Q$^{h,\tau}_{B,r}$), to the regularized 
mixed formulation (\ref{bal}) and (\ref{cond3}), where $M(\cdot)$ is replaced by 
$M_\epsilon(\cdot)$,
and prove well-posedness and stability bounds.
Here $h$ and $\tau$ are the spatial and temporal discretization
parameters, respectively. In addition, $r>1$ is a regularization parameter
in replacing the non-differentiable nonlinearity $|\cdot|$ by the strictly convex 
function $\frac{1}{r}\,|\cdot|^r$.   
The approximation (Q$^{h,\tau}_A$) is based on a continuous piecewise linear approximation
for $w$ and a piecewise constant approximation for $\vq$, whereas
(Q$^{h,\tau}_{B,r}$) is based on a piecewise constant approximation for $w$
and the lowest order Raviart--Thomas element for $\vq$. 
In Section \ref{conv} we prove subsequence convergence of both approximations to 
a solution of a weak formulation of the regularized
mixed problem. This is achieved by passing to the limit $h \rightarrow 0$ first,
then $r \rightarrow 1$ in the case of (Q$^{h,\tau}_{B,r}$),
and finally $\tau \rightarrow 0$.  
In Section \ref{snas}, we introduce iterative algorithms for solving the resulting 
nonlinear algebraic equations arising from both approximations at each time level.  
Finally, in Section \ref{numexpts} we present various numerical experiments.
These show that only the approximation (Q$^{h,\tau}_{B,r}$) leads to an efficient algorithm
to approximate both the surface $w$ and the flux $\vq$.


We end this section with a few remarks about the notation employed in this paper.
Above and throughout we adopt the standard notation for Sobolev spaces
on a bounded domain $D$ with a Lipschitz boundary,
denoting the norm of
$W^{\ell,s}(D)$ ($\ell \in {\mathbb N}$, $s\in [1, \infty]$) by
$\|. \|_{\ell,s,D}$ and the semi-norm by $|\cdot |_{\ell,s,D}$.
Of course, we have that
$|\cdot |_{0,s,D} \equiv \|\cdot \|_{0,s,D}$.
We extend these norms and semi-norms in the natural way to the corresponding
spaces of vector 
functions.
For $s=2$, $W^{\ell,2}(D)$ will be denoted by
$H^\ell(D)$ with the associated norm and semi-norm written
as, respectively, $\|\cdot\|_{\ell,D}$ and $|\cdot|_{\ell,D}$.
We set $W^{1,s}_0(D):= \{\eta \in W^{1,s}(D) : \eta = 0 \mbox{ on } \partial D \}$,
and $H^1_0(D) \equiv W^{1,2}_0(D)$.
We recall the Poincar\'{e} inequality for any $s \in [1,\infty]$
\begin{align}
|\eta|_{0,s,D} \leq C_\star(D)\,|\del \eta|_{0,s,D}\qquad \forall \ \eta \in W^{1,s}_0(D)\,,
\label{Poin}
\end{align}
where the constant $C_\star(D)$ depends on $D$, but is independent of $s$;
see e.g.\  page 164 in Gilbarg and Trudinger \cite{GT}.
In addition, $|D|$ will denote the measure of $D$
and $(\cdot,\cdot)_D$ the standard inner product on $L^2(D)$.
When $D \equiv \Omega$, for ease of notation we write $(\cdot,\cdot)$ for
$(\cdot,\cdot)_\Omega$.

For $m\in \mathbb N$, let 
(i) $C^m(\overline{D})$ 
denote the space of continuous functions with 
all derivatives up to order $m$ 
continuous on $\overline{D}$, 
(ii) $C^m_0(D)$ denote the space of continuous functions with 
compact support in $D$ with all derivatives up to order $m$ 
continuous on $D$ and (iii) 
$C^m_0(\overline{D})$ denote those functions in $C^m(\overline{D})$ 
which vanish on $\partial D$. 
In the case $m=0$, we drop the superscript $0$ for all three spaces.

As one can identify $L^{1}(D)$ as a closed subspace of
the Banach space of bounded Radon measures,
${\cal M}(\overline{D}) \equiv [C(\overline{D})]^\star$, i.e.\
the dual of $C(\overline{D})$;
it is convenient to adopt the notation 
\begin{align}
\int_{\overline{D}} |\mu| \equiv \|\mu\|_{{\cal M}(\overline{D})}
:= \sup_{\stackrel{\eta \in C(\overline{D})}{|\eta|_{0,\infty,D}\leq 1}}
\langle \mu,\eta \rangle_{C(\overline{D})}
< \infty,
\label{Mnorm2}
\end{align}
where $\langle \cdot,\cdot \rangle_{C(\overline{D})}$
denotes the duality pairing on $ [C(\overline{D})]^\star
\times C(\overline{D}).$

We introduce also the Banach spaces for a given $s \in [1,\infty]$
\eqlabstart
\begin{align}
\vV^{s}(D)
&:= \{\vv \in [L^s(D)]^d :
\del\,.\,\vv\in L^2(D)\}\,,
\label{vVs}\\
\vV^{\cal M}(D)
&:= \{\vv \in [{\mathcal M}(\overline{D})]^d :
\del\,.\,\vv\in L^2(D)\}\,.
\label{vVM}
\end{align}
\eqlabend
The condition $\del\,.\,\vv\in L^2(D)$ in (\ref{vVs},b) means that there exists
$u\in L^2(D)$ such that $\langle \vv,\del \phi \rangle_{C(\overline{D})}
=-(u,\phi)_D$ for any $\phi\in C^1_0(D)$.

We note that if
$\{\mu_n\}_{n \geq 0}$ is a bounded sequence in ${\cal M}(\overline{D})$,
then there exist a subsequence $\{\mu_{n_j}\}_{n_j \geq 0}$ and a
$\mu \in {\cal M}(\overline{D})$ such that as $n_j \rightarrow \infty$
\begin{align}
\mu_{n_j} \rightarrow \mu \quad \mbox{vaguely in } {\cal M}(\overline{D});
\quad \mbox{i.e.}
\quad \langle \mu_{n_j}-\mu,\eta\rangle_{C(\overline{D})} \rightarrow 0\quad
\forall \ \eta \in C(\overline{D})\,.
\label{Mweak1}
\end{align}
In addition,
we have that
\begin{align}
\liminf_{n_j \rightarrow \infty} \int_{\overline{D}} |\mu_{n_j}|
\geq \int_{\overline{D}} |\mu|\,;
\label{Mweak2}
\end{align}
see e.g.\ page 223 in 
Folland \cite{Folland}.

We note that $\vV^{\cal M}(D)$ and $\vV^s(D)$ for $s \in [1,2)$
are not of local type; that is, $\vv \in \vV^{\cal M}(D)\; [\vV^s(D)]$ 
and
$\phi \in C^{\infty}(\overline{D})$
does not imply that
$\phi\,\vv \in \vV^{\cal M}(D)\; [\vV^{s}(D)]$,
see e.g.\ page 22 in Temam \cite{Temam85}.
Therefore, one has to avoid cut-off functions in proving any
required density results.
If $\Omega$ is %
convex, which we shall assume for the analysis in this
paper, then it is strictly star-shaped and one can show,
using the standard techniques
of change of variable and mollification,
that
\begin{align}
[C^\infty(\overline{\Omega})]^2 \quad
\mbox{ is dense in } \quad \vV^s(\Omega)\quad
\mbox{ if } \quad s \in (1,\infty).
\label{densevVs}
\end{align}
Moreover, 
for any $\vv \in \vV^{\cal M}(\Omega)$, there exist $\{\vv_j\}_{j \geq 1}
\in  [C^\infty(\overline{\Omega})]^d$
such that
 \eqlabstart
\begin{alignat}{3}
\vv_{j}
&\rightarrow \vv \qquad &&\mbox{vaguely in }
[{\cal M}(\overline{\Omega})]^d
\qquad &&\mbox{as } j \rightarrow \infty\,,
\label{weakj}
\\
\del\,.\,\vv_{j}
&\rightarrow \del\,.\,\vv \qquad
&&\mbox{weakly in }
L^2(\Omega)\qquad &&\mbox{as } j \rightarrow \infty\,,\label{divweakj}\\
\limsup_{j \rightarrow \infty} \int_{\Omega} \rho\,|\vv_{j}|\,{\rm d}\vx
&= \int_{\overline{\Omega}} \rho \,|\vv|&&&&
\label{normj}
\end{alignat}
\eqlabend
for  any positive $\rho \in C(\overline{\Omega})$.
We briefly outline the proofs of 
(\ref{weakj}--c).
Without loss of generality, one can assume that
$\Omega$ is strictly star-shaped with respect to the origin.
Then for $\vv$ defined on $\Omega$ and $\theta>1$,
we have that $\vv_\theta(\vx)=
\vv(\theta^{-1}\vx)$ is defined on
$\Omega_{\theta}:=
\theta\,\Omega \supset \Omega$.
Applying standard Friedrich's mollifiers $J_\varepsilon$
to $\vv_\theta$, and a diagonal subsequence argument yield,
for $\theta \rightarrow 1$ and $\varepsilon \rightarrow 0$ as $j \rightarrow \infty$,
the desired sequences $\{\vv\}_{j\geq 1}$
demonstrating (\ref{densevVs}) if $\vv \in \vV^s(\Omega)$
and satisfying (\ref{weakj}--c) if $\vV^{\cal M}(\Omega)$;
see e.g.\ Lemma 2.4 in Barrett and Prigozhin \cite{BP3}, where such techniques
are used to prove similar density results.

We recall the following Sobolev interpolation theorem, see Theorem 5.8 in 
Adam and Fournier \cite{AF}. If $\eta \in W^{1,s}(D)$, with $s>d$, 
then $\eta \in C(\overline{\Omega})$ with the embedding being compact; and moreover,
\begin{align}
|\eta|_{0,\infty,D} \leq C(s,D)\, \|\eta\|_{1,s,D}^{\alpha}\, |\eta|_{0,D}^{1-\alpha}
\quad \mbox{with} \quad \alpha = \frac{d\,s}{d\,s + 2(s-d)} \in (0,1)\,. 
\label{Sobinterp}
\end{align}    
We recall also
the Aubin--Lions--Simon compactness theorem, see Corollary 4 in Simon \cite{Simon}. 
Let ${\cal B}_0$, ${\cal B}$ and
${\cal B}_1$ be Banach
spaces, ${\cal B}_i$, $i=0,1$, reflexive, with a compact embedding ${\cal B}_0
\hookrightarrow {\cal B}$ and a continuous embedding ${\cal B} \hookrightarrow
{\cal B}_1$. Then, for $\alpha>1$, the embedding
\begin{align}
&\{\,\eta \in L^{\infty}(0,T;{\cal B}_0): \frac{\partial \eta}{\partial t}
\in L^{\alpha}(0,T;{\cal B}_1)\,\} \hookrightarrow C([0,T];{\cal B})
\label{compact1}
\end{align}
is compact.

Finally, throughout $C$ denotes a generic positive constant independent of
the regularization parameter, $r \in (1,\infty)$,
the mesh parameter $h$ and the time step
parameter $\tau$. Whereas, $C(s)$ denotes a positive constant dependent on
the parameter $s$.

\section{Finite Element Approximation}
\label{fea}
\setcounter{equation}{0}

We make the following assumptions on the data.
\vspace{2mm}

\noindent
{\bf (A1)} $\Omega \subset {\mathbb R}^d$, $d=1$ or $2$,
is convex having a 
Lipschitz boundary $\partial \Omega$
with outward unit normal $\underline{\nu}$.
$f \in L^2(0,T;L^2(\Omega))$ is a nonnegative source,
and $M_\epsilon(\cdot)$ is 
given by (\ref{Meps}). 
In addition, 
the initial data $w_0^\epsilon \in W^{1,\infty}_0(\Omega)\cap 
C^1(\overline{\Omega})$ is such that 
$\del w_0^\epsilon \,.\, \underline{\nu} < k_0$.

\vspace{2mm}

For ease of exposition, we shall assume that $\Omega$ is a polygonal domain to
avoid perturbation of domain errors in the finite element approximation. We make
the following standard assumption on the partitioning. 
\vspace{2mm}

\noindent
{\bf (A2)} $\Omega$ is polygonal.
Let
 $\{{\cal T}^h\}_{h>0}$ be a regular family
 of partitionings of $\Omega$
 into disjoint open
 simplices $\sigma$ with $h_{\sigma}:={\rm diam}(\sigma)$
 and $h:=\max_{\sigma \in {\cal T}^h}h_{\sigma}$, so that
 $\overline{\Omega}=\cup_{\sigma\in{\cal T}^h}\overline{\sigma}$.

\vspace{2mm}
Let $\underline{\nu}_{\partial \sigma}$ be the outward unit normal to $\partial \sigma$,
the boundary of $\sigma$.
We then introduce the following finite element spaces
\eqlabstart
\begin{align}
S^h
&:= \{ \eta^h \in L^\infty(\Omega) : \eta^h \mid_{\sigma} = a_{\sigma} \in
{\mathbb R} \quad \forall \ \sigma \in {\cal T}^h \,\}\,,
\label{Shc}\\
S^h_{\geq 0}
&:= \{ \eta^h \in L^\infty(\Omega) : \eta^h \mid_{\sigma} = a_{\sigma} \in
{\mathbb R}_{\geq 0} \quad \forall \ \sigma \in {\cal T}^h \,\}\,,
\label{Sh+}\\
\vS^h
&:= \{ \veta^h \in [L^\infty(\Omega)]^d : \veta^h \mid_{\sigma} = \va_{\sigma} \in
{\mathbb R}^d \quad \forall \ \sigma \in {\cal T}^h \,\}\,,
\label{Sh}\\
{U}^h
&:= \{ \eta^h \in C(\overline{\Omega}) : \eta^h \mid_{\sigma} 
= \va_{\sigma}\,.\,\vx + b_{\sigma},
\ \va_{\sigma} \in {\mathbb R}^d, \ b_{\sigma} \in {\mathbb R}
\quad \forall \ \sigma \in {\cal T}^h \}
\label{Wh}\,,\\
U^h_0 &:= U^h \cap W^{1,\infty}_0(\Omega)\,,
\label{Wh0}\\
\underline{V}^h
&:= \{ \vv^h \in [L^\infty(\Omega)]^d : \vv^h \mid_{\sigma} =
\va_{\sigma} + b_{\sigma}\,\vx,
\ \va_{\sigma} \in {\mathbb R}^d, \ b_{\sigma} \in {\mathbb R}
\quad \forall \ \sigma \in {\cal T}^h 
\nonumber \\
& \hspace{0.48in} \mbox{and }
(\vv^h \mid_{\sigma} 
-  \vv^h \mid_{\sigma'}) \,.\,\underline{\nu}_{\partial \sigma}
= 0 \; \mbox{ on }  \partial \sigma \cap \partial \sigma'
\quad \forall \ \sigma,\,\sigma' \in {\cal T}^h
\}
\label{RT}\,.
\end{align}
\eqlabend
Here $\vV^h$ is the lowest order Raviart--Thomas finite element space.

Let $\pi^h : C(\overline{\Omega})  \rightarrow U^h$ denote the interpolation 
operator such that \linebreak $\pi^h \eta(\vx_j)= \eta(\vx_j)$, $j=1 \rightarrow J$,
where $\{\vx_j\}_{j=1}^J$ are the vertices of the partitioning ${\cal T}^h$.  
We note for $m=0$ and 1 that 
\eqlabstart
\begin{alignat}{2}
&|(I-\pi^h) \eta|_{m,s,\sigma} \leq C\,h^{2-m}\,|\eta|_{2,s,\sigma}
\qquad &&\forall \ \sigma \in {\cal T}^h, \quad \mbox{for any } s \in [1,\infty]\,,
\label{piconv0}\\
&\lim_{h \rightarrow 0}
\|(I-\pi^h) \eta\|_{m,\infty,\Omega} = 0
\qquad &&\forall \ \eta \in C^m(\overline{\Omega})
\,;
\label{piconv}
\end{alignat}
\eqlabend
where $I$ is the identity operator.
Let $\vP^h: [L^1(\Omega)]^d \rightarrow \vS^h$ be such that
\begin{align}
\vP^h \vv \mid_{\sigma} = \frac{1}{|\sigma|} \int_{\sigma} \vv \,{\rm d}\vx   \qquad
\forall \ \sigma \in {\cal T}^h\,.
\label{Ph}
\end{align}  
We note that 
\eqlabstart
\begin{align}
&|\vP^h \vv|_{0,s,\sigma} \leq |\vv|_{0,s,\sigma} \qquad \forall\ \vv \in [L^s(\sigma)]^d,
\quad s \in [1,\infty],\quad \forall \ \sigma \in {\cal T}^h\,,
\label{Phstab}\\
&\lim_{h \rightarrow 0}
|\,|\vv|-|\vP^h\vv|\,|_{0,\infty,\Omega}
\leq 
\lim_{h \rightarrow 0}
|\vv-\vP^h \vv|_{0,\infty,\Omega} = 0
\qquad \forall \ \vv \in [C(\overline{\Omega})]^d\,.
\label{Pconv}
\end{align}
\eqlabend
Similarly, we define $P^h : L^1(\Omega) \rightarrow S^h$ with the equivalent
to (\ref{Phstab},b) holding.

In addition, we introduce
the generalised interpolation operator
$\vI^h : 
[W^{1,s}(\Omega)]^d \rightarrow \vVh$, where $s>1$,
satisfying
\begin{align}
\int_{\partial_i \sigma} (\vv- \vI^h\vv)\,.\,\underline{\nu}_{\partial_i \sigma}\,ds =0
\qquad i=1 \rightarrow 3, \quad \forall \ \sigma \in {\cal T}^h\,;
\label{interpcon}
\end{align}
where $\partial \sigma \equiv \cup_{i=1}^{3} \partial_i \sigma$
and $\underline{\nu}_{\partial_i\sigma}$ are the corresponding outward unit normals
on $\partial_i \sigma$.
It follows that
\begin{align}
(\del\,.\,(\vv-\vI^h\vv), \eta^h)=0 \qquad \forall \ \eta^h \in S^h.
\label{interpint}
\end{align}
Moreover,
we have for all $\sigma \in {\cal T}^h$ and any $s \in (1,\infty]$ that
\begin{align}
||\vv| -|\vI^h \vv| |_{0,s,\sigma} \leq  
|\vv -\vI^h \vv |_{0,s,\sigma} &\leq C\,h_\sigma\,|\vv|_{1,s,\sigma}
\quad \mbox{and}
\quad
|\vI^h \vv|_{1,s,\sigma} \leq C \,|\vv|_{1,s,\sigma}
\,,
\label{interp1}
\end{align}
e.g.\ see Lemma 3.1 in Farhloul \cite{MF} and the proof
given there for $s\ge2$ is also valid for any $s \in (1,\infty]$.

We introduce 
$(\eta,\chi)^h :=\sum_{\sigma \in {\cal T}^h} (\eta,\chi)^h_{\sigma}$,
and 
\begin{align}
(\eta,\chi)^h_{\sigma} &:= \ts \frac{1}{d+1}\,|\sigma|\, \ds \sum_{j=1}^{d+1}
\eta(\vx_j^{\sigma})\,\chi(\vx_j^{\sigma}) = \int_{\sigma} \pi^h[\eta\,\chi]\,{\rm d}\vx
\nonumber \\
&\hspace{2in} \forall \ \eta,\,\chi \in C(\overline{\sigma}),
\quad \forall \ \sigma \in {\cal T}^h\,;
\label{ni}
\end{align}
where $\{\vx_j^{\sigma}\}_{j=1}^{d+1}$ are the vertices of $\sigma$.
Hence $(\eta,\chi)^h$ averages
the integrand $\eta\,\chi$
over each simplex $\sigma$
at its vertices,
and is exact if $\eta\,\chi$ is piecewise 
linear over the
partitioning ${\cal T}^h$.
We recall the well-known results that
\eqlabstart
\begin{align}
|\eta^h|_{0,\Omega}^2 \leq |\eta^h|_h^2 := (\eta^h,\eta^h)^h  
&\leq (d+2)\,
|\eta^h|_{0,\Omega}^2
\qquad \forall \ \eta^h \in U^h\,, 
\label{eqnorm}
\\[2mm]
\left|(\eta^h,\chi^h)-(\eta^h,\chi^h)^h \right|&=
\left|((I-\pi^h)(\eta^h\,\chi^h),1)\right| 
\leq |(I-\pi^h)(\eta^h\,\chi^h)|_{0,1,\Omega}
\nonumber \\
&\leq C\,h\,|\eta^h|_{0,\Omega}\,|\chi^h|_{1,\Omega}
\qquad \forall \ \eta^h,\,\chi^h \in U^h\,,
\label{nierr}
\end{align}
where we have noted (\ref{piconv0}).
\eqlabend

In order to prove existence of solutions to 
approximations of (\ref{cond3}),
we regularise
the non-differentiable nonlinearity $|\cdot|$ by
the strictly convex function
$ \frac{1}{r}\,|\cdot|^{r}$ for $r>1$.
We note
for all $\va, \,\vb \in {\mathbb R}^d$ that
\begin{align}
\frac{1}{r}\,\frac{\partial|\va|^r}{\partial a_i}  = |\va|^{r-2}\,a_i
\quad \Rightarrow
\quad
|\va|^{r-2}\,\va\,
.\,(\va-\vb) &\geq \ts \frac{1}{r}\,\left[\, |\va|^r -|\vb|^r\,\right]\,.
\label{aederiv}
\end{align}
Similarly to (\ref{eqnorm}), we have from the equivalence of norms and the convexity of $|\cdot|^r$
for any $r >
1$ and for any $\vv^h \in \vV^h$ that
\begin{align}
C \,(\,|\vv^h|^r,1)^h_\sigma
\leq
\int_{\sigma} |\vv^h|^r \,{\rm d}\vx \leq  (\,|\vv^h|^r,1)^h_\sigma
\qquad \forall \ \sigma \in {\cal T}^h\,.
\label{convex}
\end{align}
Furthermore,
it follows from  (\ref{ni}), (\ref{interp1}) and (\ref{aederiv})
for any $r>1$ and any $\sigma \in {\cal T}^h$ that
\begin{align}
|\,\int_{\sigma} |\vI^h \vv|^r -(|\vI^h\vv|^r,1)^h_\sigma|
&\leq
C\,h_\sigma\,|\sigma|\, |[\vI^h \vv]^r|_{1,\infty,\sigma} 
\leq C\,r\,h_\sigma\,|\sigma|\,\|\vv\|_{1,\infty,\sigma}^r
\nonumber \\
&\hspace{1.8in}
\forall \ \vv \in [W^{1,\infty}(\sigma)]^d.
\label{nih}
\end{align}

In addition, let
$0= t_0 < t_1 < \ldots < t_{N-1} < t_N = T$ be a
partitioning of $[0,T]$
into possibly variable time steps $\tau_n := t_n -
t_{n-1}$, $n=1\to N$. We set
$\tau := \max_{n=1\to N}\tau_n$ and introduce
\begin{align}
f^n(\cdot) := \frac{1}{\tau_n} \int_{t_{n-1}}^{t_n} f(\cdot,t) \,{\rm d}t
\in L^2(\Omega) \qquad n=1 \rightarrow N\,.
\label{fn}
\end{align}
We note  that
\begin{align}
\sum_{n=1}^N \tau_n \,|f^n|_{0,s,\Omega}^s 
\leq \int_0^T |f|_{0,s,\Omega}^s\,{\rm d}t
\qquad \mbox{for any } s \in [1,2]\,.
\label{fnsum}
\end{align}

Finally, on setting 
\begin{align}
w^{\eps,h}_0 = P^h [\pi^h w^\eps_0]\,,
\label{weps0h}
\end{align}
we introduce $M^h_\eps : S^h \rightarrow S^h$ approximating $M_\eps : C(\overline{\Omega})
\rightarrow C(\overline{\Omega})$, defined by (\ref{Meps}), for any $\sigma \in {\cal T}^h$ as
\begin{align}
&M_{\varepsilon}^h(\eta^h)
:=\left\{ \begin{array}{ll}  k_0 \quad& \eta^h\geq w_0^{\epsilon,h}
+ \varepsilon,\\
k_{1,\sigma}^{\epsilon,h} + (k_0-k_{1,\sigma}^{\epsilon,h}) \left(\ds\frac{\eta^h-w_0^{\epsilon,h}}{\varepsilon}\right)
\quad & \eta^h \in [w_0^{\epsilon,h},w_0^{\epsilon,h}+\varepsilon],\\[4mm]
k_{1,\sigma}^{\epsilon,h}
:= \max(k_0,|\del \pi^h w_0^\epsilon\mid_{\sigma}|)
\quad & \eta^h\le w_0^{\epsilon,h}.\end{array}\right.
\label{Meps_h}
\end{align}
We note that $M_\eps$ is also well-defined on $S^h$ with $M_\eps : S^h \rightarrow L^\infty(\Omega)$,
and we have the following result.
\begin{lem}\label{Mepshapproxlem}
For any $\eta^h \in S^h$, we have  that
\begin{align}
|M_\eps(\eta^h) - M^h_\eps(\eta^h)|_{0,\infty,\Omega} 
\leq C(\epsilon^{-1}) \left[\, |(I-P^h) w^{\eps}_0|_{0,\infty,\Omega}+
\|(I-\pi^h) w^{\eps}_0\|_{1,\infty,\Omega}\,
\right].
\label{Mepshapprox}
\end{align}
\end{lem}
\proof It is convenient to rewrite  (\ref{Meps}) and (\ref{Meps_h}) for any 
$\eta^h \in S^h$ and for a.e.\ $\vx \in \Omega$ as  
\eqlabstart
\begin{align}
M_\eps(\eta^h)(\vx) &= k_0 + \left(\ds \frac{k_{1}^\eps(\vx)-k_0}{\epsilon}\right) \min(\, 
\max( w^\eps_0(\vx)+\epsilon-\eta^h(\vx),0),\,\epsilon)\,, \label{Mepsmm} \\ 
M_\eps^h(\eta^h)(\vx) &= k_0 + \left(\ds \frac{k_{1}^{\eps,h}(\vx)-k_0}{\epsilon}\right) 
\min(\, 
\max( w^{\eps,h}_0(\vx)+\epsilon-\eta^h(\vx),0),\,\epsilon)\,; \label{Meps_hmm} 
\end{align}
\eqlabend
where
\begin{align}
k_0 \leq M_\epsilon^h (\eta^h)(\vx) \leq 
k^{\eps,h}_1(\vx) := \max(k_0,|\del \pi^h w^{\eps}_0(\vx)|)\quad \mbox{for a.e.} \ \vx \in \Omega\,.
\label{k1epshdef} 
\end{align}
Since  
\begin{align}
|\min(\max(a,0),\eps)-\min(\max(b,0),\eps)| &\leq |a-b| \nonumber \\
\mbox{and} \hspace{1.6in}  |\min(\max(a,0),\eps)| &\leq \epsilon \qquad \quad \forall
\ a,\,b \in {\mathbb R},
\label{minmax}
\end{align}
it follows from (\ref{Mepsmm},b), (\ref{k1epshdef}), (\ref{weps0h}), (\ref{Phstab}) 
and Assumption (A1) that
\begin{align}
&|M_\eps(\eta^h)-M_\eps^h(\eta^h)|_{0,\infty,\Omega} 
\nonumber \\
&\hspace{0.5in}\leq
\frac{k^{\eps}_{1,\infty} -k_0}{\epsilon} \, |w^\epsilon_0-w^{\epsilon,h}_0|_{0,\infty,\Omega}
+ |\,|\del w^\eps_0|-|\del \pi^h w^{\eps}_0|\,|_{0,\infty,\Omega} \nonumber \\
& \hspace{0.5in} \leq C(\epsilon^{-1}) \left[\, |(I-P^h) w^{\eps}_0|_{0,\infty,\Omega}+
\|(I-\pi^h) w^{\eps}_0\|_{1,\infty,\Omega}\,
\right]\,;
\label{MepsMeps_hdiff}
\end{align}
and hence the desired result (\ref{Mepshapprox}).
\endproof

\subsection{Approximation (Q$^{h,\tau}_{A}$)}
Our first fully practical finite element approximation is:

\noindent
{\bf (Q$^{h,\tau}_A$)}
For $n = 1 \rightarrow N$, find $W_A^n \in U^h_0$ and 
$\vQ_A^n \in \vS^h$ such that
\eqlabstart
\begin{alignat}{2}
\left(\frac{W_A^n-W_A^{n-1}}{\tau_n},\eta^h\right)^h - (\vQ_A^n,\del \eta^h) &= (f^n 
,\eta^h)
\quad \;&&\forall \ \eta^h \in U^h_0,\label{Qth1}\\
(M_\epsilon^h(P^h W_A^n), |\vv^h|-|\vQ_A^n|) + (\del W_A^n,\vv^h-\vQ_A^n)  &\geq0
\quad \;
&&\forall \ \vv^h\in \vS^h\,;
\label{Qth2}
\end{alignat}
\eqlabend
where $W_A^0 = \pi^h w_0^\epsilon$.

\vspace{2mm}

For any  $\chi^h \in U^h_0$, 
we introduce the closed convex non-empty set  
\begin{align}
K^h(\chi^h) &:= \{ \eta^h \in U^h_0 : |\del \eta^h| \leq M^h_\epsilon(P^h \chi^h) 
\quad \mbox{a.e.\ in } \Omega\}\,.   
\label{Kh}
\end{align}

In Theorem \ref{Qthstab} below, we will show that 
(Q$^{h,\tau}_A$), (\ref{Qth1},b), 
is equivalent to (P$^{h,\tau}_A$) and (M$^{h,\tau}_A$).
The former is the 
approximation of the primal 
quasi-variational inequality:

\noindent
{\bf (P$^{h,\tau}_A$)}
For $n = 1 \rightarrow N$, find $W_A^n \in K^h(W_A^n)$ such that
\begin{align}
&\left(\frac{W_A^n-W_A^{n-1}}{\tau_n},\eta^h-W_A^n\right)^h  
\geq (f^n ,\eta^h-W_A^n)
\qquad \forall \ \eta^h \in K^h(W_A^n)\,,\label{Pth}
\end{align}
where $W_A^0 = \pi^h w_0^\epsilon$.

The latter, having obtained $\{W_A^n\}_{n=1}^N$ from (P$^{h,\tau}_A$), 
is the minimization 
problem:

\noindent 
{\bf (M$^{h,\tau}_A$)}
For $n= 1 \rightarrow N$, find $\vQ_A^n \in \vZ^{h,n}$ such that
\begin{align}
(M^h_\epsilon(P^h W^n_A), |\vQ_A^n|) &\leq (M^h_\epsilon(P^h W_A^n), |\vv^h|)
\qquad \forall \ \vv^h \in \vZ^{h,n}\,,\label{Mth}
\end{align}
where
\begin{align}
\vZ^{h,n} &:= \Bigl\{ \vv^h \in \vS^h : 
(\vv^h, \del \eta^h) =\left(\frac{W_A^n-W_A^{n-1}}{\tau_n},\eta^h\right)^h-(f^n,\eta^h)
\nonumber \\
& \hspace{3in}\quad \forall \ \eta^h \in U^h_0\Bigr\}\,.   
\label{Zh}
\end{align}
As $\del U^h_0$ is a strict subset of $\vS^h$, it follows that the affine manifold $\vZ^{h,n}$,
$n=1 \rightarrow N$,
is non-empty.

We consider the following regularization 
of (Q$^{h,\tau}_A$)
for a given $r>1$:

\noindent
{\bf (Q$^{h,\tau}_{A,r}$)}
For $n = 1 \rightarrow N$, find $W^n_{A,r} \in U^h_0$ and
$\vQ^n_{A,r} \in \vS^h$ such that
\eqlabstart
\begin{alignat}{2}
\left(\frac{W^n_{A,r}-W^{n-1}_{A,r}}{\tau_n},\eta^h\right)^h 
- (\vQ^n_{A,r},\del \eta^h) &= (f^n 
,\eta^h)
\quad \;&&\forall \ \eta^h \in U^h_0,\label{Qthr1}\\
(M_\epsilon^h(P^h W^n_{A,r})\,|\vQ^n_{A,r}|^{r-2} \vQ^n_{A,r},
\vv^h) + (\del W^n_{A,r},\vv^h)  &=0
\quad \;
&&\forall \ \vv^h\in \vS^h\,;
\label{Qthr2}
\end{alignat}
\eqlabend
where $W^0_{A,r} = \pi^h w_0^\epsilon$.

Associated with (Q$^{h,\tau}_{A,r}$) is the corresponding 
approximation of a generalised $p$-Laplacian 
problem for $p>1$,
where, here and throughout the paper, $\frac{1}{r}+ \frac{1}{p}=1$:

\noindent
{\bf (P$^{h,\tau}_{A,p}$)}
For $n = 1 \rightarrow N$, find $W^n_{A,r} \in U^h_0$ such that
\begin{align}
&\left(\frac{W^n_{A,r}-W^{n-1}_{A,r}}{\tau_n},\eta^h\right)^h  
+ \left( [M_\epsilon^h(P^h W^n_{A,r})]^{-(p-1)}\, |\del W^n_{A,r}|^{p-2} \, 
\del W^n_{A,r}, \del \eta^h \right)
\nonumber 
\\ &\hspace{2in}= (f^n ,\eta^h)
\qquad \forall \ \eta^h \in U^h_0\,,\label{Pthr}
\end{align}
where $W^0_{A,r} = \pi^h w_0^\epsilon$.

\begin{thm}\label{Qthrstab}
Let the Assumptions (A1) and (A2) hold.
Then for all $r \in (1,2)$,
for all regular partitionings ${\cal T}^h$ of $\Omega$,
and for all $\tau_n >0$,
there exists a solution, 
$W^n_{A,r} \in U^h_0$ and $\vQ^n_{A,r} \in \vS^h$
to the $n^{\rm th}$ step
of (Q$^{h,\tau}_{A,r}$).
In addition, we have that
\eqlabstart
\begin{align}
\max_{n=0 \rightarrow N}
|W^n_{A,r}|_{0,\Omega} + \sum_{n=1}^N |W^n_{A,r}-W^{n-1}_{A,r}|_{0,\Omega}^2 +
\sum_{n=1}^N \tau_n\,|\vQ^n_{A,r}|_{0,r,\Omega}^r
&\leq C\,,
\label{energy3}\\
\left( \sum_{n=1}^N \tau_n |\del W^n_{A,r}|_{0,p,\Omega}^p\right)^{\frac{1}{p}} 
&\leq C 
\,
\label{energy2}
\end{align}
\eqlabend
where $\frac{1}{r}+ \frac{1}{p}=1$.
Moreover, 
(Q$^{h,\tau}_{A,r}$), (\ref{Qthr1},b), is equivalent to (P$^{h,\tau}_{A,p}$), (\ref{Pthr}). 
\end{thm}
\proof
It follows immediately from (\ref{Qthr2}) that 
\begin{align}
\del W^n_{A,r} &=- M^h_\epsilon(P^h W^n_{A,r})\,|\vQ^n_{A,r}|^{r-2} \vQ^n_{A,r} \nonumber \\ 
\Leftrightarrow \quad \vQ^n_{A,r} &= -
[M_\epsilon^h(P^h W^n_{A,r})]^{-(p-1)}\, |\del W^n_{A,r}|^{p-2} \, \del W^n_{A,r}
\quad \mbox{ on } \sigma,
\quad \forall \ \sigma \in {\cal T}^h\,.  
\label{QthrW}
\end{align}
Substituting this expression for $Q^n_{A,r}$ into (\ref{Qthr1}) yields (\ref{Pthr}).
Hence (P$^{h,\tau}_{A,p}$), with (\ref{QthrW}), is equivalent to (Q$^{h,\tau}_{A,r}$). 

We now apply the Brouwer fixed point theorem to prove existence of a solution to  (P$^{h,\tau}_{A,p}$),
and therefore to (Q$^{h,\tau}_{A,r}$). 
Let $F^h : U^h_0 \rightarrow U^h_0$ be such that for any $\varphi^h \in 
U^h_0$, $F^h \varphi^h \in U^h_0$ solves
\begin{align}
&\left(\frac{F^h \varphi^h -W^{n-1}_{A,r}}{\tau_n},\eta^h\right)^h  
+ \left( [M_\epsilon^h(P^h\varphi^h)]^{-(p-1)}\, |\del F^h \varphi^h |^{p-2} \, 
\del  F^h \varphi^h, \del \eta^h \right)
\nonumber 
\\ &\hspace{2in}= (f^n ,\eta^h)
\qquad \forall \ \eta^h \in U^h_0\,.\label{Fh}
\end{align}
The well-posedness of the mapping $F^h$ follows from noting that
(\ref{Fh}) is the Euler--Lagrange system associated with the strictly convex minimization
problem:
\eqlabstart
\begin{align}
\min_{\eta^h \in U^h_0} E^{h,n}_p(\eta^h)\,,
\label{pmin}
\end{align}
where $E^{h,n}_p : U^h_0 \rightarrow {\mathbb R}$ is defined by 
\begin{align} E^{h,n}_p(\eta^h):=
\frac{1}{2 \tau_n}\, |\eta^h-W^{n-1}_{A,r}|_h^2 + \frac{1}{p} \int_{\Omega} [M_\epsilon^h(P^h \varphi^h)]^{-(p-1)}\,
|\del \eta^h|^p \,{\rm d}\vx -(f^n,\eta^h)\,;
\label{pminE}
\end{align}
\eqlabend
that is, there exists a unique element $(F^h \varphi^h) \in U^h_0$ solving (\ref{Fh}). 
It follows immediately from (\ref{pmin},b) that
\begin{align}
\frac{1}{2 \tau_n}\, |F^h \varphi^h-W^{n-1}_{A,r}|_h^2  -(f^n,F^h \varphi^h) \leq 
E^{h,n}_p(F^h \varphi^h) \leq E^{h,n}_p(0) = \frac{1}{2 \tau_n}\, |W^{n-1}_{A,r}|_h^2\,.
\label{pminEb}
\end{align}
It is easily deduced from (\ref{pminEb}) and (\ref{eqnorm}) that
\begin{align} 
F^h \varphi^h \in B_\gamma := \{\eta^h \in U^h_0 : |\eta^h|_{0,\Omega} \leq \gamma \}\,,
\label{Bgamma}
\end{align}
where $\gamma \in {\mathbb R}_{>0}$ 
depends on $|W^{n-1}_{A,r}|_{0,\Omega}$, $|f^n|_{0,\Omega}$ and 
$\tau_n$. 
Hence $F^h : B_\gamma \rightarrow B_\gamma$. In addition, it is easily verified 
that the mapping $F^h$ is continuous, as $M_\eps^h: S^h \rightarrow S^h$ is continuous.
Therefore, the Brouwer fixed point theorem yields that the mapping $F^h$ has at least one 
fixed point in $B_\gamma$. Hence,
there exists a solution to  (P$^{h,\tau}_{A,p}$), (\ref{Pthr}),
and therefore to (Q$^{h,\tau}_{A,r}$), (\ref{Qthr1},b).    

It follows from (\ref{QthrW}) and (\ref{k1epshdef}) that for $n=1 \rightarrow N$
\begin{align}
|\del W^n_{A,r}|_{0,p,\Omega}^p &= 
|[M^h_\epsilon(P^h W^n_{A,r})]^{p-1}\vQ^n_{A,r}|_{0,r,\Omega}^r 
\nonumber \\ &
\leq (k^{\epsilon,h}_{1,\infty})^{p-1}\,
(M^h_\epsilon(P^h W^n_{A,r}),|\vQ^n_{A,r}|^r) 
\,;
\label{WnQn}
\end{align}
where, on noting (\ref{k1epshdef}), (\ref{piconv}) and Assumption (A1), 
\begin{align}
k^{\epsilon,h}_{1,\infty} := \max_{\vx \in \Omega} k^{\epsilon,h}_1(\vx) \leq C\,.
\label{kepsh1max}
\end{align}
Choosing $\eta^h = W^n_{A,r}$, $\vv^h= \vQ^{n}_{A,r}$ in (\ref{Qthr1},b), 
combining and noting 
the simple 
identity
\begin{align}
(a-b)\,a = \frac{1}{2} \left[ a^2 + (a-b)^2 - b^2 \right]
\qquad \forall \ a,\,b \in {\mathbb R}
\,,
\label{simpid}
\end{align} 
we obtain for $n=1 \rightarrow N$, on applying a Young's inequality and (\ref{Poin}), that
for all $\delta >0$
\begin{align}
&|W^n_{A,r}|_h^2 + |W^n_{A,r}-W^{n-1}_{A,r}|_h^2 + 2\tau_n\, (M_{\epsilon}^h(P^h W^n_{A,r}),|\vQ^n_{A,r}|^r)
\nonumber \\
& \hspace{0.3in} = |W^{n-1}_{A,r}|_h^2 + 2 \tau_n\,(f^n,W^n_{A,r})
\nonumber \\
& \hspace{0.3in} \leq |W^{n-1}_{A,r}|_h^2 + 2 \tau_n\,\left[\,\frac{1}{r}\,\delta^{-r}\,|f^n|_{0,r,\Omega}^r
+ \frac{1}{p}\,\delta^p |W^n_{A,r}|_{0,p,\Omega}^p\,\right]
\nonumber \\
& \hspace{0.3in} \leq |W^{n-1}_{A,r}|_h^2 + 2 \tau_n\,\left[\,\frac{1}{r}\,\delta^{-r}\,|f^n|_{0,r,\Omega}^r
+ \frac{1}{p}\,[\delta\,C_\star(\Omega)]^p |\del W^n_{A,r}|_{0,p,\Omega}^p\,\right]\,.
\label{Wnr1}
\end{align}
It follows on summing (\ref{Wnr1}) from $n=1$ to $m$, with $\delta = 1/ 
(C_\star(\Omega)\,[k_{1,\infty}^{\epsilon,h}]^{\frac{1}{r}})$,
and noting (\ref{WnQn}) and (\ref{kepsh1max}) that for $m=1 \rightarrow N$
\begin{align}
&|W^m_{A,r}|_h^2 + \sum_{n=1}^m |W^n_{A,r}-W^{n-1}_{A,r}|_h^2 + \sum_{n=1}^m
\tau_n\, (M_{\epsilon}^h(P^h W^n_{A,r}),|\vQ^n_{A,r}|^r)
\nonumber \\
& \hspace{1in} \leq |W^0_{A,r}|_h^2
+ 2\,[C_\star(\Omega)]^r\, k^{\epsilon,h}_{1,\infty}\,
\sum_{n=1}^m \tau_n\,|f^n|_{0,r,\Omega}^r\,.
\label{Wnr2}
\end{align}  
The desired results (\ref{energy3},b) follow immediately from (\ref{Wnr2}), 
(\ref{eqnorm}), (\ref{fnsum}), (\ref{k1epshdef}), (\ref{kepsh1max})  and   (\ref{WnQn}).  
\endproof

\begin{thm}\label{Qthstab}
Let the Assumptions (A1) and (A2) hold.
Then 
for all regular partitionings ${\cal T}^h$ of $\Omega$,
and for all $\tau_n >0$,
there exists a solution, 
$W_A^n \in U^h_0$ and $\vQ_A^n \in \vS^h$
to the $n^{\rm th}$ step
of (Q$^{h,\tau}_A$).
In addition, we have that
\eqlabstart
\begin{align}
\max_{n=0 \rightarrow N}
|W^n_A|_{0,\Omega} + \sum_{n=1}^N |W_A^n-W_A^{n-1}|_{0,\Omega}^2 +
\sum_{n=1}^N \tau_n\,|\vQ_A^n|_{0,1,\Omega}
&\leq C\,,
\label{energy3a}\\
\max_{n=0 \rightarrow N} \|W_A^n\|_{1,\infty,\Omega} &\leq C\,.
\label{energy2a}
\end{align}
\eqlabend
Moreover, 
(Q$^{h,\tau}_A$), (\ref{Qth1},b), is equivalent to 
(P$^{h,\tau}_A$), (\ref{Pth}),
and  
(M$^{h,\tau}_A$), \linebreak (\ref{Mth}).
Furthermore, for $n=1 \rightarrow N$, having obtained $W_A^n$, then 
$\vQ_A^n = -\lambda_A^n \, \del W_A^n$, where $\lambda_A^n \in S^h_{\geq 0}$
is the Lagrange multiplier associated with the gradient inequality 
constraint in (P$^{h,\tau}_A$). 
\end{thm}
\proof
It follows immediately from (\ref{energy3}), 
on noting that $|\cdot|_{0,1,\Omega} \leq |\Omega| + |\cdot|_{0,r,\Omega}^r$,
that for fixed ${\cal T}^h$ 
and $\{\tau_n\}_{n=1}^N$, 
there exists for $n=1 \rightarrow N$
a subsequence of $\{W^n_{A,r},\vQ^n_{A,r}\}_{r>1}$ (not indicated) and 
 $W^n_A \in U^h_0$ and $\vQ_A^n \in \vS^h$ such that
\begin{align}
W^n_{A,r} \rightarrow W_A^n, \qquad \vQ^n_{A,r} \rightarrow \vQ_A^n \qquad 
\mbox{as } r \rightarrow 1\,. 
\label{rlim1}
\end{align}
We now need to establish that $\{W^n_A,\vQ^n_A\}_{n=1}^N$ solves (Q$^{h,\tau}_A$), (\ref{Qth1},b).
Noting (\ref{rlim1}), one can pass to the limit $r \rightarrow 1$ in (\ref{Qthr1}) to obtain
(\ref{Qth1}). Choosing $\vv^h = \vQ^n_{A,r} - \vpsi^h$ in (\ref{Qthr2}) and noting (\ref{aederiv}),
(\ref{k1epshdef})
and (\ref{kepsh1max}),
one obtains that
\begin{align}
(\del W^n_{A,r},\vpsi^h-\vQ^n_{A,r}) &= (M_\epsilon^h(P^h W^n_{A,r})\,|\vQ^n_{A,r}|^{r-2} \vQ^n_{A,r},
\vQ^n_{A,r}-\vpsi^h)  
\nonumber \\
&\geq \frac{1}{r}\,(M_\epsilon^h(P^h W^n_{A,r}),|\vQ^n_{A,r}|^{r} -|\vpsi^h|^r) 
\nonumber \\
&\geq (M_\epsilon^h(P^h W^n_{A,r}),|\vQ^n_{A,r}| -\frac{1}{r}\,|\vpsi^h|^r)
+ \frac{1-r}{r} \, k_{1,\infty}^{\epsilon,h}\,|\Omega|
\nonumber \\
& \hspace{2in}\qquad \forall \ \vpsi^h \in \vS^h\,.
\label{rlim2}
\end{align}
Noting (\ref{rlim1}), (\ref{Meps_hmm}) and  
(\ref{minmax}), one can pass to the limit $r \rightarrow 1$ in (\ref{rlim2}) to obtain
(\ref{Qth2}). Hence there exists a solution to (Q$^{h,\tau}_A$), (\ref{Qth1},b).

In addition, one can pass to limit $r \rightarrow 1$ in (\ref{energy3}), on noting (\ref{rlim1}),
to obtain the desired result (\ref{energy3a}).

Choosing $\vv^h = \vzero$ and $2\,\vQ^n_A$ in (\ref{Qth2}),
yields 
for $n=1 \rightarrow N$ that
\eqlabstart
\begin{alignat}{2}
&(M_\epsilon^h(P^h W^n_A), |\vQ^n_A|) + (\del W^n_A,\vQ^n_A)  &&= 0 \label{Wneq}\\
\mbox{and hence that} \quad 
&(M_\epsilon^h(P^h W^n_A), |\vv^h|) + (\del W^n_A,\vv^h)  &&\geq 0 
\quad 
\forall \ \vv^h\in \vS^h\,.
\label{Wnineq}
\end{alignat}
\eqlabend 
Choosing
\begin{align}
\vv^h = \left\{\begin{array}{cl}
- \del W^n_A \mid_{\sigma_\star} & \qquad \sigma=\sigma_\star\,,\\
\vzero &  \qquad \sigma \neq \sigma_\star\,;
\end{array} \right.
\label{vvhch}
\end{align}
in (\ref{Wnineq}), and repeating for all $\sigma_\star \in {\cal T}^h$, yields 
for $n=1 \rightarrow N$ that
\begin{align}
|\del W^n_A| \leq M^h_{\varepsilon}(P^h W^n_A) \qquad \mbox{a.e.\ on } \Omega\,.  
\label{gradcon}
\end{align}
As $W^n_A \in U^h_0$, it follows from (\ref{gradcon}), (\ref{k1epshdef}), (\ref{kepsh1max}), 
(\ref{Poin}), (\ref{piconv})
and our choice of $W^0_A$
that the desired result (\ref{energy2a}) holds.

It follows from (\ref{gradcon}) that $W^n_A \in K^h(W^n_A)$.
Choosing $\eta^h = \varphi^h -W^n_A$ for any $\varphi^h \in K^h(W^n_A)$ in (\ref{Qth1}), 
we obtain, on noting (\ref{Wneq}) and (\ref{Kh}), that  
\begin{align}
&\left(\frac{W^n_A-W^{n-1}_A}{\tau_n},\varphi^h-W^n_A\right)^h  - (f^n 
,\varphi^h-W^n_A) \nonumber \\
& \hspace{0.3in}= (\vQ^n_A,\del (\varphi^h-W^n_A))
= (M^h_\epsilon(P^h W^n_A),|\vQ^n_A|)+(\del \varphi^h,\vQ^n_A)\geq 0\,.  
\label{Qth1ineq}
\end{align} 
Hence $\{W^n_A\}_{n=1}^N$ solves (P$^{h,\tau}_A$), (\ref{Pth}). 
It follows from (\ref{Qth1}) that $\vQ^n_A \in \vZ^{h,n}$, $n=1 \rightarrow N$.
Therefore 
(\ref{Qth2}) immediately yields (\ref{Mth})
on choosing $\vv^h \in \vZ^{h,n}$.
Hence $\{\vQ^n_A\}_{n=1}^N$ solves (M$^{h,\tau}_A$), (\ref{Mth}). 
Therefore a solution $\{W^n_A,\vQ^n_A\}_{n=1}^N$ of (Q$^{h,\tau}_A$) 
solves (P$^{h,\tau}_A$) and (M$^{h,\tau}_A$). 

We now prove the reverse.
If $\{W^n_A\}_{n=1}^N$ solves (P$^{h,\tau}_A$), then, for $n=1 \rightarrow N$,
$W^n_A$ is the unique solution to the strictly convex minimization
problem:
\eqlabstart
\begin{align}
\min_{\eta^h \in K^h(W^n_A)} E^{h,n}(\eta^h)\,,
\label{Kmin}
\end{align}
where $E^{h,n} : U^h_0 \rightarrow {\mathbb R}$ is defined by
\begin{align} E^{h,n}(\eta^h):=
\frac{1}{2 \tau_n}\, |\eta^h-W_A^{n-1}|_h^2 -(f^n,\eta^h)\,.
\label{KminE}
\end{align}
\eqlabend 
Next we introduce the Lagrangian $L^{h,n}: U^h_0 \times S^h_{\geq 0} \rightarrow {\mathbb R}$ 
defined by  
\begin{align}
L^{h,n}(\eta^h,\mu^h):=  E^{h,n}(\eta^h) + \tfrac{1}{2}\,(\mu^h,|\del \eta^h|^2 - 
[M^h_\epsilon(P^h W^n_A)]^2)\,.
\label{Lag}  
\end{align}
As $k_0 >0$, we note 
that the Slater constraint qualification hypothesis,
see e.g.\ (5.34) on page 69 in Ekeland and Temam \cite{ET}, is obviously satisfied; that is,
there exists an $\eta^h_0 \in U^h_0$ such that $|\del \eta^h_0| 
< M^h_\epsilon(P^h W^n_A)$\,.    
Hence it follows from the Kuhn--Tucker theorem, see e.g.\ Theorem 5.2 on page 69 
in Ekeland and Temam \cite{ET}, that there exists a $\lambda^n_A \in S^h_{\geq 0}$ such that
\begin{align}
L^{h,n}(W^n_A,\mu^h) \leq L^{h,n}(W^n_A,\lambda^n_A)
\leq L^{h,n}(\eta^h,\lambda^n_A) \qquad \forall \ \eta^h \in U^h_0, \quad \forall \ \mu^h
\in S^h_{\geq0}\,.
\label{Lagineq}
\end{align}
The first inequality in (\ref{Lagineq}) yields for $\mu^h =0$ and $2\lambda^n_A$ that
\begin{align}
&(\lambda^n_A,|\del W^n_A|^2-[M^h_\epsilon(P^h W^n_A)]^2)=0  
\nonumber \\
\Rightarrow \qquad  
&(\lambda^n_A\,|\del W^n_A|,|\del W^n_A|-M^h_\epsilon(P^h W_A^n))=0\,,
\label{Lagineq1}
\end{align}
as $W_A^n \in K^h(W_A^n)$.
The second inequality in (\ref{Lagineq}) yields that
\begin{align}
\left(\frac{W_A^n-W_A^{n-1}}{\tau_n},\eta^h\right)^h + (\lambda_A^n\,\del W_A^n,\del \eta^h) 
&= (f^n 
,\eta^h)
\qquad \forall \ \eta^h \in U^h_0\,.
\label{Qth1near}
\end{align}
It follows that (\ref{Qth1}) holds on 
setting $\vQ_A^n = -\lambda_A^n\,\del W_A^n$, and $\vQ_A^n \in \vZ^{h,n}$.
Furthermore, we have from this definition for $\vQ_A^n \in \vZ^{h,n}$ and (\ref{Lagineq1}) that
for all $\vv^h \in \vZ^{h,n}$
\begin{align}
(M^h_\epsilon(P^h W^n_A),|\vQ^n_A|) = -( \vQ^n_A,\del W^n_A) =-(\vv^h,\del W^n_A)
\leq (M^h_\epsilon(P^h W^n_A),|\vv^h|)\,,
\label{Qthmin}
\end{align}
where we have recalled that $W_A^n \in K^h(W_A^n)$ for the last inequality.
Hence, for $n=1 \rightarrow N$, $\vQ^n_A = - \lambda^n_A\,\del W_A^n \in \vZ^{h,n}$
solves
the minimization problem (M$^{h,\tau}_A$), (\ref{Mth}). 
Since the inequality in (\ref{Qthmin}) holds for all $\vv^h \in \vS^h$, 
it follows from this and 
the first equality in (\ref{Qthmin}) that (\ref{Qth2}) holds.  
Therefore a solution $\{W_A^n,\vQ_A^n\}_{n=1}^N$ of (P$^{h,\tau}_A$) and (M$^{h,\tau}_A$)
solves (Q$^{h,\tau}_A$). 
\endproof

\subsection{Approximation (Q$^{h,\tau}_{B}$)}

Our second fully practical finite element approximation is:

\noindent
{\bf (Q$^{h,\tau}_B$)}
For $n = 1 \rightarrow N$, find $W_B^n \in S^h$ and 
$\vQ_B^n \in \vV^h$ such that
\eqlabstart
\begin{alignat}{2}
\left(\frac{W_B^n-W_B^{n-1}}{\tau_n},\eta^h\right) + (\del\,.\,\vQ_B^n,\eta^h) &= (f^n 
,\eta^h)
\quad \;&&\forall \ \eta^h \in S^h,\label{Qth1B}\\
(M_\epsilon^h (W_B^n), |\vv^h|-|\vQ_B^n|)^h - (W_B^n,\del\,.\,(\vv^h-\vQ_B^n))  &\geq0
\quad \;
&&\forall \ \vv^h\in \vV^h\,;
\label{Qth2B}
\end{alignat}
\eqlabend
where $W_B^0 = P^h[\pi^h w_0^\epsilon]$.

For computational and theoretical purposes, it is convenient 
to consider the following regularization 
of (Q$^{h,\tau}_B$)
for a given $r>1$:

\noindent
{\bf (Q$^{h,\tau}_{B,r}$)}
For $n = 1 \rightarrow N$, find $W^n_{B,r} \in S^h$ and
$\vQ^n_{B,r} \in \vV^h$ such that
\eqlabstart
\begin{alignat}{2}
\left(\frac{W^n_{B,r}-W^{n-1}_{B,r}}{\tau_n},\eta^h\right) 
+ (\del\,.\,\vQ^n_{B,r},\eta^h) &= (f^n 
,\eta^h)
\quad \;&&\forall \ \eta^h \in S^h,\label{Qthr1B}\\
(M_\epsilon^h(W^n_{B,r})\,|\vQ^n_{B,r}|^{r-2} \vQ^n_{B,r},
\vv^h)^h - (W^n_{B,r},\del\,.\,\vv^h)  &=0
\quad \;
&&\forall \ \vv^h\in \vV^h\,;
\label{Qthr2B}
\end{alignat}
\eqlabend
where $W^0_{B,r} = P^h[\pi^h w_0^\epsilon]$.

\begin{thm}\label{QthrstabB}
Let the Assumptions (A1) and (A2) hold.
Then for all $r \in (1,2)$,
for all regular partitionings ${\cal T}^h$ of $\Omega$,
and for all $\tau_n >0$,
there exists a solution, 
$W^n_{B,r} \in S^h$ and $\vQ^n_{B,r} \in \vV^h$
to the $n^{\rm th}$ step
of (Q$^{h,\tau}_{B,r}$).
In addition, we have for $\tau \in (0,\frac{1}{2}]$ that
\eqlabstart
\begin{align}
\max_{n=0 \rightarrow N}
|W^n_{B,r}|_{0,\Omega} + \sum_{n=1}^N |W^n_{B,r}-W^{n-1}_{B,r}|_{0,\Omega}^2 +
\sum_{n=1}^N \tau_n\,|\vQ^n_{B,r}|_{0,r,\Omega}^r
&\leq C\,,
\label{energy3B} \\
\sum_{n=1}^N \tau_n^2\,|\del\,.\,\vQ^n_{B,r}|_{0,\Omega}^2
&\leq C\,.
\label{energy2B}
\end{align}
\eqlabend
\end{thm}
\proof
It follows from (\ref{Qthr1B}) and (\ref{Ph}) that
\begin{align}
W^n_{B,r}=g^n-\tau_n \,\del\,.\,\vQ^n_{B,r}\,, \quad \mbox{where}\quad
g^n = W^{n-1}_{B,r} + \tau_n\,P^h f^n\,.
\label{gn}
\end{align}
Substituting (\ref{gn}) into (\ref{Qthr2B}) yields that the $n^{\rm th}$ step of (Q$^{h,\tau}_{B,r}$)
can be rewritten as find $\vQ^{n}_{B,r} \in \vV^h$ such that
\begin{align}
&(M_\epsilon^h(g^n-\tau_n\,\del\,.\,\vQ^n_{B,r})\,|\vQ^n_{B,r}|^{r-2} \vQ^n_{B,r},
\vv^h)^h + \tau_n\, (\del\,.\,\vQ^n_{B,r},\del\,.\,\vv^h)  
\nonumber \\
&\hspace{2.5in}=(g^n,\del\,.\,\vv^h) 
\qquad
\forall \ \vv^h\in \vV^h\,.
\label{Qthr2Balt}
\end{align}
One can apply the Brouwer fixed point theorem to prove existence of a solution to
(\ref{Qthr2Balt}), and therefore to (Q$^{h,\tau}_{B,r}$).
Let $\vG^h : \vV^h \rightarrow \vV^h$ be such that 
for any $\vpsi^h \in \vV^h$, $\vG^h \vpsi^h \in \vV^h$ solves
\begin{align}
&(M_\epsilon^h(g^n-\tau_n\,\del\,.\,\vpsi^h)\,|\vG^h \vpsi^h|^{r-2} \vG^h \vpsi^h,
\vv^h)^h + \tau_n\, (\del\,.\,(\vG^h \vpsi^h),\del\,.\,\vv^h)  
\nonumber \\
&\hspace{2.5in}=(g^n,\del\,.\,\vv^h) 
\qquad
\forall \ \vv^h\in \vV^h\,.
\label{Qthr2BaltBFT}
\end{align}
The well-posedness of the mapping $\vG^h$ follows from noting that 
(\ref{Qthr2BaltBFT}) is the Euler-Lagrange system associated with the strictly 
convex minimization problem:
\eqlabstart
\begin{align}
\min_{\vv^h \in \vV^h} J^{h,n}_r(\vv^h)\,,
\label{rmin}
\end{align}
where $J^{h,n}_r : \vV^h \rightarrow \mathbb{R}$ is defined by  
\begin{align}
J^{h,n}_r(\vv^h) :=
\frac{1}{r}\,(M_\epsilon^h(g^n-\tau_n\,\del\,.\,\vpsi^h),|\vv^h|^{r})^h + \frac{\tau_n}{2}\,
|\del\,.\,\vv^h|_{0,\Omega}^2-(g^n,\del\,.\,\vv^h);
\label{rminJ}
\end{align} 
\eqlabend
that is, there exists a unique element $(\vG^h \vpsi^h) 
\in \vV^h$ solving (\ref{Qthr2BaltBFT}). 
It follows immediately from (\ref{rmin},b) that $J^{h,n}_r(\vG^h \vpsi^h) \leq J^{h,n}_r(\vzero)$,
and this yields, on noting (\ref{k1epshdef}) 
and (\ref{convex}), that
\begin{align}
\frac{k_0}{r}\,(\,|\vG^h \vpsi^h|^{r},1) + \frac{\tau_n}{2}\,
|\del\,.\,(\vG^h \vpsi^h)|_{0,\Omega}^2 &\leq (g^n, \del\,.\,(\vG^h \vpsi^h)\,)
\,. 
\label{rminJb}
\end{align} 
It is easily from (\ref{rminJb}) that
\begin{align} 
\vG^h \vpsi^h \in \vB_\gamma := \{\vv^h \in \vV^h : |\vv^h|_{0,r,\Omega} \leq \gamma \}\,,
\label{BgammaB}
\end{align}
where $\gamma \in {\mathbb R}_{>0}$ 
depends on $|g^n|_{0,\Omega}$, $r$ and 
$\tau_n$. 
Hence $\vG^h : \vB_\gamma \rightarrow \vB_\gamma$. In addition, it is easily verified 
that the mapping $\vG^h$ is continuous, as $M_\eps^h: S^h \rightarrow S^h$ is continuous.
Therefore, the Brouwer fixed point theorem yields that the mapping $\vG^h$ has at least one 
fixed point in $\vB_\gamma$. Hence,
there exists a solution to  (Q$^{h,\tau}_{B,r}$), (\ref{Qthr1B},b). 

Choosing $\eta^h = W^n_{B,r}$, $\vv^h = \vQ^n_{B,r}$ in (\ref{Qthr1B},b), 
combining and noting (\ref{simpid}) yields, similarly to (\ref{Wnr1}), that
\begin{align}
&|W^n_{B,r}|_{0,\Omega}^2 + |W^n_{B,r}-W^{n-1}_{B,r}|_{0,\Omega}^2 + 
2\tau_n\, (M_{\epsilon}^h(W^n_{B,r}),|\vQ^n_{B,r}|^r)^h
\nonumber \\
& \hspace{0.1in} = |W^{n-1}_{B,r}|_{0,\Omega}^2 + 2 \tau_n\,(f^n,W^n_{B,r})
\leq |W^{n-1}_{B,r}|_{0,\Omega}^2 + \tau_n\,\left[
|W^n_{B,r}|_{0,\Omega}^2 +
|f^n|_{0,\Omega}^2 \right] 
\,.
\label{Wnr1B}
\end{align}
It follows from (\ref{Wnr1B}),
on noting that $(1-\tau_n)^{-1} \leq (1+ 2\tau_n) \leq e^{2 \tau_n}$ 
as $\tau_n \in (0,\frac{1}{2}]$ and (\ref{fnsum}),
that 
for $n=1 \rightarrow N$
\begin{align}
|W^n_{B,r}|_{0,\Omega}^2 
&\leq 
e^{2\tau_n}\left[\, |W^{n-1}_{B,r}|_{0,\Omega}^2 + \tau_n\,|f^n|_{0,\Omega}^2 \right] 
\nonumber \\
&\leq e^{2t_n} \left[ 
|W^{0}_{B,r}|_{0,\Omega}^2 + \sum_{m=1}^N \tau_m\,|f^m|_{0,\Omega}^2\right] 
\leq C
\,,
\label{Wnr1Ba}
\end{align}
which yields the first bound in (\ref{energy3B}).
Summing (\ref{Wnr1B}) from $n=1 \rightarrow N$ yields, on noting
(\ref{k1epshdef}), (\ref{convex}) and (\ref{Wnr1Ba}), the second and third bounds in 
(\ref{energy3B}). 

Choosing $\eta^h = \del\,.\,Q^n_{B,r}$ in (\ref{Qthr1B}) yields that
\begin{align}
\tau_n^2\, |\del\,.\vQ^n_{B,r}|^2_{0,\Omega} 
&=
\tau_n \,(W^{n-1}_{B,r}-W^{n}_{B,r} + \tau_n\, f^n, \del\,.\, \vQ^n_{B,r})
\nonumber \\
&\leq 2\left[\,|W^n_{B,r}-W^{n-1}_{B,r}|_{0,\Omega}^2 + \tau_n^2\,|f^n|_{0,\Omega}^2\right]\,.
\label{Wnr1Bb}
\end{align}
Summing (\ref{Wnr1Bb}) from $n =1 \rightarrow N$, and noting (\ref{energy3B}) and (\ref{fnsum}),
yields the desired result (\ref{energy2B}).
\endproof

\begin{thm}\label{QthstabB}
Let the Assumptions (A1) and (A2) hold.
Then 
for all regular partitionings ${\cal T}^h$ of $\Omega$,
and for all $\tau_n >0$,
there exists a solution, 
$W_B^n \in S^h$ and $\vQ_B^n \in \vV^h$
to the $n^{\rm th}$ step
of (Q$^{h,\tau}_B$).
In addition, we have for $\tau \in (0,\frac{1}{2}]$ that
\eqlabstart
\begin{align}
\max_{n=0 \rightarrow N}
|W^n_B|_{0,\Omega} + \sum_{n=1}^N |W_B^n-W_B^{n-1}|_{0,\Omega}^2 +
\sum_{n=1}^N \tau_n\,|\vQ_B^n|_{0,1,\Omega}
&\leq C\,.
\label{energy3aB}
\\
\sum_{n=1}^N \tau_n^2\,|\del\,.\,\vQ_B^n|_{0,\Omega}^2
&\leq C\,.
\label{energy2aB}
\end{align}
\eqlabend
\end{thm}
\proof
Similarly to (\ref{rlim1}), 
on noting that $|\cdot|_{0,1,\Omega} \leq |\Omega| + |\cdot|_{0,r,\Omega}^r$,
it follows from (\ref{energy3B},b), that for fixed ${\cal T}^h$ and 
$\{\tau_n\}_{n=1}^N$, there exists a subsequence of $\{W^{n}_{B,r},Q^{n}_{B,r}\}_{r>1}$
(not indicated)
and $W^n_B \in S^h$ and $\vQ^n_{B} \in \vV^h$ such that
\begin{align}
W^n_{B,r} \rightarrow W^n_B, \qquad \vQ^n_{B,r} \rightarrow \vQ^n_B \qquad \mbox{as } 
r \rightarrow 1,
\label{rlim1B}
\end{align} 
and the bounds (\ref{energy3aB},b) hold.
One can now immediately pass to the limit $r \rightarrow 1$ in 
(\ref{Qthr1B}) to obtain (\ref{Qth1B}).
Similarly to (\ref{rlim2}), choosing $\vv^h = \vQ^n_{B,r}- \vpsi^h$ in (\ref{Qthr2B})
and noting (\ref{aederiv}), (\ref{k1epshdef})
and (\ref{kepsh1max}), one obtains that
\begin{align}
(W^n_{B,r}, \del\,.\,(\vQ^{n}_{B,r}-\vpsi^h)\,)
&\geq (M_\epsilon^h(W^n_{B,r}),|\vQ^n_{B,r}| -\frac{1}{r}\,|\vpsi^h|^r)^h
+ \frac{1-r}{r} \, k_{1,\infty}^{\epsilon,h}\,|\Omega|
\nonumber \\
& \hspace{2.2in}
\forall \ \vpsi^h \in \vV^h\,.
\label{rlim2B}
\end{align}
Noting (\ref{rlim1B}), one can pass to the limit $r \rightarrow 1$ in (\ref{rlim2B})
to obtain (\ref{Qth2B}). Hence there exists a solution to (Q$^{h,\tau}_B$), (\ref{Qth1B},b). 
\endproof

\section{Convergence}
\label{conv}
\setcounter{equation}{0}

We introduce the following discrete approximation of the mixed formulation:
\vspace{2mm}

\noindent
{\bf (Q$^{\tau}$)}
For $n = 1 \rightarrow N$, find $w^n \in W^{1,\infty}_0(\Omega)$ and 
$\vq^n \in \vV^{\cal M}(\Omega)$ such that
\eqlabstart
\begin{alignat}{2}
\left(\frac{w^n-w^{n-1}}{\tau_n},\eta\right) + (\del\,.\,\vq^n,\eta) &= (f^n 
,\eta)
\quad \;&&\forall \ \eta \in L^2(\Omega),\label{Qt1}\\
\langle |\vv|-|\vq^n|,M_\epsilon(w^n)\rangle_{C(\overline{\Omega})} - (\del\,.\,(\vv-\vq^n),
w^n)  &\geq0
\quad \;
&&\forall \ \vv\in \vV^{\cal M}(\Omega)\,;
\label{Qt2}
\end{alignat}
\eqlabend
where $w^0 = w_0^\epsilon$.

For any  $\chi \in W^{1,\infty}_0(\Omega)$, we introduce 
the closed convex non-empty set
\begin{align}
K(\chi) := \{ \eta \in W^{1,\infty}_0(\Omega) : |\del \eta| \leq M_\epsilon(\chi) 
\quad \mbox{a.e.\ on } \Omega\}\,.   
\label{K}
\end{align}
Associated with (Q$^{\tau}$) is the corresponding approximation of the primal 
quasi-variational inequality:

\noindent
{\bf (P$^{\tau}$)}
For $n = 1 \rightarrow N$, find $w^n \in K(w^n)$ such that
\begin{align}
&\left(\frac{w^n-w^{n-1}}{\tau_n},\eta-w^n\right)  
\geq (f^n ,\eta-w^n)
\qquad \forall \ \eta \in K(w^n)\,,\label{Pt}
\end{align}
where $w^0 = w_0^\epsilon$.

\vspace{2mm}

In Section \ref{secconvA} we show, 
for a fixed time partition $\{\tau_n\}_{n=1}^N$, that
a subsequence of $\{ \{W^n_A,\vQ^n_A\}_{n=1}^N \}_{h>0}$, where
$\{W^n_A,\vQ^n_A\}_{n=1}^N$ solves (Q$^{h,\tau}_A$),
converges, as $h \rightarrow 0$ to $\{w^n,\vq^n\}_{n=1}^N$ solving (Q$^{\tau}$).
In Section \ref{secconvB} we show,
for a fixed time partition $\{\tau_n\}_{n=1}^N$, that
a subsequence of $\{ \{W^n_{B,r},\vQ^n_{B,r}\}_{n=1}^N \}_{h>0}$, where
$\{W^n_{B,r},\vQ^n_{B,r}\}_{n=1}^N$ solves (Q$^{h,\tau}_{B,r}$),
converges, as $h \rightarrow 0$ and $r \rightarrow 1$, 
to $\{w^n,\vq^n\}_{n=1}^N$ solving (Q$^{\tau}$).
For our final convergence result in Section \ref{secconvtau}, 
we need an extra assumption on the data.

\vspace{2mm}

\noindent
{\bf (A3)}
$w^\epsilon_0 \geq 0$ and 
$f \in L^\infty(0,T;L^2(\Omega))$.  

\vspace{2mm}

Under this further assumption, we will show that
a subsequence of \linebreak $\{ \{w^n,\vq^n\}_{n=1}^N \}_{\tau >0}$, 
where $\{w^n,\vq^n\}_{n=1}^N$ solves (Q$^\tau$),
converges, as $\tau \rightarrow 0$, to $\{w,\vq\}$ solving

\vspace{2mm}

\noindent
{\bf (Q)}
Find $w \in L^\infty(0,T;W^{1,\infty}_0(\Omega))\cap W^{1,\infty}(0,T;
[C^1_0(\overline{\Omega})]^*)$ and 
$\vq \in L^\infty(0,T;$ $[{\cal M}(\overline{\Omega})]^d)$ such that
\eqlabstart
\begin{align}
&\int_0^T \left[ 
\langle\frac{\partial w}{\partial t},\eta \rangle_{C^1_0(\overline{\Omega})} - 
\langle \vq,\del \eta\rangle_{C(\overline{\Omega})} -(f 
,\eta) \right] \,{\rm d}t =0
\nonumber \\ &\hspace{3.1in}
\forall \ \eta \in L^1(0,T;C^1_0(\overline{\Omega})),\label{Q1}\\
&\int_0^T \left[ \langle |\vv|-|\vq|,M_\epsilon(w)\rangle_{C(\overline{\Omega})} - 
(\del\,.\,\vv-f,
w)\right] \,{\rm d}t  \nonumber \\
& \hspace{0.6in} \geq \tfrac{1}{2} \left[\, |w(\cdot,T)|^2_{0,\Omega} - |w^\epsilon_0(\cdot)|^2_{0,\Omega}
\, \right]
\qquad \forall \ \vv\in L^1(0,T;\vV^{\cal M}(\Omega))\,;
\label{Q2}
\end{align}
\eqlabend
where $w(\cdot,0) = w_0^\epsilon(\cdot)$.

Associated with (Q) is the corresponding approximation of the primal 
quasi-variational inequality:

\noindent
{\bf (P)}
Find $w \in L^\infty(0,T;K(w))\cap W^{1,\infty}(0,T;[C^1_0(\overline{\Omega})]^*)$ such that
\begin{align}
&\int_0^T \left[ 
\langle\frac{\partial w}{\partial t},\eta \rangle_{C^1_0(\overline{\Omega})} - 
(f ,\eta-w) \right] \,{\rm d}t 
\geq 
\tfrac{1}{2} \left[ \,|w(\cdot,T)|^2_{0,\Omega} - |w^\epsilon_0(\cdot)|^2_{0,\Omega}
\,\right]
\nonumber \\ &\hspace{2.7in}
\forall \ \eta \in L^1(0,T;K(w)\cap C^1_0(\overline{\Omega})),\label{P}
\end{align}
where $w(\cdot,0) = w_0^\epsilon(\cdot)$.

\begin{rem}\label{remqvi}
One might expect the inequality in the primal quasi-variational inequality (P) to be such that
\begin{align}
&\int_0^T \left[ 
\langle\frac{\partial w}{\partial t},\eta -w \rangle_{C^1_0(\overline{\Omega})} - 
(f ,\eta-w) \right] \,{\rm d}t 
\geq 0\,.
\label{Porig}
\end{align}
However, the term 
\begin{align}
\int_0^T 
\langle\frac{\partial w}{\partial t},w \rangle_{C^1_0(\overline{\Omega})}  
\,{\rm d}t
\label{wnotdef}
\end{align}
is not well-defined for 
$w \in L^\infty(0,T;W^{1,\infty}_0(\Omega))\cap W^{1,\infty}(0,T;[C^1_0(\overline{\Omega})]^*)$,
and has been rewritten to yield (\ref{P}), which is well defined.
This follows from (\ref{compact1}) with ${\mathcal B} =L^2(\Omega)$,
and, for example, the reflexive Banach spaces ${\mathcal B}_0 =H^1_0(\Omega)$
and ${\mathcal B}_1 =[W^{2,s}_0(\Omega)]^\star$ with $s \in (d,\infty)$;
see the proof of Theorem \ref{conthm} below.
In addition, the test space has been smoothed to make the first term on the 
left-hand side of (\ref{P}) well-defined.

Similar remarks apply to (Q), where one might expect the inequality in 
(\ref{Q2}) to take the form
\begin{align}
&\int_0^T \left[ \langle |\vv|-|\vq|,M_\epsilon(w)\rangle_{C(\overline{\Omega})} - 
(\del\,.\,(\vv-\vq),
w)\right] \,{\rm d}t \geq 0\,.  
\label{Q2orig}
\end{align}
However, the term 
\begin{align}
\int_0^T 
(\del\,.\, \vq, w )
\,{\rm d}t = - \int_0^T 
\langle \vq, \del w \rangle_{C(\overline{\Omega})}
\,{\rm d}t
\label{qnotdef}
\end{align} 
is not well-defined for 
$w \in L^\infty(0,T;W^{1,\infty}_0(\Omega))$
and 
$\vq \in L^\infty(0,T;[{\cal M}(\overline{\Omega})]^d)$.
This term has been rewritten using (\ref{Q1}) formally with $\eta = w$,
and the rewrite of the term (\ref{wnotdef}) employed in (\ref{P}), 
to yield (\ref{Q2}).
\end{rem}

\subsection{Convergence of (Q$^{h,\tau}_A$) to (Q$^{\tau}$)}
\label{secconvA}
\begin{thm}\label{hconthm}
Let the Assumptions (A1) and (A2) hold.
For any fixed time partition $\{\tau_n\}_{n=1}^N$
and for all regular partitionings ${\cal T}^h$ of $\Omega$,
there exists a subsequence of $\{\{W^n_A,\vQ^n_A\}_{n=1}^N\}_{h>0}$ (not indicated),
where $\{W^n_A,\vQ^n_A\}_{n=1}^N$ solves (Q$^{h,\tau}_A$), 
such that as $h \rightarrow 0$
\eqlabstart
\begin{alignat}{3}
\del W^{n}_A&\rightarrow \del w^n  \qquad &&
\mbox{weak$^\star$ in } [L^\infty(\Omega)]^d,
\qquad && n=0 \rightarrow N,
\label{wgradwcon}
\\
W^{n}_A&\rightarrow w^n  \qquad &&
\mbox{strongly in } C(\overline{\Omega}),
\qquad && n=0 \rightarrow N,
\label{wscon}
\\
M_\epsilon^h( P^h W^{n}_A) &\rightarrow M_\epsilon(w^n)  \qquad &&
\mbox{strongly in } L^\infty(\Omega),
\qquad && n=0 \rightarrow N,
\label{Mwscon}
\\
\vQ^n_A &\rightarrow \vq^n \qquad &&
\mbox{vaguely in } [{\cal M}(\overline{\Omega})]^d,
\qquad && n=1 \rightarrow N;
\label{qwcon}
\end{alignat}
\eqlabend
where $\{w^n,\vq^n\}_{n=1}^N$ is a solution of (Q$^\tau$), (\ref{Qt1},b).
\end{thm}
\proof
The desired subsequence convergence results (\ref{wgradwcon},b,d) 
for a fixed time partition $\{\tau_n\}_{n=1}^N$
follow immediately from (\ref{energy3a},b),                      
on noting that $W^{1,\infty}(\Omega)$ is compactly embedded in 
$C(\overline{\Omega})$ 
and (\ref{Mweak1}).
Next we note 
that
\begin{align}
&|M^h_\epsilon(P^h W^n_A)-M_\epsilon(w^n)|_{0,\infty,\Omega}
\nonumber \\
&\hspace{0.4in}\leq
|M_\epsilon^h(P^h W^n_A)-M_\epsilon(P^h W^n_A)|_{0,\infty,\Omega}
+|M_\epsilon(P^h W^n_A)-M_\epsilon(w^n)|_{0,\infty,\Omega}\,.
\label{Qth2cone}
\end{align}
It follows from (\ref{Mepsmm}), (\ref{minmax}) and (\ref{Phstab}) that
\begin{align}
&|M_\epsilon(P^h W^n_A)-M_\epsilon(w^n)|_{0,\infty,\Omega}
\nonumber \\
&\hspace{1.2in} \leq C(\epsilon^{-1})\,|w^n-P^h W^n_A|_{0,\infty,\Omega} 
\nonumber \\
& \hspace{1.2in}
\leq C(\epsilon^{-1})\left[\,|(I-P^h) w^n|_{0,\infty,\Omega}
+ |w^n-W^n_A|_{0,\infty,\Omega}\right]. 
\label{Qth2conh}
\end{align}
Hence, the desired result (\ref{Mwscon}) 
follows from (\ref{Qth2cone}), (\ref{Mepshapprox}), (\ref{Qth2conh}), 
(\ref{Pconv}), (\ref{piconv}) and (\ref{wscon}). 

We now need to establish that $\{w^n,\vq^n\}_{n=1}^N$ solve (Q$^\tau$), (\ref{Qt1},b).
For any $\eta \in C^\infty_0(\Omega)$, we choose  $\eta^h = \pi^h \eta$ in (\ref{Qth1})
and now pass to the limit $h \rightarrow 0$ for the subsequence
to obtain, on noting (\ref{wscon}),
(\ref{nierr}), (\ref{piconv}) and (\ref{qwcon}),  
for $n=1 \rightarrow N$ that
\begin{align}
\left(\frac{w^n-w^{n-1}}{\tau_n},\eta\right) 
- \langle \vq^n,\del \eta \rangle_{C(\overline{\Omega})} 
&= (f^n 
,\eta)
\qquad \forall \ \eta \in C^\infty_0(\Omega)\,.\label{Qt1sm}
\end{align}
It follows from (\ref{Qt1sm}), (\ref{fnsum}) and as $w^n \in C(\overline{\Omega})$ that
\begin{align}
\left| \langle \vq^n,\del \eta \rangle_{C(\overline{\Omega})}
\right| \leq C(\tau_n^{-1})\,|\eta|_{0,\Omega}
\qquad \forall \ 
\eta \in C^\infty_0(\Omega)\,.\label{Qt1sma}
\end{align}
We deduce from (\ref{Qt1sma}) that the distributional divergence of $\vq^n$
belongs to $L^2(\Omega)$, and hence $\vq^n \in \vV^{\cal M}(\Omega)$, $n=1 \rightarrow N$,
and so (\ref{Qt1sm})   
can be rewritten as 
\begin{align}
\left(\frac{w^n-w^{n-1}}{\tau_n},\eta\right) + (\del\,.\,\vq^n,\eta) 
&= (f^n 
,\eta)
\qquad \forall \ \eta \in C^\infty_0(\Omega)\,.\label{Qt1smb}
\end{align}
Noting that $C^\infty_0(\Omega)$ is dense in $L^2(\Omega)$ 
and that $w^n,\,\del\,.\,\vq^n,\,f^n \in L^2(\Omega)$ 
yields the desired (\ref{Qt1}).

For any $\vv \in [C^\infty(\overline{\Omega})]^d$, we choose  $\vv^h = \vP^h \vv$ in (\ref{Qth2})
and now try to pass to the limit for the subsequence as $h \rightarrow 0$. 
First we note from (\ref{wgradwcon}), (\ref{Pconv}) and as $w^n \in W^{1,\infty}_0(\Omega)$ that
for $n=1 \rightarrow N$
\begin{align}
\lim_{h \rightarrow 0} (\del W^n_A, \vP^h\vv) = (\del w^n, \vv) = - (w^n, \del\,.\,\vv) 
\,.
\label{Qth2cona}
\end{align}
It follows  from (\ref{Qth1}) with $\eta^h = W^n_A$,
(\ref{wscon}),
(\ref{nierr}), (\ref{energy2a}) and (\ref{Qt1}) with $\eta = w^n$  
that for $n=1 \rightarrow N$
\begin{align}
\lim_{h \rightarrow 0} (\del W^n_A, \vQ^n_A) &= 
\lim_{h \rightarrow 0} \left[ 
\left(\frac{W^n_A-W^{n-1}_A}{\tau_n},W^n_A\right)^h - (f^n 
,W^n_A)\right]
\nonumber \\
=&\left(\frac{w^n-w^{n-1}}{\tau_n},w^n\right) - (f^n 
,w^n)
 = - (w^n,\del\,.\vq^n)\,. 
\label{Qth2conb}
\end{align}
Next we note 
that
\begin{align}
&(M_\epsilon^h(P^h W^n_A),|\vQ^n_A|-|\vP^h \vv|) 
\nonumber \\
& \hspace{0.2in}= (M_\epsilon(w^n),|\vQ^n_A|-|\vP^h \vv|)
+(M^h_\epsilon(P^h W^n_A)-M_\epsilon(w^n), |\vQ^n_A|-|\vP^h \vv|)\,.
\label{Qth2conc}
\end{align}
As $M_\epsilon(w^n) \in C(\overline{\Omega})$ is positive, 
it follows from (\ref{qwcon}), (\ref{Mweak2}) and  
(\ref{Pconv}) that
\begin{align}
\liminf_{h \rightarrow 0}\, (M_\epsilon(w^n),|\vQ^n_A|-|\vP^h\vv|) \geq
\langle |\vq^n| -|\vv|,M_\epsilon(w^n)\rangle_{C(\overline{\Omega})}\,.
\label{Qth2cong}
\end{align}
It follows from (\ref{energy3a}) and (\ref{Phstab}) that
\begin{align}
&\left|(M^h_\epsilon(P^h W^n_A)-M_\epsilon(w^n), |\vQ^n_A|-|\vP^h \vv|)\right|
\nonumber \\&\hspace{1in}
\leq |M^h_\epsilon(P^h W^n_A)-M_\epsilon(w^n)|_{0,\infty,\Omega}\,
\left[ C\,\tau_n^{-1} + |\vv|_{0,1,\Omega}\right]\,.
\label{Qth2cond}
\end{align}
Combining (\ref{Qth2cona})--(\ref{Qth2cond}) and (\ref{Mwscon}),
we can pass to the limit for the subsequence as $h \rightarrow 0$
in (\ref{Qth2}), with $\vv^h = \vP^h \vv$
for any fixed $\vv \in [C^\infty(\overline{\Omega})]^d$, to obtain
for $n=1 \rightarrow N$ that
\begin{align}
(\del\,.\,(\vq^n-\vv),w^n) \geq 
\langle |\vq^n|-|\vv|,M_\epsilon(w^n)\rangle_{C(\overline{\Omega})} \qquad 
\forall \ \vv \in [C^\infty(\overline{\Omega})]^d\,.
\label{Qth2conf}
\end{align}
Recalling the density results (\ref{weakj}--c) 
and that $w^n,\,M_\epsilon(w^n) \in C(\overline{\Omega})$, we obtain the
desired result (\ref{Qt2}).   
\endproof

\subsection{Convergence of (Q$^{h,\tau}_{B,r}$) to (Q$^{\tau}$)}
\label{secconvB}

For the purposes of the convergence analysis in this subsection, it is convenient to introduce
the following regularization of (Q$^{\tau}$) for a given $r>1$: 

\noindent
{\bf (Q$^{\tau}_r$)}
For $n = 1 \rightarrow N$, find $w^n_r \in W^{1,p}_0(\Omega)$ and 
$\vq^n_r \in \vV^{r}(\Omega)$ such that
\eqlabstart
\begin{alignat}{2}
\left(\frac{w^n_r-w^{n-1}_r}{\tau_n},\eta\right) + (\del\,.\,\vq^n_r,\eta) &= (f^n 
,\eta)
\quad \;&&\forall \ \eta \in L^2(\Omega),\label{Qt1r}\\
( M_\eps(w^n_r)\,|\vq^n_r|^{r-2}\vq^n_r,\vv) - (w^n_r,\del\,.\,\vv)  &=0
\quad \;
&&\forall \ \vv\in \vV^{r}(\Omega)\,;
\label{Qt2r}
\end{alignat}
\eqlabend
where $w^0_r = w_0^\epsilon$.

\begin{thm}\label{hconthmB}
Let the Assumptions (A1) and (A2) hold.
For any fixed $r \in (1,2)$ and fixed time partition $\{\tau_n\}_{n=1}^N$ with $\tau \in (0,\frac{1}{2}]$,
and for all regular partitionings ${\cal T}^h$ of $\Omega$,
there exists a subsequence of $\{\{W^n_{B,r},\vQ^n_{B,r}\}_{n=1}^N\}_{h>0}$ (not indicated),
where $\{W^n_{B,r},\vQ^n_{B,r}\}_{n=1}^N$ solves (Q$^{h,\tau}_{B,r}$), 
such that as $h \rightarrow 0$, for any $s \in [1,\infty)$,
\eqlabstart
\begin{alignat}{3}
W^{n}_{B,r}&\rightarrow w^n_r  \qquad &&
\mbox{strongly in } L^2(\Omega),
\qquad && n=0 \rightarrow N,
\label{wsconBr}
\\
M_\epsilon^h( W^{n}_{B,r}) &\rightarrow M_\epsilon(w^n_r)  \qquad &&
\mbox{strongly in } L^s(\Omega),
\qquad && n=0 \rightarrow N,
\label{MsconBr}
\\
\vQ^n_{B,r} &\rightarrow \vq^n_r \qquad &&
\mbox{weakly in } [L^r(\Omega)]^d,
\qquad && n=1 \rightarrow N,
\label{qwconBr}\\
\del\,.\,\vQ_{B,r}^n &\rightarrow \del\,.\,\vq^n_{r} \qquad &&
\mbox{weakly in } L^2(\Omega),
\qquad && n=1 \rightarrow N;
\label{hdelqconBr}
\end{alignat}
\eqlabend
where $\{w^n_r,\vq^n_r\}_{n=1}^N$ is a solution of (Q$^\tau_r$), (\ref{Qt1r},b).
\end{thm}
\proof
The desired subsequence weak convergence results (\ref{qwconBr},d) follow immediately from the
bounds on $\{Q^n_{B,r}\}_{n=1}^N$ in (\ref{energy3B},b), on noting that the time
partition $\{\tau_n\}_{n=1}^N$ is fixed.
In addition, we obtain from the first bound in (\ref{energy3B}) that
\begin{alignat}{3}  
W^{n}_{B,r}&\rightarrow w^n_r  \qquad &&
\mbox{weakly in } L^2(\Omega),
\qquad && n=0 \rightarrow N.
\label{wwconBr}
\end{alignat}
Furthermore, we obtain from  
(\ref{Qthr2B}), (\ref{k1epshdef}), (\ref{kepsh1max}),
(\ref{convex}) and (\ref{energy3B}) for $n=1 \rightarrow N$ that
\begin{align}
|(W^n_{B,r}, \del\,.\,\vv^h)| &=  
|(M_\epsilon^h(W^n_{B,r})\,|\vQ^n_{B,r}|^{r-2} \vQ^n_{B,r},
\vv^h)^h| \leq 
k^{\eps,h}_{1,\infty}
(|\vQ^n_{B,r}|^{r-1},
|\vv^h|)^h
\nonumber \\
& \leq C\,
[(|\vQ^n_{B,r}|^r,1)^h]^{\frac{r-1}{r}}\,[(|\vv^h|^r,1)^h]^{\frac{1}{r}}
\leq C
\,|\vQ^n_{B,r}|_{0,r,\Omega}^{r-1}\,|\vv^h|_{0,r,\Omega}
\nonumber \\
& \leq C(\tau_n^{-1})\,|\vv^h|_{0,r,\Omega}
\qquad \hspace{1in} \forall \ \vv^h \in \vV^h
\,.
\label{gradhWnBr}
\end{align}
For any fixed $\vv \in [C^\infty(\overline{\Omega})]^d$, on choosing $\vv^h = \vI^h \vv$ in 
(\ref{gradhWnBr}), letting $h \rightarrow 0$ and noting (\ref{interpint}), 
(\ref{wwconBr}) and
(\ref{interp1}), we obtain that 
\begin{align}
|(w^n_r, \del\,.\,\vv)| \leq C(\tau_n^{-1})\,|\vv|_{0,r,\Omega}\,, \qquad n=1 \rightarrow N\,.
\label{gradwnr}
\end{align}  
Repeating (\ref{gradwnr}) for all $\vv \in [C^\infty(\overline{\Omega})]^d$
and as $C^\infty(\overline{\Omega})$ is dense in $L^r(\Omega)$, we obtain
that 
\begin{align}
w^n_r \in W^{1,p}_0(\Omega) \qquad \mbox{with} \qquad \|w^n_r\|_{1,p,\Omega} \leq C(\tau_n^{-1}),
\qquad n=1 \rightarrow N\,.   
\label{wnrW1p}
\end{align}
The fact that $w^n_r$ vanishes on $\partial \Omega$ can be deduced from (\ref{gradwnr})
by using an argument similar to that in \cite[p.\ 699]{BP5}. 

Next, for $n=1 \rightarrow N$, we introduce $\del_h W^n_{B,r} \in \vV^h$ such that
\begin{align}
(\del_h W^n_{B,r}, \vv^h) = - (W^n_{B,r}, \del\,.\,\vv^h) \qquad \forall \ \vv^h \in \vV^h\,.
\label{delhWnBr}
\end{align}
It follows from (\ref{delhWnBr}) and (\ref{gradhWnBr}) that 
\begin{align}
|\del_h W^n_{B,r}|_{0,\Omega} \leq C(\tau_n^{-1}), \qquad n=1 \rightarrow N\,.
\label{delhWnBrL2}
\end{align}
For $n=1 \rightarrow N$, we now introduce $\hat W^n_{B,r} \in U^h_0$ such that
\begin{align}
(\del \hat W^n_{B,r}, \del \eta^h) 
= (\del_h W^n_{B,r}, \del \eta^h) \qquad \forall \ \eta^h \in U^h_0\,.
\label{hatWnBr}
\end{align}
It follows from (\ref{Poin}), (\ref{hatWnBr}) and (\ref{delhWnBrL2}) that for 
$n=1 \rightarrow N$ 
\begin{align}
\| \hat W^n_{B,r}\|_{1,\Omega} \leq C\,|\del \hat W^n_{B,r}|_{0,\Omega}
\leq C\,|\del_h W^n_{B,r}|_{0,\Omega} \leq C(\tau_n^{-1})\,.
\label{hatWnBrH1}
\end{align}
We deduce from (\ref{delhWnBrL2}) and (\ref{hatWnBrH1}) that there exists 
a further subsequence of $\{\{\del_h W^n_{B,r}, \hat W^n_{B,r}\}_{n=1}^N\}_{h>0}$
(not indicated) such that as $h \rightarrow 0$, for any $s \in [1,\infty)$, 
\eqlabstart
\begin{alignat}{3}  
\del_h W^{n}_{B,r}&\rightarrow \vd^n_r  \qquad &&
\mbox{weakly in } [L^2(\Omega)]^d,
\qquad && n=1 \rightarrow N\,,
\label{wwcondelBr}\\ 
\del \hat W^{n}_{B,r}&\rightarrow \del \hat w^n_r  \qquad &&
\mbox{weakly in } [L^2(\Omega)]^d,
\qquad && n=1 \rightarrow N\,,
\label{wwcondelhatBr}\\ 
\hat W^{n}_{B,r}&\rightarrow \hat w^n_r  \qquad &&
\mbox{strongly in } L^s(\Omega),
\qquad && n=1 \rightarrow N\,;
\label{wsconhatBr}
\end{alignat}
\eqlabend
where $\hat w^n_r \in H^1_0(\Omega)$.
For any fixed $\vv \in [C^\infty(\overline{\Omega})]^d$, on choosing $\vv^h = \vI^h \vv$ in 
(\ref{delhWnBr}), letting $h \rightarrow 0$ for the subsequence
and noting (\ref{interpint}),
(\ref{wwcondelBr}), 
(\ref{wwconBr}) and
(\ref{interp1})  
yields that
\begin{align}
(\vd^n_r, \vv) = -(w^n_r, \del\,.\,\vv) \qquad n=1 \rightarrow N\,. 
\label{defdnr}
\end{align} 
Repeating (\ref{defdnr}) for all $\vv \in [C^\infty(\overline{\Omega})]^d$ yields that 
$\vd^n_r = \del w^n_r$.
Similarly, for any fixed $\eta \in C^\infty_0(\Omega)$, on choosing $\eta^h = \pi^h \eta$
in (\ref{hatWnBr}), letting $h \rightarrow 0$ for the subsequence and noting (\ref{piconv}), 
({\ref{wwcondelBr},b) and $\vd^n_r = \del w^n_r$ yields that   
\begin{align}
(\del \hat w^n_r , \del \eta ) = (\vd^n_r, \del \eta) = (\del w^n_r, \del \eta) 
\qquad n=1 \rightarrow N\,. 
\label{defhatwnr}
\end{align}
Repeating (\ref{defhatwnr}) for all $\eta \in C^\infty_0(\Omega)$ yields that $\hat w^n_r = 
w^n_r$.    

For $n=1 \rightarrow N$, let $z^n$ be such that 
\begin{align}
- \Delta z^n = \hat W^n_{B,r}- W^n_{B,r} \quad \mbox{in } \Omega, 
\qquad z^n=0 \quad \mbox{on } \partial \Omega. 
\label{zn}
\end{align}
As $\Omega$ is convex 
polygonal, elliptic regularity yields that 
\begin{align}
\|z^n\|_{2,\Omega} \leq C\,| \hat W^n_{B,r}- W^n_{B,r}|_{0,\Omega}.
\label{znreg}
\end{align}
It follows from (\ref{zn}), (\ref{hatWnBr}), (\ref{interpint}),  (\ref{delhWnBr}), 
(\ref{hatWnBrH1}), (\ref{piconv0}),
(\ref{interp1}) and (\ref{znreg}) that for $n=1 \rightarrow N$ 
\begin{align}
&| \hat W^n_{B,r}-W^n_{B,r}|_{0,\Omega}^2 
\nonumber \\
&\qquad = (\del  \hat W^n_{B,r}, \del z^n) +
( W^n_{B,r}, \Delta z^n) \nonumber \\
&\qquad =(\del  \hat W^n_{B,r}, \del (z^n-\pi^h z^n)\,) + (\del_h W^n_{B,r}, \del [\pi^h z^n])
+( W^n_{B,r}, \Delta z^n) \nonumber \\
&\qquad =(\del  \hat W^n_{B,r} -\del_h  W^n_{B,r}, \del (z^n-\pi^h z^n)\,) 
+ (\del_h W^n_{B,r}, \del z^n)
+( W^n_{B,r}, \Delta z^n) \nonumber \\
&\qquad =(\del  \hat W^n_{B,r} -\del_h  W^n_{B,r}, \del (z^n-\pi^h z^n)\,) 
+ (\del_h W^n_{B,r}, \del z^n -\vI^h(\del z^n)\,) \nonumber \\
&\qquad \leq C(\tau_n^{-1})\,\left[\,|\del (z^n-\pi^h z^n)|_{0,\Omega} 
+ |\del z^n -\vI^h(\del z^n)|_{0,\Omega} \right] \nonumber \\
&\qquad \leq C(\tau_n^{-1})\,h\,\|z\|_{2,\Omega} \leq C(\tau_n^{-1})\,h^2.
\label{znerr}
\end{align}
As $\hat w^n_r =w^n_r$, $n=1 \rightarrow N$,
it follows from (\ref{znerr}) and (\ref{wsconhatBr}) that
the desired result (\ref{wsconBr}) holds.
 
We deduce from (\ref{wsconBr}), (\ref{Mepsmm}) and (\ref{minmax})
for a further subsequence of $\{\{W^n_{B,r}\}_{n=0}^N\}_{h>0}$ (not indicated)
that as $h \rightarrow 0$, for $n=0 \rightarrow N$,
\begin{align}
W^n_{B,r} \rightarrow w^n_r \quad \mbox{a.e.\ in } \Omega
\quad \Rightarrow \quad M_\epsilon(W^n_{B,r}) \rightarrow M_\epsilon(w^n_r) 
\quad \mbox{a.e.\ in } \Omega\,.
\label{MWwae}
\end{align} 
It follows from (\ref{MWwae}), (\ref{Mepsineq}), (\ref{k1max}) and the Lebesgue's general 
convergence theorem that as $h \rightarrow 0$  
for any $s \in [1,\infty)$
\begin{align}
M_\eps(W^n_{B,r}) \rightarrow M_\epsilon(w^n_r) \qquad \mbox{strongly in } L^s(\Omega),
\qquad n=0 \rightarrow N\,.
\label{MsconBrs}
\end{align} 
Combining (\ref{Mepshapprox}), (\ref{Pconv}),  
(\ref{piconv})
and (\ref{MsconBrs}) yields the desired result (\ref{MsconBr}).
 
We now need to establish that $\{w^n_r,\vq^n_r\}_{n=1}^N$ solve (Q$^\tau_r$), (\ref{Qt1r},b).
For any $\eta \in C^\infty_0(\Omega)$, we choose  $\eta^h = P^h \eta$ in (\ref{Qthr1B})
and now pass to the limit $h \rightarrow 0$ for the subsequence,
on noting (\ref{wsconBr},d)
and (\ref{Pconv}),
to obtain (\ref{Qt1r}) for all $\eta \in C^\infty_0(\Omega)$.    
Noting that $C^\infty_0(\Omega)$ is dense in $L^2(\Omega)$ 
and that $w^n_r,\,\del\,.\,\vq^n_r,\,f^n \in L^2(\Omega)$ 
yields the desired (\ref{Qt1r}).

For any $\vv \in [C^\infty(\overline{\Omega})]^d$, we choose  $\vv^h = \vQ^n_{B,r}-\vI^h \vv$ 
in (\ref{Qthr2B})
and now try to pass to the limit for the subsequence as $h \rightarrow 0$. 
First, we note from (\ref{aederiv}) and (\ref{convex}) that
for $n=1 \rightarrow N$ 
\begin{align}
(W^n_{B,r}, \del \,.\,(\vQ^n_{B,r}- \vI^h \vv)\,) &= 
(M^h_\eps(W^n_{B,r})\,|\vQ^n_{B,r}|^{r-2}\,\vQ^n_{B,r},\vQ^n_{B,r}-\vI^h \vv)^h 
\nonumber \\
& \geq \frac{1}{r} \, (M^h_\eps(W^n_{B,r}),|\vQ^n_{B,r}|^r-|\vI^h \vv|^r)^h 
\nonumber \\
& \geq \frac{1}{r} \, (M^h_\eps(W^n_{B,r}),|\vQ^n_{B,r}|^r-|\vI^h \vv|^r)
\nonumber \\
& \quad + \frac{1}{r} \, \left[ (M^h_\eps(W^n_{B,r}),|\vI^h \vv|^r)-
(M^h_\eps(W^n_{B,r}),|\vI^h \vv|^r)^h \right]
\,. 
\label{Qthr2BtoQtr2a}
\end{align}
Once again, it follows from (\ref{aederiv}) that 
\begin{align}
\frac{1}{r} \, (M^h_\eps(W^n_{B,r}),|\vQ^n_{B,r}|^r-|\vI^h \vv|^r)
& \geq  
(M^h_\eps(W^n_{B,r}),|\vI^h \vv|^{r-2} \vI^h \vv,
\vQ^n_{B,r}-\vI^h \vv)\,.
\label{Qthr2BtoQtr2b}
\end{align}
In addition, it follows from (\ref{k1epshdef}), (\ref{kepsh1max}) and (\ref{nih}) that 
\begin{align}
\frac{1}{r} \, \left|(M^h_\eps(W^n_{B,r}),|\vI^h \vv|^r)-
(M^h_\eps(W^n_{B,r}),|\vI^h \vv|^r)^h \right| \leq C\,h\,\|\vv\|_{1,\infty,\Omega}\,.
\label{Qthr2BtoQtr2c}
\end{align} 
Combining (\ref{Qthr2BtoQtr2a}) and (\ref{Qthr2BtoQtr2b}),
and passing to the limit $h \rightarrow 0$ for the subsequence yields,
on noting (\ref{interpint}), (\ref{wsconBr}--d), (\ref{interp1}) 
and (\ref{Qthr2BtoQtr2c}), yields for $n=1 \rightarrow N$ that
\begin{align}
(w^n_r, \del\,.\,(\vq^n_r-\vv) \geq (M_\eps(w^n_r) \,|\vv|^{r-2}\,\vv, \vq^n_r - \vv)
\qquad \forall \ \vv \in [C^\infty(\overline{\Omega})]^d\,.
\label{Qthr2BtoQtr2d}
\end{align}  
As $w^n_r,\,M_\eps(w^n_r) \in C(\overline{\Omega})$, $\vq^n_r \in \vV^r(\Omega)$ and 
$f^n \in L^2(\Omega)$, it follows from (\ref{densevVs}) that 
(\ref{Qthr2BtoQtr2d}) holds true for all $\vv \in \vV^r(\Omega)$.  
For any fixed $\vz \in \vV^r(\Omega)$, choosing $\vv = \vq^n_r \pm \alpha \vz$ with 
$\alpha \in {\mathbb R}_{>0}$ in (\ref{Qthr2BtoQtr2d}) and letting $\alpha \rightarrow 0$
yields the desired result (\ref{Qt2r}) on repeating the above for any 
$\vz \in \vV^r(\Omega)$. Hence $\{w^n_r,\vq^n_r\}_{n=1}^N$ is a solution of (Q$^\tau_r$),
(\ref{Qt1r},b).
\endproof    
 
\begin{thm}\label{hconthmr}
Let the Assumptions (A1) and (A2) hold.
For any fixed time partition $\{\tau_n\}_{n=1}^N$ with $\tau \in (0,\frac{1}{2}]$,
there exists a subsequence of $\{\{w^n_{r},\linebreak
\vq^n_{r}\}_{n=1}^N \}_{r>1}$ (not indicated),
where $\{w^n_{r},\vq^n_{r}\}_{n=1}^N$ solves (Q$^{\tau}_{r}$), 
such that as $r \rightarrow 1$ 
\eqlabstart
\begin{alignat}{3}
w^{n}_{r}&\rightarrow w^n  \qquad &&
\mbox{strongly in } C(\overline{\Omega}),
\qquad && n=0 \rightarrow N,
\label{wsconr}
\\
M_\epsilon( w^{n}_{r}) &\rightarrow M_\epsilon(w^n)  \qquad &&
\mbox{strongly in } C(\overline{\Omega}),
\qquad && n=0 \rightarrow N,
\label{Msconr}
\\
\vq^n_{r} &\rightarrow \vq^n \qquad &&
\mbox{vaguely in } [{\cal M}(\overline{\Omega})]^d,
\qquad && n=1 \rightarrow N,
\label{qwconr}\\
\del\,.\,\vq_{r}^n &\rightarrow \del\,.\,\vq^n \qquad &&
\mbox{weakly in } L^2(\Omega),
\qquad && n=1 \rightarrow N;
\label{hdelqconr}
\end{alignat}
\eqlabend
where $\{w^n,\vq^n\}_{n=1}^N$ is a solution of (Q$^\tau$), (\ref{Qt1},b).
\end{thm}
\proof 
It follows immediately from (\ref{energy3B},b), (\ref{aederiv}) and (\ref{wsconBr},c,d) that
\eqlabstart
\begin{align}
\max_{n=0 \rightarrow N}
|w^n_{r}|_{0,\Omega} + \sum_{n=1}^N |w^n_{r}-w^{n-1}_{r}|_{0,\Omega}^2 +
\sum_{n=1}^N \tau_n\,|\vq^n_{r}|_{0,r,\Omega}^r
&\leq C\,,
\label{energy3Br} \\
\sum_{n=1}^N \tau_n^2\,|\del\,.\,\vq^n_{r}|_{0,\Omega}^2
&\leq C\,.
\label{energy2Br}
\end{align}
\eqlabend
The desired convergence results (\ref{wsconr}--d)
follow immediately from (\ref{energy3Br},b) and (\ref{wnrW1p}) on recalling that 
the embedding $W^{1,p}(\Omega) \hookrightarrow C(\overline{\Omega})$  
is compact for $p>d$, $M_\epsilon : C(\overline{\Omega}) \rightarrow C(\overline{\Omega})$.
One can immediately pass to the limit $r \rightarrow 1$ for the subsequence 
in (\ref{Qt1r}), on noting 
(\ref{wsconr},d), 
to obtain (\ref{Qt1}).
Similarly to (\ref{rlim2}),
choosing $\vv = \vq^n_{r} - \vpsi$ in (\ref{Qt2r}) and noting (\ref{aederiv}),
(\ref{Mepsineq}) and (\ref{energy3Br}),
one obtains for $n = 1 \rightarrow N$ that
\begin{align}
(w^n_{r},\del\,.\,(\vq^n_{r}-\vpsi)\,) &= (M_\epsilon(w^n_{r})\,|\vq^n_{r}|^{r-2} \vq^n_{r},
\vq^n_{r}-\vpsi)  
\nonumber \\
&\geq (M_\epsilon(w^n),|\vq^n_r|) 
- C(\tau_n^{-1})\, |M_\epsilon(w_n)-M_\epsilon(w^n_r)|_{0,\infty,\Omega}
\nonumber \\
& \quad -\frac{1}{r}\,(M_\epsilon(w^n_{r}),|\vpsi|^r)
+ \frac{1-r}{r} \, k_{1,\infty}^{\epsilon}\,|\Omega|
\qquad \forall \ \vpsi \in \vV^r(\Omega)\,.
\label{rlim2r}
\end{align}
Noting (\ref{wsconr}--d) and (\ref{Mweak2}), 
one can pass to the limit $r \rightarrow 1$ for the subsequence in (\ref{rlim2r}) to obtain
(\ref{Qt2}). Hence $\{w^n,\vq^n\}_{n=1}^N$ solves (Q$^{\tau}$), (\ref{Qt1},b).
\endproof 

\begin{rem} \label{splitproof}
It appears necessary to split the convergence proof 
of solutions of (Q$^{h,\tau}_{B,r}$) to solutions of (Q$^\tau$), as $h \rightarrow 0$
and $r \rightarrow 1$, by first considering the limit $h \rightarrow 0$ to solutions
of (Q$^\tau_r$), then 
the limit $r \rightarrow 1$ to solutions of (Q$^\tau$). 
Similarly, it does not appear possible to directly  
prove convergence of solutions of (Q$^{h,\tau}_{B}$) to solutions of (Q$^\tau$), 
as $h \rightarrow 0$.
For example, if we attempted the latter, we would still only be able to show
$M_\eps(W^n_B) \rightarrow M_\eps(w^n)$  
strongly in $L^s(\Omega)$ for $s \in [1,\infty)$, as $h \rightarrow 0$;
and this is not adequate to pass to the limit $h \rightarrow 0$ in (Q$^{h,\tau}_B$), 
(\ref{Qth1B},b). 
\end{rem}

\subsection{Convergence of (Q$^{\tau}$) to (Q)}
\label{secconvtau}
First we note the following result.

\begin{thm}
\label{QPtaueq}
Let the Assumptions (A1) and (A2) hold.
If $\{w^n,\vq^n\}_{n=1}^N$ is a solution of (Q$^{\tau}$), (\ref{Qt1},b), then
$\{w^n\}_{n=1}^N$ solves (P$^{\tau}$), (\ref{Pt}), and
\begin{align}
w^n \geq w^{n-1} \qquad n=1 \rightarrow N\,.
\label{wnmono}
\end{align}
\end{thm}
\proof
Similarly to (\ref{Wneq},b), we deduce on choosing $\vv=\vzero$ and $2\,\vq^n$ in (\ref{Qt2})
that
\eqlabstart
\begin{align}
\langle |\vq^n|, M_\epsilon(w^n) \rangle_{C(\overline{\Omega})} 
&= (\del\,.\,\vq^n,w^n)
\label{wneq}\\
\mbox{and hence that} \quad \langle |\vv|, M_\epsilon(w^n) \rangle_{C(\overline{\Omega})} 
&\geq (\del\,.\,\vv,w^n) \quad \forall \ \vv \in \vV^{\cal M}(\Omega)\,.
\label{wnineq}
\end{align}
\eqlabend
Noting that $[C^\infty(\overline{\Omega})]^d \subset \vV^{\cal M}(\Omega)$ and $w^n \in W^{1,\infty}_0(\Omega)
\subset C(\overline{\Omega})$, we deduce from (\ref{wnineq}) that
\begin{align}
(M_\epsilon(w^n),|\vv| ) 
\geq -(\del w^n, \vv) \qquad \forall \ \vv \in [C^\infty(\overline{\Omega})]^d\,,
\label{wnineqa}
\end{align}
It follows from (\ref{wnineqa}) that for $n=1 \rightarrow N$
\begin{align}
|\del w^n| \leq M_\epsilon(w^n) \quad \mbox{a.e.\ on } \Omega
\qquad \Rightarrow \qquad w^n \in K(w^n)
\,;
\label{wnineqb} 
\end{align}
see, for example, the argument in \cite[p.\ 698]{BP5}.

Similarly to (\ref{Qth1ineq}), choosing $\eta=\varphi -w^n$ for any $\varphi \in  K(w^n)$ 
in (\ref{Qt1}), we obtain, on noting (\ref{wneq}), employing a sequence of the 
form (\ref{weakj}--c) (with $\vv$ and $\{\vv_j\}_{j \geq1}$
replaced by $\vq^n$ and $\{\vq^n_j\}_{j\geq 1}$, respectively) and (\ref{K}), that 
\begin{align}
&\left(\frac{w^n-w^{n-1}}{\tau_n}-f^n,\varphi-w^n\right) 
\nonumber \\
& \hspace{0.3in}= -(\del\,.\,\vq^n,\varphi-w^n) 
= \langle |\vq^n|, M_\epsilon(w^n)\rangle_{C(\overline{\Omega})}
- \lim_{j \rightarrow \infty} (\del\,.\,\vq^n_j,\varphi)
\nonumber \\
&\hspace{0.3in}
= \langle |\vq^n|, M_\epsilon(w^n)\rangle_{C(\overline{\Omega})}
+ \lim_{j \rightarrow \infty} (\vq^n_j,\del \varphi)
\nonumber \\
&\hspace{0.3in}
\geq 
\langle |\vq^n|, M_\epsilon(w^n)\rangle_{C(\overline{\Omega})}
-\lim_{j \rightarrow \infty}\langle |\vq^n_j|,M_\epsilon(w^n)
\rangle_{C(\overline{\Omega})}= 0 \,.
\label{wnineqc}
\end{align} 
Hence $\{w^n\}_{n=0}^N$ solves (P$^{\tau}$), (\ref{Pt}). 

Let $\eta = w^n + [w^{n-1}-w^n]_+$, where $[s]_+ : =\max(s,0)$ for any $s\in \mathbb{R}$.
It follows from (\ref{wnineqb}) and (\ref{Mepsineq}) that for a.e.\ $\vx \in \Omega$ 
\begin{alignat}{2}
w^n(\vx) \geq w^{n-1}(\vx) 
\quad &\Rightarrow \quad
|\del \eta (\vx)| &&= |\del w^n(\vx)| \leq M_\epsilon(w^n(\vx))\,,
\nonumber \\
 w^{n-1}(\vx) \geq w^{n}(\vx) 
\quad &\Rightarrow \quad
|\del \eta (\vx)| &&= |\del w^{n-1}(\vx)| \leq M_\epsilon(w^{n-1}(\vx))
\nonumber \\
&&&\leq M_\epsilon(w^n(\vx))
\,.
\label{wnineqd}
\end{alignat}
Hence $\eta =w^n + [w^{n-1}-w^n]_+ \in K(w^n)$. Substituting this into (\ref{Pt}), and 
recalling that the source $f^n \geq 0$ yields for $n=1 \rightarrow N$ that 
\begin{align}
| [w^{n-1}-w^{n}]_+\,|_{0,\Omega}^2 \leq - \tau_n \,( f^n, [w^{n-1}-w^{n}]_+) \leq 0\,, 
\label{wnineqe}
\end{align}
and hence the desired result (\ref{wnmono}).
\endproof

\begin{rem} \label{remmono}
We note that the monotonicity result (\ref{wnmono}) for $\{w^n\}_{n=1}^N$ solving 
(Q$^\tau) \equiv ($P$^\tau$)  does not hold for 
$\{W^n_{A,r}\}_{n=1}^N$ solving (Q$^{h,\tau}_{A,r}) \equiv ($P$^{h,\tau}_{A,p}$),
\linebreak
$\{W^n_A\}_{n=1}^N$ solving (Q$^{h,\tau}_A) \equiv ($P$^{h,\tau}_A$),
$\{W^n_{B,r}\}_{n=1}^N$ solving (Q$^{h,\tau}_{B,r}$)
and $\{W^n_{B}\}_{n=1}^N$ solving (Q$^{h,\tau}_{B}$).
\end{rem}

We introduce the following notation for $t \in (t_{n-1},t_n]$, $n = 1 \rightarrow N$,
\begin{align}
f^{\tau,+}(\cdot,t) &:= f^n(\cdot)\, \qquad
w^{\tau}(\cdot,t) := \frac{(t-t_{n-1})}{\tau_n}\,w^n(\cdot)
+\frac{(t_n-t)}{\tau_n}\,w^{n-1}(\cdot)\,,
\nonumber \\
w^{\tau,+}(\cdot,t) &:= w^n(\cdot),
\qquad
w^{\tau,-}(\cdot,t) := w^{n-1}(\cdot),
\qquad
\vq^{\tau,+}(\cdot,t) := \vq^n(\cdot)\,.
 \label{wqtau+}
\end{align}
In addition, we write $w^{\tau(,\pm)}$ to mean with or without the superscripts $\pm$.
We note from (\ref{wqtau+}) and (\ref{fn}) that    
\begin{align}
f^{\tau,+} \rightarrow f \qquad \mbox{strongly in } L^s(0,T;L^2(\Omega))
\mbox{ as } \tau \rightarrow 0\,;
\label{ftauconv}
\end{align}
where $s=2$ if Assumption (A1) holds, and $s=\infty$ if (A3) holds.

Adopting the notation (\ref{wqtau+}), (Q$^{\tau}$), (\ref{Qt1},b), can be rewritten as: 
Find $w^{\tau} 
\in L^\infty(0,T;W^{1,\infty}_0(\Omega))\cap W^{1,\infty}(0,T;L^2(\Omega)
)$ and 
$\vq^{\tau,+} \in L^\infty(0,T;$ $[{\cal M}(\overline{\Omega})]^d)$ such that
\eqlabstart
\begin{align}
&\int_0^T \left[ 
\left(\frac{\partial w^{\tau}}{\partial t},\eta \right) 
- \langle \vq^{\tau,+},\del \eta\rangle_{C(\overline{\Omega})} -(f^{\tau,+} 
,\eta) \right] \,{\rm d}t =0
\nonumber \\ &\hspace{3.1in}
\forall \ \eta \in L^1(0,T;C^1_0(\overline{\Omega})),\label{Q1tau}\\
&\int_0^T \left[ \langle |\vv|-|\vq^{\tau,+}|,M_\epsilon(w^{\tau,+})
\rangle_{C(\overline{\Omega})} - 
(\del\,.\,\vv-f^{\tau,+},
w^{\tau,+})\right] \,{\rm d}t  \nonumber \\
& \hspace{0.49in} \geq \tfrac{1}{2} \left[\, |w^{\tau}(\cdot,T)|^2_{0,\Omega} - |w^\epsilon_0(\cdot)|^2_{0,\Omega}
\,\right]
\qquad \forall \ \vv\in L^1(0,T;\vV^{\cal M}(\Omega))\,;
\label{Q2tau}
\end{align}
\eqlabend
where $w^\tau(\cdot,0) = w_0^\epsilon(\cdot)$.

(\ref{Q1tau}) is obtained from (\ref{Qt1}) by choosing $\eta(\cdot) = 
\int^{t_n}_{t_{n-1}} 
\chi(\cdot,t) \,{\rm d}t$ in (\ref{Qt1}) and summing from $n=1 \rightarrow N$,
for any $\chi \in L^1(0,T;C^1_0(\overline{\Omega}))$, and 
noting that      
\begin{align}
\int_0^T (\del\,.\,\vq^{\tau,+},\chi)\,{\rm d}t = - 
\int_0^T \langle \vq^{\tau,+},\del \chi\rangle_{C(\overline{\Omega})}\,{\rm d}t
\qquad \forall \ \chi \in L^1(0,T;C^1_0(\overline{\Omega}))\,.
\label{Qt1equiv}
\end{align}
Similarly, (\ref{Q2tau}) is obtained from (\ref{Qt2}) by choosing $\vv(\cdot) = 
\frac{1}{\tau_n}\,\int^{t_n}_{t_{n-1}} 
\vpsi(\cdot,t) \,{\rm d}t$ in (\ref{Qt2}), multiplying by $\tau_n$ 
and summing from $n=1 \rightarrow N$,
for any $\vpsi \in L^1(0,T;\vV^{{\cal M}}(\Omega))$, and 
noting from (\ref{Qt1}) and (\ref{simpid}) that
\begin{align}
&- \sum_{n=1}^N \tau_n \,(\del\,.\,\vq^n,w^n) 
\nonumber \\
& \qquad = 
\tfrac{1}{2} \left[ \,|w^N|_{0,\Omega}^2 - |w^\eps_0|_{0,\Omega}^2 \,\right]
+ \sum_{n=1}^N \left[ 
\tfrac{1}{2} 
\,|w^n-w^{n-1}|_{0,\Omega}^2 - \tau_n \,(f^n,w^n)  \right]
\nonumber \\
& \qquad \geq  
\tfrac{1}{2} \left[\, |w^N|_{0,\Omega}^2 - |w^\eps_0|_{0,\Omega}^2 \,\right]
- \sum_{n=1}^N \tau_n \,(f^n,w^n)
\,. 
\label{Qt2equiv}
\end{align}

\begin{thm}\label{conthm}
Let the Assumptions (A1), (A2) and (A3) hold.
For all time partitions $\{\tau_n\}_{n=1}^N$,
there exists a subsequence of $\{\{w^n,\vq^n\}_{n=1}^N\}_{\tau>0}$ (not indicated),
where $\{w^n,\vq^n\}_{n=1}^N$ solves (Q$^{\tau}$), 
such that as $\tau \rightarrow 0$
\eqlabstart
\begin{alignat}{2}
w^{\tau},\,w^{\tau,\pm}&\rightarrow w  \qquad &&
\mbox{weak$^\star$ in } L^\infty(0,T;W^{1,\infty}(\Omega)),
\label{wgradwcont}\\
\frac{\partial w^{\tau}}{\partial t}&\rightarrow \frac{\partial w}{\partial t}  \qquad &&
\mbox{vaguely in } L^\infty(0,T;[C^{1}_0(\overline{\Omega})]^\star),
\label{wtwcont}\\
w^{\tau} &\rightarrow w  \qquad &&
\mbox{strongly in } C([0,T];C(\overline{\Omega})),
\label{wscont}
\\
w^{\tau,\pm}&\rightarrow w  \qquad &&
\mbox{strongly in } L^2(0,T;C(\overline{\Omega})),
\label{wsL2}
\\
M_\epsilon(w^{\tau})&\rightarrow M_\epsilon(w)  \qquad &&
\mbox{strongly in } C([0,T];C(\overline{\Omega})),
\label{Mwscont}
\\
M_\epsilon(w^{\tau,\pm}) &\rightarrow M_\epsilon(w)  \qquad &&
\mbox{strongly in } L^2(0,T;C(\overline{\Omega})),
\label{MwsL2}
\\
\vq^{\tau,+} &\rightarrow \vq \qquad &&
\mbox{vaguely in } L^\infty(0,T;[{\cal M}(\overline{\Omega})]^d);
\label{qwcont}
\end{alignat}
\eqlabend
where $\{w,\vq\}$ is a solution of (Q), (\ref{Q1},b).
Moreover, $w$ solves (P), (\ref{P}).
\end{thm}
\proof
It follows from (\ref{Poin}), (\ref{wnineqb}), (\ref{Mepsineq}) and (A1) that
\begin{align}
\max_{n=0 \rightarrow N} \|w^n\|_{1,\infty,\Omega} \leq C\,.
\label{wn1inf}
\end{align}
Choosing $\eta=w^n$ in (\ref{Qt1}), summing from $n=1 \rightarrow N$ and  
noting (\ref{wneq}), (\ref{simpid}), (\ref{fnsum}) and (\ref{wn1inf}) yields that 
\begin{align}
&|w^N|_{0,\Omega}^2 + \sum_{n=1}^N |w^n-w^{n-1}|_{0,\Omega}^2 +
2 \,\sum_{n=1}^N \tau_n\,\langle |\vq^n|,M_\eps(w^n) \rangle_{C(\overline{\Omega})}
\nonumber \\
&\hspace{0.5in} = |w_0^\eps|^2_{0,\Omega} + 2 \sum_{n=1}^N \tau_n\,(f^n,w^n)
\nonumber \\
& \hspace{0.5in} \leq |w_0^\eps|^2_{0,\Omega}+
2 \left(\sum_{n=1}^N \tau_n\,|f^n|_{0,\Omega}^2\right)^{\frac{1}{2}}
\left(\sum_{n=1}^N \tau_n\,|w^n|_{0,\Omega}^2 \right)^{\frac{1}{2}}
\leq C\,.
\label{wnbasic}
\end{align}
The bounds (\ref{wn1inf}) and (\ref{wnbasic}) only assume the Assumptions (A1) and (A2).

Choosing $\eta=w^n$ in (\ref{Qt1}), and 
noting (\ref{wneq}), 
(\ref{wnmono}) and (A3), yields for $n=1 \rightarrow N$ that
\begin{align}
\langle |\vq^n|, M_\epsilon(w^n) \rangle_{C(\overline{\Omega})} &= (\del\,.\,\vq^n,w^n)=
\left( f^n - \frac{w^n-w^{n-1}}{\tau_n}, w^n \right)
\nonumber \\
& \leq |f^n|_{0,\Omega}\,|w^n|_{0,\Omega} \leq C\,. 
\label{vqnmax}
\end{align}
Therefore (\ref{vqnmax}) and (\ref{Mepsineq}) 
yield that
\begin{align}
\max_{n=1 \rightarrow N} \int_{\overline{\Omega}} |\vq^n| \leq C\,.
\label{vqnmaxa}
\end{align}
We obtain from (\ref{Qt1}), (\ref{vqnmaxa}) and (A3) for $n=1 \rightarrow N$ that
\begin{align}
\left|\left( \frac{w^n-w^{n-1}}{\tau_n} ,\eta\right) \right|
&= \left| (f^n- \del\,.\,\vq^n, \eta) \right|
= \left| (f^n,\eta)  + \langle \vq^n, \del \eta\rangle_{C(\overline{\Omega})} \right|
\nonumber \\
& \leq \left[ |f^n|_{0,1,\Omega} + \int_{\overline{\Omega}} |\vq^n| \right]  
\|\eta\|_{1,\infty,\Omega}
\nonumber \\
& \leq C\,\|\eta\|_{1,\infty,\Omega}
\qquad \qquad \qquad \forall \ \eta \in C^1_0(\overline{\Omega})\,.
\label{wtbd1}
\end{align}

Combining the bounds (\ref{wn1inf}), (\ref{wnbasic}), (\ref{vqnmaxa}) and (\ref{wtbd1}),
we obtain, on adopting the notation (\ref{wqtau+}),  that
\begin{align}
&\|w^{\tau(,\pm)}\|_{L^\infty(0,T;W^{1,\infty}(\Omega))}
+\left \|\frac{\partial w^{\tau}}{\partial t} \right \|_{L^\infty(0,T;[C^1_0(\overline{\Omega})]^{\star})}
\nonumber \\
& \hspace{0.5in} +\frac{1}{\tau} \, \|w^{\tau,+}-w^{\tau,-}\|_{L^2(0,T;L^2(\Omega))}^2
+ \|\vq^{\tau,+}\|_{L^\infty(0,T;[{\cal M}(\overline{\Omega})]^{d})}
\leq C\,.
\label{wtbd2}
\end{align}
It follows from (\ref{wqtau+}) and 
the third bound in (\ref{wtbd2}) that
\begin{align}
\|w^{\tau,+}-w^{\tau,-}\|_{L^2(0,T;L^2(\Omega))}^2 + 
\|w^{\tau}-w^{\tau,\pm}\|_{L^2(0,T;L^2(\Omega))}^2 \leq C\,\tau\,.
\label{wdiff}
\end{align}
The subsequence convergence results 
(\ref{wgradwcont},b) and (\ref{qwcont}) follow immediately from
the bounds (\ref{wtbd2}) and (\ref{wdiff}).
To apply (\ref{compact1}) to $w^\tau$, we first note that $C^1_0(\overline{\Omega})$ 
is not a reflexive 
Banach space. However, $W^{2,s}_0(\Omega)$, the closure of $C^\infty_0(\Omega)$ for the norm 
$\|\cdot\|_{2,s,\Omega}$, with $s \in (d,\infty)$ is a reflexive Banach space such that
$W^{2,s}_0(\Omega) \subset C^1_0(\overline{\Omega})$. 
Hence, the first two bounds in (\ref{wtbd2}) 
yield for $s \in (d,\infty)$ that
\begin{align}
&\|w^{\tau}\|_{L^\infty(0,T;W^{1,\infty}(\Omega))}
+\left \|\frac{\partial w^{\tau}}{\partial t} \right \|_{L^\infty(0,T;[W^{2,s}_0(\Omega)]^{\star})}
\leq C\,.
\label{wtbd3}
\end{align}
Next we note that the reflexive Banach space $W^{2,s}_0(\Omega)$ is dense in $L^2(\Omega)$.
It follows that 
$[L^2(\Omega)]^{\star} \equiv L^2(\Omega)$ is continuously embedded and dense in
$[W^{2,s}_0(\Omega)]^{\star}$;
see, for  example, the first two remarks in \S5 in Simon \cite{SimonPres}. 
Furthermore, we have that $C(\overline{\Omega})$ is continuously embedded and dense in
$[W^{2,s}_0(\Omega)]^{\star}$.
Hence, on recalling the compact embedding of $W^{1,s}(\Omega)$ into 
$C(\overline{\Omega})$ for $s>d$ and that $M_\eps: C(\overline{\Omega})
\rightarrow C(\overline{\Omega})$, we obtain from    
(\ref{wtbd3}), (\ref{compact1}) and (\ref{Mepscont}) the strong convergence results (\ref{wscont},e).
It follows from (\ref{Sobinterp})
for $s>d$ and $\alpha(s,d) \in (0,1)$ that
\begin{align}
&\|w^{\tau}-w^{\tau,\pm}\|_{L^2(0,T;C(\overline{\Omega}))}^2
\nonumber \\
& \qquad \qquad \leq C(T,\|w^{\tau(,\pm)}\|_{L^\infty(0,T;W^{1,s}(\Omega))})
\,\|w^{\tau}-w^{\tau,\pm}\|_{L^2(0,T;L^2(\Omega))}^{2(1-\alpha)}.
\label{wdifft}
\end{align} 
Therefore the strong convergence results (\ref{wsL2},f) follow immediately from
(\ref{wdifft}), (\ref{wtbd2}), (\ref{wdiff}), (\ref{wscont}) and (\ref{Mepscont}). 

On noting (\ref{wtwcont},g), Assumption (A3) and (\ref{ftauconv}), we can pass to the limit 
$\tau \rightarrow 0$ for the subsequences in (\ref{Q1tau}) to obtain (\ref{Q1}). 

We now consider passing to the limit $\tau \rightarrow 0$ for the subsequences in 
(\ref{Q2tau}), where at first we fix $\vv \in C^\infty(0,T; [C^\infty(\overline{\Omega})]^d)$.
Noting (\ref{wscont}--g), Assumption (A3) and (\ref{ftauconv}) we immediately obtain 
(\ref{Q2}) for the fixed $\vv \in C^\infty(0,T; [C^\infty(\overline{\Omega})]^d)$.
The only term that requires some comment is the one involving $\vq^{\tau,+}$,
which, similarly to (\ref{Qth2conc})--(\ref{Qth2cond}), we now discuss. First we note that
\begin{align}
&\int_{0}^T \langle |\vq^{\tau,+}|, M_\epsilon (w^{\tau,+}) \rangle_{C(\overline{\Omega})}
\,{\rm d}t
\nonumber \\
& \quad = \int_{0}^T \langle |\vq^{\tau,+}|, M_\epsilon (w) \rangle_{C(\overline{\Omega})}
\,{\rm d}t
+
\int_{0}^T \langle |\vq^{\tau,+}|, M_\epsilon (w^{\tau,+}) - M_\epsilon(w) \rangle_{C(\overline{\Omega})}
\,{\rm d}t.
\label{T1+T2}
\end{align} 
As $M_\epsilon(w) \in C([0,T],C(\overline{\Omega}))$ is positive, 
it follows from (\ref{qwcont}) and (\ref{Mweak2}) that
\begin{align}
\liminf_{\tau \rightarrow 0}
\int_{0}^T \langle |\vq^{\tau,+}|, M_\epsilon (w) \rangle_{C(\overline{\Omega})}
\,{\rm d}t
\geq
\int_{0}^T \langle |\vq|, M_\epsilon (w) \rangle_{C(\overline{\Omega})}
\,{\rm d}t.
\label{T1}
\end{align}
It follows from (\ref{wtbd2}) and (\ref{MwsL2}) that
\begin{align}
&\lim_{\tau \rightarrow 0 }\left|\int_{0}^T \langle |\vq^{\tau,+}|, 
M_\epsilon (w^{\tau,+}) - M_\epsilon(w) \rangle_{C(\overline{\Omega})}
\,{\rm d}t\right| 
\nonumber \\
& \hspace{1in} \leq \lim_{\tau \rightarrow 0 } C\,
\|M_\epsilon (w^{\tau,+}) - M_\epsilon(w)\|_{L^2(0,T;C(\overline{\Omega}))}
=0.
\label{T2}
\end{align}
Combining (\ref{T1+T2})--(\ref{T2}) yields the desired result.
Finally, we obtain the desired result (\ref{Q2}) by noting that
any $\vv \in L^1(0,T;{\vV}^{\cal M}(\Omega))$ can be approximated
by a sequence $\{\vv_j\}_{j \geq 1}$, with $\vv_j \in   
C^\infty(0,T; [C^\infty(\overline{\Omega})]^d)$, on recalling (\ref{weakj}--c);
and that all the terms in (\ref{Q2}) are well-defined.
Hence we have shown that $\{w,\vq\}$ solves (Q), (\ref{Q1},b).

We now show that $w$ solves (P), (\ref{P}). 
Choosing $\vv=\vzero$ in (\ref{Q2}) yields that
\begin{align}
&-\int_0^T \langle |\vq|,M_\epsilon(w)\rangle_{C(\overline{\Omega})} \,{\rm d}t  
\;\geq\; 
-\int_0^T  
(f,w)\,{\rm d}t 
+ \tfrac{1}{2} \left[\, |w(\cdot,T)|^2_{0,\Omega} - |w^\epsilon_0(\cdot)|^2_{0,\Omega}
\,\right]\,.
\label{Q2v=0}
\end{align}
Then for any $\eta \in L^1(0,T;K(w) \cap C^1_0(\overline{\Omega}))$,
on noting (\ref{Q2v=0}),
we have that
\begin{align}
\int_0^T \langle \vq ,\del \eta \rangle_{C(\overline{\Omega})} \,{\rm d}t  
& \geq 
-\int_0^T \langle |\vq|,M_\epsilon(w)\rangle_{C(\overline{\Omega})} \,{\rm d}t  
\nonumber \\
& \geq 
-\int_0^T  
(f,w)\,{\rm d}t 
+ \tfrac{1}{2} \left[\, |w(\cdot,T)|^2_{0,\Omega} - |w^\epsilon_0(\cdot)|^2_{0,\Omega}
\,\right]\,.
\label{MQ2v=0}
\end{align}
Using the relationship (\ref{MQ2v=0}) in (\ref{Q1}), we obtain (\ref{P}).
Finally, we need to show that $w \in L^\infty(0,T;K(w))$, as opposed to just 
$w \in L^\infty(0,T;W^{1,\infty}_0(\Omega))$.  
It follows from (\ref{Q2}) that
\eqlabstart
\begin{align}
&\int_0^T \left[ \langle |\vq|,M_\epsilon(w)\rangle_{C(\overline{\Omega})} - 
(f,w)\right] \,{\rm d}t  + \tfrac{1}{2} \left[\, |w(\cdot,T)|^2_{0,\Omega} - |w^\epsilon_0(\cdot)|^2_{0,\Omega}
\, \right]
\nonumber \\ 
& \hspace{2.5in} \leq {\cal J} :=
\inf_{\vv\in L^1(0,T;\vV^{\cal M}(\Omega))} J(\vv)\,, 
\label{Q2infa}
\end{align}
where
\begin{align}
J(\vv):=
\int_0^T \left[ \langle |\vv|,M_\epsilon(w)\rangle_{C(\overline{\Omega})} - 
(\del\,.\,\vv,w)\right] \,{\rm d}t\,.
\label{Q2infb}
\end{align}
\eqlabend
Choosing $\vv=\vzero$ yields that ${\cal J} \leq 0$.
If ${\cal J}<0$ then, for any minimizing sequence $\{ \vv_j\}_{j \geq 1}$, we
obtain that $J(2\,\vv_j)=2\,J(\vv_j)\rightarrow 2\,{\cal J}<{\cal J}$, which is a
contradiction. Hence ${\cal J}=0$, and so we have that
$J(\vv) \geq 0$ 
for any  $\vv \in L^1(0,T;{\vV}^{\cal M}(\Omega))$.
Since this is true also for $-\vv$, 
and as $w \in L^\infty(0,T;W^{1,\infty}_0(\Omega))$,
we obtain that
\begin{align*}
&\left| \int_0^T ( \vv, \del w) \,{\rm d}t \right|
= \left| \int_0^T ( \del\,.\,\vv, w) \,{\rm d}t \right|
\leq \left| \int_0^T (|\vv|,M_\epsilon(w))\,{\rm d}t \right|
\nonumber \\
& \hspace{2.8in}\qquad
\forall \ \vv \in L^1(0,T;W^{1,1}(\Omega))
\,;
\end{align*}
and therefore by a density result that 
\begin{align}
&\left| \int_0^T ( \vv, \del w) \,{\rm d}t \right|
\leq \left| \int_0^T (|\vv|,M_\epsilon(w))\,{\rm d}t \right|
\qquad
\forall \ \vv \in L^1(0,T;L^1(\Omega))
\,.
\label{Mgradwb}
\end{align}
For any $p \in (2,\infty)$,  choosing $\vv = |[M_\eps(w)]^{-1}\,\del w|^{p-2}\, 
[M_\eps(w)]^{-2}\,\del w$ 
in (\ref{Mgradwb}), and noting the continuity of the $p$ norm for $p\in[1,\infty]$,
we obtain that  
\begin{align}
\|[M_\epsilon(w)]^{-1}\,\del w\,\|_{L^\infty(0,T;L^\infty(\Omega))} \leq 1.
\label{Mgradw}
\end{align}
Hence we have that $w \in L^\infty(0,T;K(w))$, and so $w$ solves (P), (\ref{P}).
\endproof  

\begin{rem}
\label{wtL1rem}
Under only Assumptions (A1) and (A2), we obtain in place of the second and fourth bounds
in (\ref{wtbd2}) that 
\begin{align}
\left \|\frac{\partial w^{\tau}}{\partial t} \right \|_{L^1(0,T;[C^1_0(\overline{\Omega})]^{\star})}
+ \|\vq^{\tau,+}\|_{L^1(0,T;[{\cal M}(\overline{\Omega})]^{d})}
\leq C\,.
\label{wtbd4}
\end{align}
The second bound in (\ref{wtbd4}) follows from (\ref{wnbasic}) and (\ref{Mepsineq}), whilst the
first follows from (\ref{wtbd1}) and the second bound in (\ref{wtbd4}).   
Unfortunately, the first bound in (\ref{wtbd1}) is not enough to obtain strong convergence
of $w^{\tau}$ using (\ref{compact1}), as we require $\alpha >1$.
Hence, the need for the additional Assumption (A3).

We believe the assumption $w_0^\epsilon \geq 0$ in (A3) is not really required
to prove (\ref{vqnmax}), and the assumption 
$\del w^\epsilon_0 \,. \, \underline{\nu} <k_0$ on $\partial \Omega$
in (A1) should be sufficient. 
For $n=1 \rightarrow N$, as $w^n - w^0 \in C_0(\overline{\Omega})$ is nonnegative,  
it follows from (A1) that 
$\del w^n \,. \, \underline{\nu} <k_0$ on $\partial \Omega$. 
Formally, $\vq^n = - \lambda^n \, \del w^n$ in $\Omega$ with $\lambda^n 
\geq 0$, and  as $\vq^n = \vzero$ on subcritical slopes,
this yields that $\int_{\Omega} \del\,.\, \vq^n\,{\rm d}\vx
= \int_{\partial \Omega} \vq^n\,.\, \underline{\nu} \, {\rm d}s 
= - \int_{\partial \Omega}
\lambda^n \, \del w^n  \,.\, \underline{\nu} \, {\rm d}s
\geq 0$. This can then be exploited in (\ref{vqnmax}) by noting that
\begin{align*}
\langle |\vq^n|, M_\epsilon(w^n) \rangle_{C(\overline{\Omega})} 
&\leq (\del\,.\,\vq^n,w^n+{\mathfrak M})=
\left( f^n - \frac{w^n-w^{n-1}}{\tau_n}, w^n + {\mathfrak M}\right)
\nonumber \\
& \leq |f^n|_{0,\Omega}\,|w^n+ {\mathfrak M}|_{0,\Omega} \leq C\,, 
\end{align*}  
where ${\mathfrak M} = \max_{n=1 \rightarrow N} \|w^n\|_{0,\infty,\Omega} \leq C.$
Unfortunately, we are not able to make rigorous the formal argument above
establishing that $\int_{\Omega} \del\,.\, \vq^n\,{\rm d}\vx \geq 0$
under Assumption (A1).
\end{rem}

\section{The Nonlinear Algebraic Systems}
\label{snas}
\setcounter{equation}{0}

\subsection{Solution of (Q$^{h,\tau}_A$)}
\label{QAalg}

To solve the nonlinear algebraic system arising from (Q$^{h,\tau}_A$), we recall Theorem 
\ref{Qthstab} and use an extension of the splitting algorithm, ALG2, see page 
170 in Glowinski \cite{Glow} from the variational to the quasivariational case.    
We introduce the Lagrangian ${\cal L}^{h,n} : U^h_0 \times \vS^h \times \vS^h \rightarrow
{\mathbb R}$ defined by
\begin{align}
{\cal L}^{h,n}(\eta^h,\vpsi^h,\vv^h) :=
E^{h,n}(\eta^h) - (\vv^h,\del \eta^h - \vpsi^h)\,,
\label{cLag}
\end{align} 
where $E^{h,n}(\cdot)$ is defined by (\ref{KminE}).
For a given $\rho \in {\mathbb R}_{>0}$, we introduce the augmented 
Lagrangian ${\cal L}^{h,n}_{\rho} : U^h_0 \times \vS^h \times \vS^h \rightarrow
{\mathbb R}$ defined by
\begin{align}
{\cal L}^{h,n}_{\rho}(\eta^h,\vpsi^h,\vv^h) :=
{\cal L}^{h,n}(\eta^h,\vpsi^h,\vv^h) + \frac{\rho}{2}\,
|\del \eta^h - \vpsi^h|^2_{0,\Omega}\,.
\label{cLagr}
\end{align} 
For any $\chi^h \in U^h_0$, we introduce the closed convex non-empty set  
\begin{align}
\vR^h(\chi^h) &:= \{ \vpsi^h \in \vS^h : |\vpsi^h| \leq M^h_\epsilon(P^h \chi^h) 
\quad \mbox{a.e.\ on } \Omega\}\,.   
\label{Rh}
\end{align}

Set $W_A^{n,0} = W^{n-1}_A\in U^h_0$,  
$\vphi_A^{n,0} = \vphi_A^{n-1} \in \vS^h$
and
$\vQ_A^{n,0} = \vQ_A^{n-1} \in \vS^h$, where we choose $ \vphi_A^0 =\vQ_A^0 =\vzero$. 

For $m\geq 1$, given iterates $W_A^{n,m-1} \in U^h_0$,   
$\vphi^{n,m-1}_A \in \vS^h$
and $\vQ_A^{n,m-1} \in \vS^h$, then
\eqlabstart
\begin{itemize}
\item[(i)] Find $W^{n,m}_A \in U^h_0$ such that
\begin{align}
{\cal L}^{h,n}_{\rho}(W^{n,m}_A,\vphi^{n,m-1}_A,\vQ^{n,m-1}_A) 
\leq {\cal L}^{h,n}_{\rho}(\eta^h,\vphi^{n,m-1}_A,\vQ^{n,m-1}_A)
\quad \forall \ \eta^h \in U^h_0\,.
\label{Wnm}
\end{align}
\item[(ii)] Find $\vphi^{n,m}_A \in \vR^h(W^{n,m}_A)$ such that
\begin{align}
{\cal L}^{h,n}_{\rho}(W^{n,m}_A,\vphi^{n,m}_A,\vQ^{n,m-1}_A) 
&\leq {\cal L}^{h,n}_{\rho}(W^{n,m}_A,\vpsi^{h},\vQ^{n,m-1}_A)
\nonumber \\
& \hspace{1in}\quad \forall \ \vpsi^h \in \vR^h(W^{n,m}_A)\,.
\label{phinm}
\end{align}
\item[(iii)] 
Set
\begin{align}
\vQ^{n,m}_A = \vQ^{n,m-1}_A - \rho\,(\del W^{n,m}_A-\vphi^{n,m}_A)\,.
\label{Qnm}
\end{align}
\end{itemize}
\eqlabend

Step (i) is equivalent to finding the unique $W^{n,m}_A \in U^h_0$ solving the linear problem
\begin{align}
&\left(\frac{W^{n,m}_A-W^{n-1}_A}{\tau_n},\eta^h\right)^h + 
\rho\,(\del W^{n,m}_A-\vphi^{n,m-1}_A,\del \eta^h) 
- (\vQ^{n,m-1}_A,\del \eta^h) 
\nonumber \\
&\hspace{2.8in}= (f^n,\eta^h)
\qquad \forall \ \eta^h \in U^h_0\,.
\label{Wnmlin}
\end{align} 
Step (ii) decouples to solving the problem on each element $\sigma \in {\cal T}^h$.
Let ${\cal L}^{h,n}_{\rho,\sigma} : U^h_0 \mid_\sigma \times {\mathbb R}^d \times  
{\mathbb R}^d \rightarrow {\mathbb R}$ be such that
\begin{align}
{\cal L}^{h,n}_{\rho} (\eta^h,\vpsi^h,\vv^h) = \sum_{\sigma \in {\cal T}^h}
{\cal L}^{h,n}_{\rho,\sigma} (\eta^h_{\sigma},\vpsi^h_{\sigma},\vv^h_{\sigma}) 
\,,
\label{calLrsig}
\end{align} 
where the subscript $\sigma$ denotes restriction to the element $\sigma$.
Hence, for all $\sigma \in {\cal T}^h$, first find $\widehat \vphi^{n,m}_{A,\sigma} 
\in {\mathbb R}$
such that
\begin{align}
{\cal L}^{h,n}_{\rho,\sigma} (W^{n,m}_{A,\sigma},
\widehat \vphi^{n,m}_{A,\sigma},\vQ^{n,m-1}_{A,\sigma})
\leq {\cal L}^{h,n}_{\rho,\sigma} (W^{n,m}_{A,\sigma},
\va,\vQ^{n,m-1}_{A,\sigma}) \quad \forall   \ \va \in {\mathbb R}^d\,,
\label{calLsigmin}
\end{align}
then project $\widehat \vphi^{n,m}_{A,\sigma}$ to the ball of radius 
$[M^h_\epsilon(P^h W^{n,m}_A)]_\sigma$ centred at the origin
to yield $\vphi^{n,m}_{A,\sigma}$. The minimization (\ref{calLsigmin}) leads
to
\eqlabstart
\begin{align}
\widehat \vphi^{n,m}_{A,\sigma} = \frac{ - \vQ^{n,m-1}_{A,\sigma} + \rho\, \del W^{n,m}_A \mid_{\sigma}}
{\rho}\,; 
\label{hatphinm}
\end{align}
and we then set
\begin{align}
\vphi_{A,\sigma}^{n,m} = \left\{\begin{array}{ll}
\widehat \vphi_{A,\sigma}^{n,m} &\mbox{if } |\widehat \vphi_{A,\sigma}^{n,m}| 
\leq [M^h_\epsilon(P^h W^{n,m}_A)]_\sigma\,,\\[3mm]
\displaystyle \frac{\widehat \vphi_{A,\sigma}^{n,m}}
{|\widehat \vphi_{A,\sigma}^{n,m}|} \;  [M^h_\epsilon(P^h W^{n,m}_A)]_\sigma
\qquad & \mbox{otherwise.} \end{array}
\right.
\label{phinmex}
\end{align}
\eqlabend
So overall, $\vphi^{n,m}_A \in \vR^h(W^{n,m}_A)$ is such that
\begin{align}
( \rho\,(\vphi^{n,m}_A - \del W^{n,m}_A) + \vQ^{n,m-1}_A, \vpsi^h - \vphi^{n,m}_A) \geq 0
\qquad \forall \ \vpsi^h \in \vR^h(W^{n,m}_A)\,.
\label{vphinmineq}
\end{align}

On noting that $\vphi^n_A = \del W^n_A$, $n=1 \rightarrow N$,
in the variational case, $M_\eps(\cdot) \equiv k_0 \in {\mathbb R}_{>0}$, then following 
the abstract framework in \S 5.1 in \cite{Glow} one can show that for $n=1 \rightarrow N$ 
and $m \geq 1$ that 
\begin{align}
&\left[ |\widetilde \vQ^{n,m-1}_A|^2_{0,\Omega} + \rho^2\,|\widetilde 
\vphi^{n,m-1}_A|^2_{0,\Omega} \right]
- \left[ |\widetilde \vQ^{n,m}_A|^2_{0,\Omega} + \rho^2\,|\widetilde 
\vphi^{n,m}_A|^2_{0,\Omega} \right]
\nonumber \\
& \qquad \geq \frac{2\rho}{\tau_n} |\widetilde W^{n,m}_A|^2_h
+ \rho^2\, | \del \widetilde W^{n,m}_A - \widetilde \vphi^{n,m}_A|_{0,\Omega}^2
+ \rho^2 \,|\widetilde \vphi^{n,m}_A  - \widetilde \vphi^{n,m-1}_A|^2_{0,\Omega}\,;
\label{ALG2conv1}
\end{align}
where $\widetilde W^{n,m}_A :=W^{n}_A-W^{n,m}_A$, 
$\widetilde \vphi^{n,m}_A:= \vphi^{n}_A-\vphi^{n,m}_A$ 
and $\widetilde \vQ^{n,m}_A:=\vQ^{n}_A-\vQ^{n,m}_A$. 
Hence, one can deduce 
from (\ref{ALG2conv1}), (\ref{Qnm}), (\ref{Wnmlin}) and (\ref{vphinmineq}) 
for $n=1 \rightarrow N$ that as $m \rightarrow \infty$ 
\begin{align}
W^{n,m}_A \rightarrow W^n_A, \qquad \vphi^{n,m}_A \rightarrow \vphi^n_A = \del W^n_A,\qquad  
\vQ^{n,m}_A \rightarrow \vQ^n_A\,. 
\label{ALG2conv2}
\end{align}

Although we have no convergence proof of the iterative algorithm (\ref{Wnm}--c) 
in the quasi-variational inequality case, in practice it worked well.

\subsection{Solution of (Q$^{h,\tau}_B$)}

Adopting the notation (\ref{gn}), we find $\vQ^n_{B,r}$ solving (\ref{Qthr2Balt}),
and hence $\{W^n_{B,r},$ $\vQ^n_{B,r}\}$ 
solving the $n^{\rm th}$ step of (Q$^{h,\tau}_{B,r}$), 
(\ref{Qthr1B},b), using the following iteration: 

Set $\vQ^{n,0}_{B,r} = \vQ^{n-1}_{B,r}\in \vV^h$. For $m \geq 1$, given iterate 
$\vQ^{n,m-1}_{B,r} \in \vV^h$, find
$\vQ^{n,m}_{B,r} \in \vV^h$ such that
\begin{align}
&(M_\epsilon^h(g^n-\tau_n\,\del\,.\,\vQ^{n,m-1}_{B,r})
\,|\vQ^{n,m-1}_{B,r}|_{\delta}^{r-2} \vQ^{n,m}_{B,r},
\vv^h)^h  + \tau_n\, (\del\,.\,\vQ^{n,m}_{B,r},\del\,.\,\vv^h)
\nonumber \\ 
& \qquad =
(M_\epsilon^h(g^n-\tau_n\,\del\,.\,\vQ^{n,m-1}_{B,r})\left[|\vQ^{n,m-1}_{B,r}|^{r-2}_{\delta} 
-|\vQ^{n,m-1}_{B,r}|^{r-2} \right]\vQ^{n,m-1}_{B,r},
\vv^h)^h
\nonumber \\
& \hspace{2in} + (g^n,\del\,.\,\vv^h) 
\qquad
\forall \ \vv^h\in \vV^h\,,
\label{Qthr2it}
\end{align}
where $|\vv|_{\delta} := (|\vv|^2 +\delta^2)^{\frac{1}{2}}$ with $\delta^2 \ll 1$.
Clearly, the linear system (\ref{Qthr2it}) is well-posed.
We note that similar algorithms have been used in \cite{BP3,BP5}.
Although we have no convergence proof of (\ref{Qthr2it}), in practice it worked well.

\section{Numerical Experiments}
\label{numexpts}
\setcounter{equation}{0}

In this section we perform numerical experiments for our finite element approximations 
(Q$^{h,\tau}_A$), (\ref{Qth1},b), and (Q$^{h,\tau}_{B,r}$), (\ref{Qthr1B},b), as stated in 
Section \ref{fea}; except for ease of implementation 
we replaced $w^\epsilon_0$ and $w^{\epsilon,h}_0$, (\ref{weps0h}), by $w_0$
and $w_0^h = P^h [\pi^h w_0]$, respectively, in $M^h_\eps(\cdot)$, (\ref{Meps_h}),
and in the initial data for both approximations.

The approximation (Q$^{h,\tau}_{A}$) is simpler and is easier to implement than 
(Q$^{h,\tau}_{B,r}$), which is based on the lowest order Raviart--Thomas element.
We refer to \cite{Carst} for a Matlab implementation of 
the lowest Raviart--Thomas element.
Although, 
both approximations lead to an efficient numerical approximation of the evolving 
sand surface, our numerical experiments, see below, show that only the (Q$^{h,\tau}_{B,r}$) 
approximation yields a useful approximation to the surface sand flux for
a reasonable choice of discretization parameters.

The simulations have been performed in Matlab R2011a (64 bit) on a PC with Intel Core i5-2400 
3.10Hz processor with 4Gb RAM. 

In all of our examples, we set the sand internal friction coefficient $k_0=0.4$
and chose $r=1+10^{-7}$ for the approximation (Q$^{h,\tau}_{B,r}$), (\ref{Qthr1B},b).

The stopping criterion for the splitting iterative algorithm, (\ref{Wnm}--c),
 for (Q$^{h,\tau}_A$) 
was chosen as 
\begin{align}
&\frac{  \| W_A^{n,m} - W_A^{n,m-1} \|_{L^1(\Omega)}} { \| W_A^{n,m} \|_{L^1(\Omega)}} 
<10^{-6} \quad \;\mbox{and}\quad \; 
\frac{  \| \underline{\phi}_A^{n,m} - \underline{\phi}_A^{n,m-1} \|_{[L^1(\Omega)]^2}} 
{ \| \underline{\phi}_A^{n,m} \|_{[L^1(\Omega)]^2}} <5\,.\,10^{-4}\,.
\label{stopA}
\end{align}

For solving the nonlinear algebraic system, arising from    
(Q$^{h,\tau}_{B,r}$) at each time level, 
we chose $\delta=10^{-9}$ for the iterative method (\ref{Qthr2it}) with stopping criterion
\begin{align}
\frac{\ds \sum_{e \in {\cal E}^h}|e||Q^{n,m}_{B,r}(e)-Q^{n,m-1}_{B,r}(e)|}
{\ds \sum_{e\in {\cal E}^h} |e||Q^{n,m}_{B,r}(e)|}<3\,.\, 10^{-4}\,.
\label{stopB}
\end{align}
Here ${\cal E}^h$ is the collection of edges associated with the partitioning ${\cal T}^h$
so that any $\vv^h \in \vV^h$ can be written as $\vv^h(\vx) = \sum_{e \in {\cal E}^h} 
v^h(e) \,\vphi_e(\vx)$, where $\{\vphi_e\}_{e \in {\cal E}^h}$ are the standard 
lowest order Raviart--Thomas basis
functions, see \cite{Carst}. 

\subsection{Approximation (Q$^{h,\tau}_{A}$)}
We start with a simple variational inequality example. Let sand be disposed onto a flat, 
$w_0=0$, open circular platform, $\Omega=\{\vx\ :\ |\vx|<1\}$. The source $f$ is 
uniform in its support $|\vx|\leq R_0=0.2$ with $\int_{\Omega} f(\vx,t) \,{\rm d} \vx =1$
for all $t \geq 0$. 
Due to the radial symmetry, 
the analytical solution to this problem, $\{w,\vq\}$, is easy to find. The growing pile 
starts as a cut-off cone having critical slopes, volume $t$ and height $t/(\pi\, R_0^2)$. 
Then, at $t^*=\pi\, k_0\,R_0^3\,\sqrt{3}\approx 0.0174$, the pile turns into a cone 
$w(\vx,t)=k_0\,\max(R_c(t)-|\vx|,0)$. This cone grows until its base, a circle of radius 
$R_c(t)=(3t/(\pi \,k_0))^{\frac{1}{3}}$, fills the domain $\Omega$. The flux can be found 
as 
$\vq(\vx,t)=q(|\vx|,t)\,\vx\,/|\vx|$, where $q(R,t)$ is a solution to the balance equation
\beq \frac{1}{R}\,\frac{\partial}{\partial R}\,(R\,q)=f-\frac{\partial w}{\partial t}
\quad \mbox{for } R \in (0,1), 
\quad q(0,t)=0, 
\qquad \mbox{for } t>0.\label{qan}\eeq
 
The iterations of the augmented Lagrangian method with splitting, recall Section \ref{QAalg}, 
converged 
quickly with $\rho=1$. For $t=0.1$ we compared our numerical approximations obtained for
different finite element meshes  and a constant time step $\tau$
with the analytical solution. 
The approximate surfaces, $W^n_A$ with $n\, \tau= t$, were close 
to the exact surface, $w(\cdot,t)$; see Figure \ref{FigVI}, left. 
We checked that, for the meshes employed and  
$\tau \in (0,0.01]$, the error in $w$ 
was dominated by the spatial discretization.  For meshes with maximal 
element sizes $h=$0.01, 0.02, 0.04 the relative errors of $w(\cdot,t)$ 
in the $L^1$ norm were, respectively, 0.3\%, 0.9\%, and 1.8\%.

Although in our simulations the approximate flux iterates $Q^{n,m}_A$ 
also converged on every mesh, 
no pointwise convergence of $Q^n_A$, with $n\,\tau=t$, 
to the analytical flux $\vq(\cdot,t)$ was observed; see Figure \ref{FigVI}, middle.
The approximate flux $Q^n_A$ has a fine structure in the region where the exact flux 
is not zero. There  
elements $\sigma \in {\cal T}^h$
with zero numerical flux were intermixed with elements in which the numerical flux 
was much stronger than the exact one; see Figure \ref{FigVI}, middle and right. 
Such a behavior of the numerical 
solution does not contradict our proof of its vague convergence to $\vq$; but, 
clearly, a different 
method should be employed for approximating the flux even in the variational inequality case.

\begin{figure}[h!]
\begin{center}
\includegraphics[width=8cm,height=6cm]{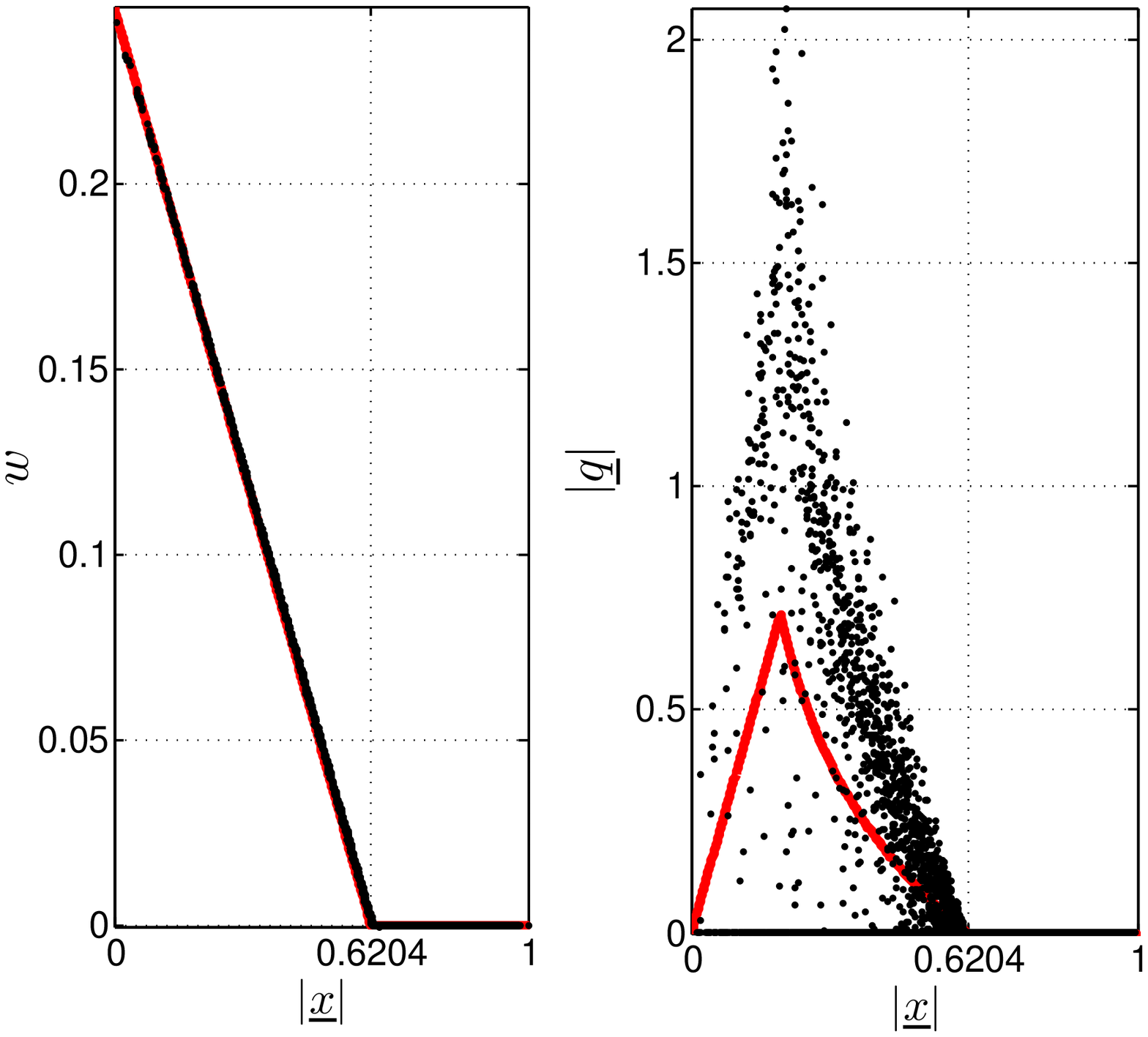}\includegraphics[width=6cm,height=6cm]{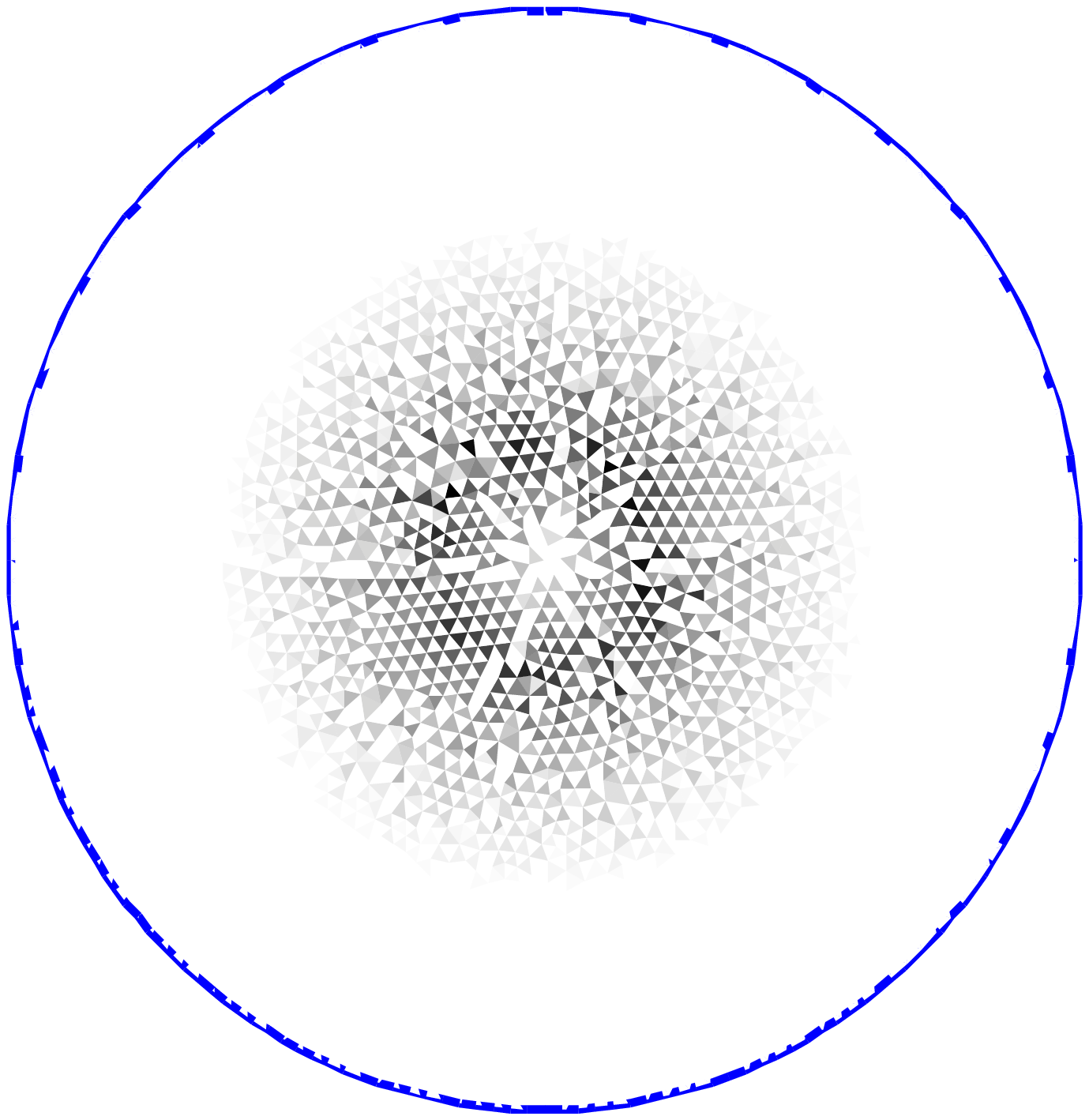}
\end{center}
\caption{Variational inequality, $(Q^{h,\tau}_A)$ approximation, simulation results for  
$h=0.04$, $\tau=0.01,\ t=0.1$. Left -- exact surface $w(|\vx|,t)$ (red line) and its 
approximation $W^n_A$ at the mesh nodes (black dots). Middle -- exact flux modulus 
$|\vq(|\vx|,t)|$ (red line) and its approximation $|\vQ^n_A|$ in the elements (black dots). 
Right -- 
levels of $|\vQ^n_A|$ showing the
mosaic structure of the approximate flux $\vQ^n_A$.}\label{FigVI}
\end{figure}
\begin{figure}[h!]
\centerline{\includegraphics[width=5cm]{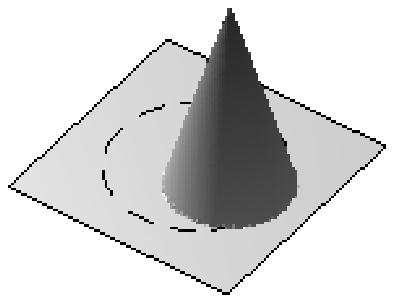}\hspace{-.3cm}
\includegraphics[width=5cm]{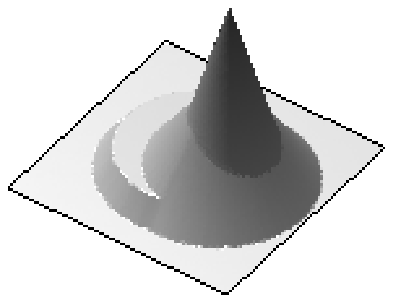}\hspace{-.3cm}
\includegraphics[width=5cm]{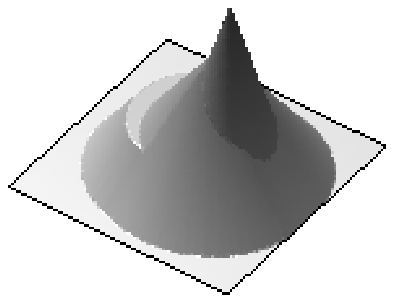}}
\caption{Regularized  quasi-variational inequality with $\epsilon=0.01$, 
$(Q^{h,\tau}_A)$ approximation, 
simulation results for  $h=0.02$, $\tau=0.01$. Left -- the support surface $w_0$, 
the dashed line indicates the boundary of the support of $f$. Middle and Right -- 
approximate sandpile surface, $W^n_A$, corresponding to $t=0.1$ and $t=0.2$, respectively. }
\label{FigVI2}
\end{figure}

The situation is similar for the quasi-variational inequality case, with 
only the evolving pile surface 
being approximated  well using this method. 
In our second example, see Figure \ref{FigVI2},  we set
$w_0=\max(0.5-|\vx-\vx_0|,0)$, where $\vx_0=(0.3,0)$, for the square 
$\Omega=(-1,1)\times(-1,1)$. The source $f$ is uniform in its support 
$|\vx|\leq 0.7$ 
with  $\int_{\Omega}f(\vx,t)=1$ for all $t \geq 0$.
For $h=0.02$ the generated mesh contained approximately 
34,000 elements. Since the gradient constraint 
is now updated after each iteration of the splitting algorithm, existing theory does not 
guarantee its convergence. We found that, for the regularization parameter $\eps=0.01$, good 
convergence of this algorithm is achieved for a smaller value of the augmented Lagrangian
 parameter $\rho$. In this example we chose 
$\rho=0.05$ with stopping criterion (\ref{stopA}) 
and obtained the solution with $\tau=0.005$ and twenty time steps in 
11 minutes of CPU time.

\subsection{Approximation (Q$^{h,\tau}_{B,r}$)}

This approximation performed better in the variational inequality example from the previous 
section, recall Figure \ref{FigVI}. 
Using the time step, $\tau=0.005$, and meshes, $h$=0.02 
and 0.04, 
we obtained, for $t=0.1$, the pile surface with smaller relative errors in the $L^1$ 
norm; 0.1\% and 0.6\%, respectively. Furthermore, for this approximation the fluxes 
$Q^n_{B,r}$ also converged to the exact solution.
Comparing the approximate and exact fluxes at the element centers we estimated the relative 
flux error in the $L^1$ norm. For the two meshes chosen these errors were, correspondingly, 
2.8\% and 5.2\%.

Choosing a different initial support, $w_0(\vx)=\max(0.4-|\vx|,0)$, and keeping the same 
source  $f$ and domain $\Omega$ from the variational inequality example, 
we arrive at a quasi-variational inequality problem that 
can be solved analytically for the unregularized $M(\cdot)$.  
Being discharged from the source, sand now pours down the steep 
conical part of the support surface and forms a pile around this cone. The volume of the pile 
is $t$ and its surface $w(|\vx|,t)=k_0\,(R_2(t)-|\vx|)$ for $|\vx| \in [R_1(t),R_2(t)]$, 
see Figure \ref{qvi_ill}.
Using simple geometric arguments, we first find the two variables, $R_1(t)$ and $R_2(t)$, 
determining this surface from the equations
$$t=\frac{\pi}{3}\left[(R_2^3-R_1^3)\,k_0-(0.4^3-R_1^3)\right]\quad
\mbox{and} \quad 
R_2=R_1+\frac{1}{k_0}(0.4-R_1)\,.$$
We then find the flux using the balance equation (\ref{qan}).

We solved the problem numerically, see Figure \ref{QVI_EX_fig}, with the regularization  
parameter 
$\eps=0.005$ and estimated the errors of $W_{B,r}^n$ and $\vQ_{B,r}^n$ at $t=0.1$ 
using the analytical solution.
As could be expected, the smaller the value of $\eps$,
the more difficult it is to obtain convergence of the iterations (\ref{Qthr2it}) in the 
quasi-variational inequality case. We were, however, able to achieve convergence of 
these iterations by decreasing the time step $\tau$.
For the stated value of $\eps$, we chose $\tau=0.0005$ yielding 200 time steps on 
the time interval $[0,0.1]$ for two different meshes. 
For a mesh generated with $h=0.04$ 
the relative errors in the $L^1$ norm were
0.6\% for the pile surface and 5\% for the surface flux. For  a finer mesh, $h=0.02$,
the corresponding errors were 0.1\% and 2\%.
These results confirm the validity of our regularization, $M_\epsilon(\cdot)$, of $M(\cdot)$.

\begin{figure}[h!]
\psfrag{Y}{$y$}\psfrag{f}{$f$} \psfrag{R1}{{$R_1(t)$}}
\psfrag{R2}{$R_2(t)$} \psfrag{R}{$|\underline{x}|$}
\psfrag{0.4}{0.4}\psfrag{0}{0}
\begin{center}
\includegraphics[width=10.0cm,height=5cm]{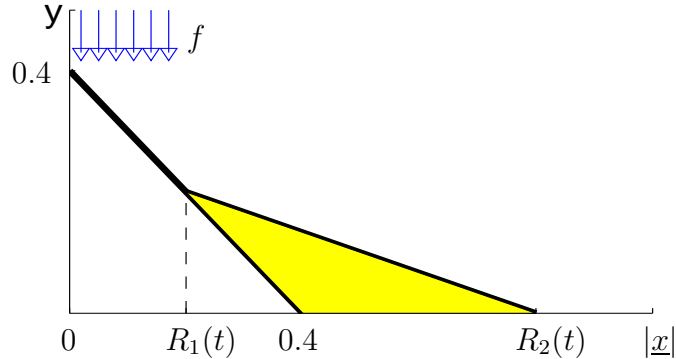}
\end{center}
\caption{Sandpile, yellow region, forming upon the support platform with a steep cone.}\label{qvi_ill}
\end{figure}
\begin{figure}[h!]
\begin{center}
\includegraphics[width=8cm,height=6cm]{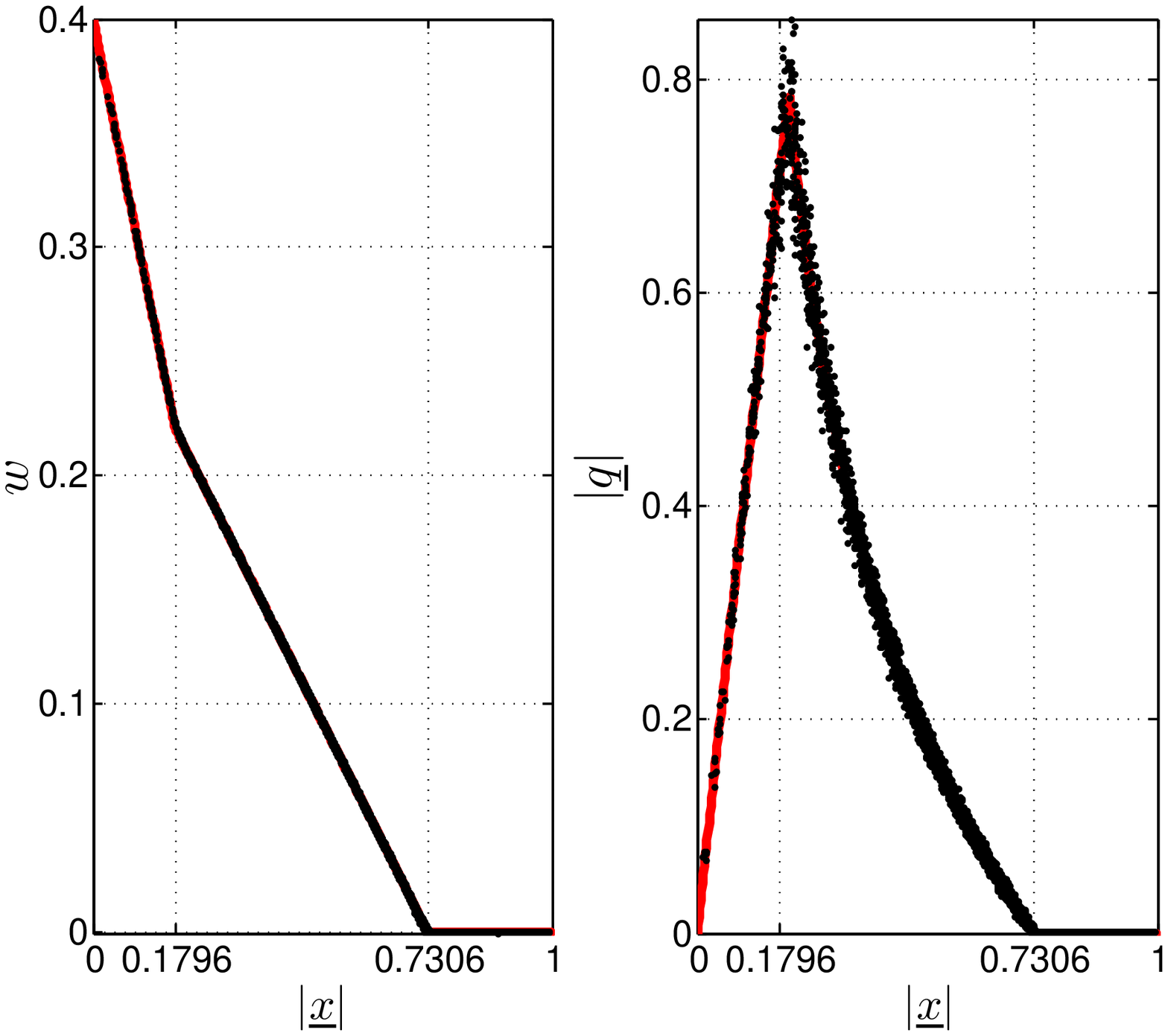}\includegraphics[width=6cm,height=6cm]{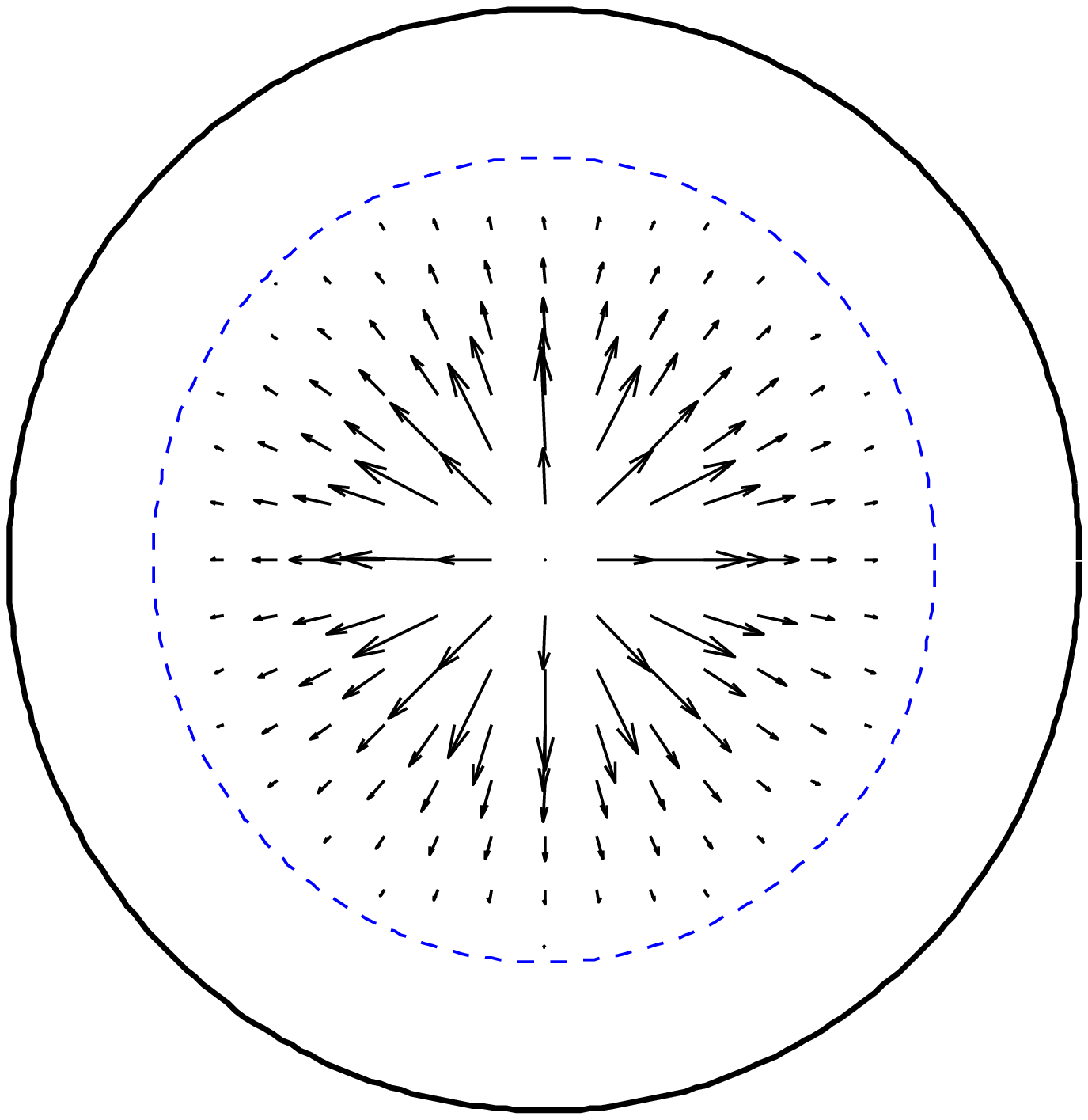}
\end{center}
\caption{Regularized quasi-variational inequality with $\eps=0.005$, 
$(Q^{h,\tau}_B)$ approximation, simulation 
results for  $h=0.04$, $\tau=0.0005$ and $t=0.1$. Left -- exact surface $w(|\vx|,t)$ 
(red line) and its approximation $W^n_{B,r}$ in the elements (black dots). Middle -- 
exact flux modulus $|\vq(|\vx|,t)|$ (red line) and its approximation $|\vQ^n_{B,r}|$ 
at the element centers (black dots). Right -- the $\vQ^n_{B,r}$ vector field, where the
dashed line 
indicates the exact pile boundary.}\label{QVI_EX_fig}
\end{figure}
We solved again, now using 
$(Q^{h,\tau}_{B,r})$, 
the quasi-variational problem 
considered in Figure \ref{FigVI2} above, using the same mesh, time step, and the value of 
regularization parameter. Now we were able to find good approximations, not only to 
the pile surface but the surface 
flux as well, see Figure \ref{QVI_coneB_fig}; and the computation time was about the same.
We note that as the surface $W^n_{B,r}$ touches the support boundary $\partial \Omega$
at some time in the interval $(0.1,0.2)$ sand flows out of the system, which can be seen
from the flux $Q^n_{B,r}$. 

\begin{figure}[h!]
\centerline{\includegraphics[width=5cm]{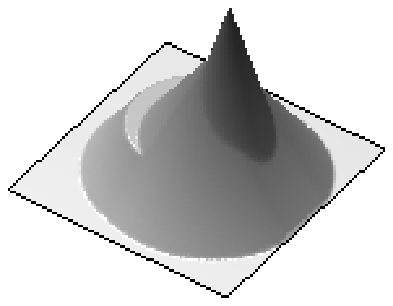}
\includegraphics[width=5cm]{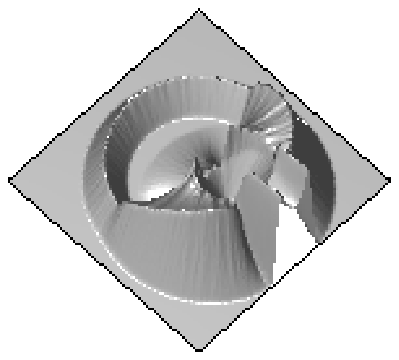}\includegraphics[width=6cm]{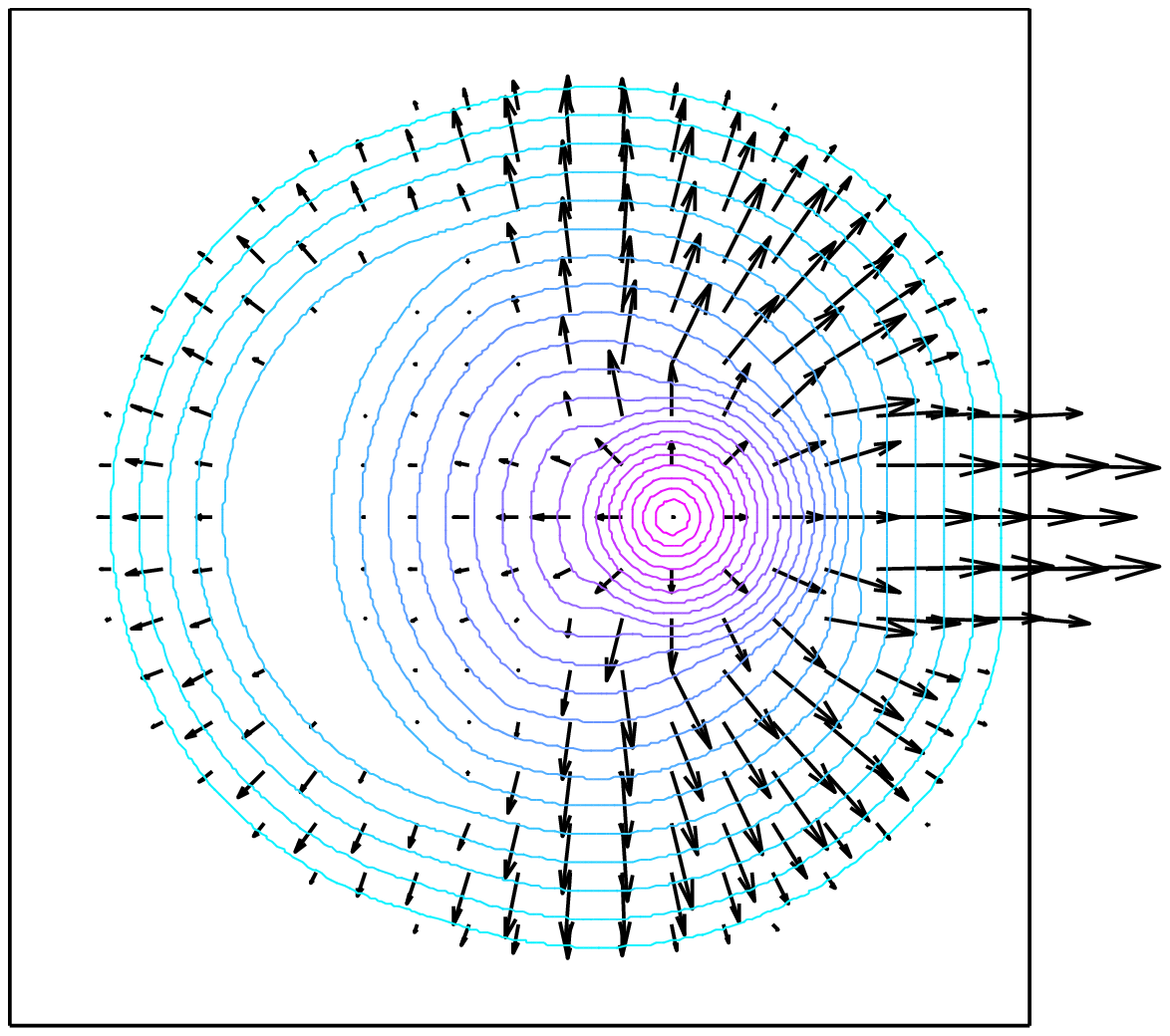}\
}
\caption{Regularized quasi-variational inequality with $\eps=0.01$ 
as in Figure \ref{FigVI2}, 
$(Q^{h,\tau}_B)$ approximation, simulation results for  $h=0.02$, $\tau=0.01$ and  
$t=0.2$. Left -- the calculated surface  $W^n_{B,r}$. Middle -- the flux modulus 
$|\vQ^n_{B,r}|$ at the element centers. 
Right -- the $\vQ^n_{B,r}$ vector field and level contours of $W^n_{B,r}$.}
\label{QVI_coneB_fig}
\end{figure}

In our last example $\Omega=(-1,1)\times(-1,1)$ and $w_0=\min(\max(|x_1|-0.9,|x_2|-0.9),0)$ 
is the surface 
of an inverted pyramid supplemented, to satisfy the no-influx condition (\ref{influx}), 
by a narrow horizontal margin. The uniform source is $f(\vx,t)\equiv 0.25$.
Sand, discharged from the source, flows down the pyramid faces until it reaches a pyramid 
edge; then it pours down along the edge and forms a pile above the apex of the 
inverted pyramid,  see Figure \ref{Fig_pyr}. Our numerical solution clearly shows the 
singularity of the edge fluxes.

\begin{figure}[h!]
\centerline{\includegraphics[width=6cm,height=5cm]{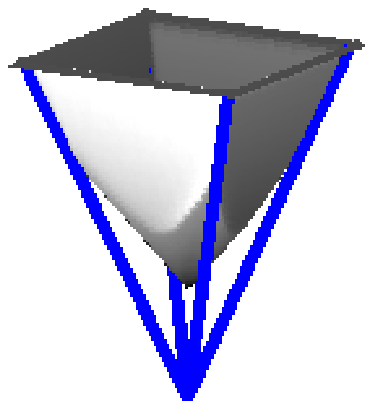}
\includegraphics[width=5cm,height=7cm]{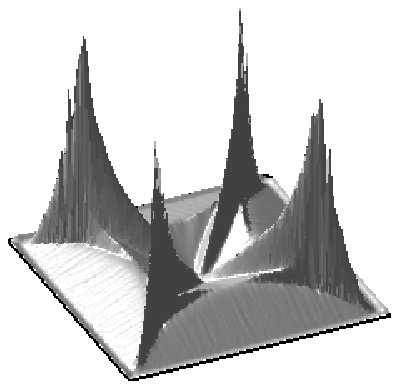}
\includegraphics[width=5cm]{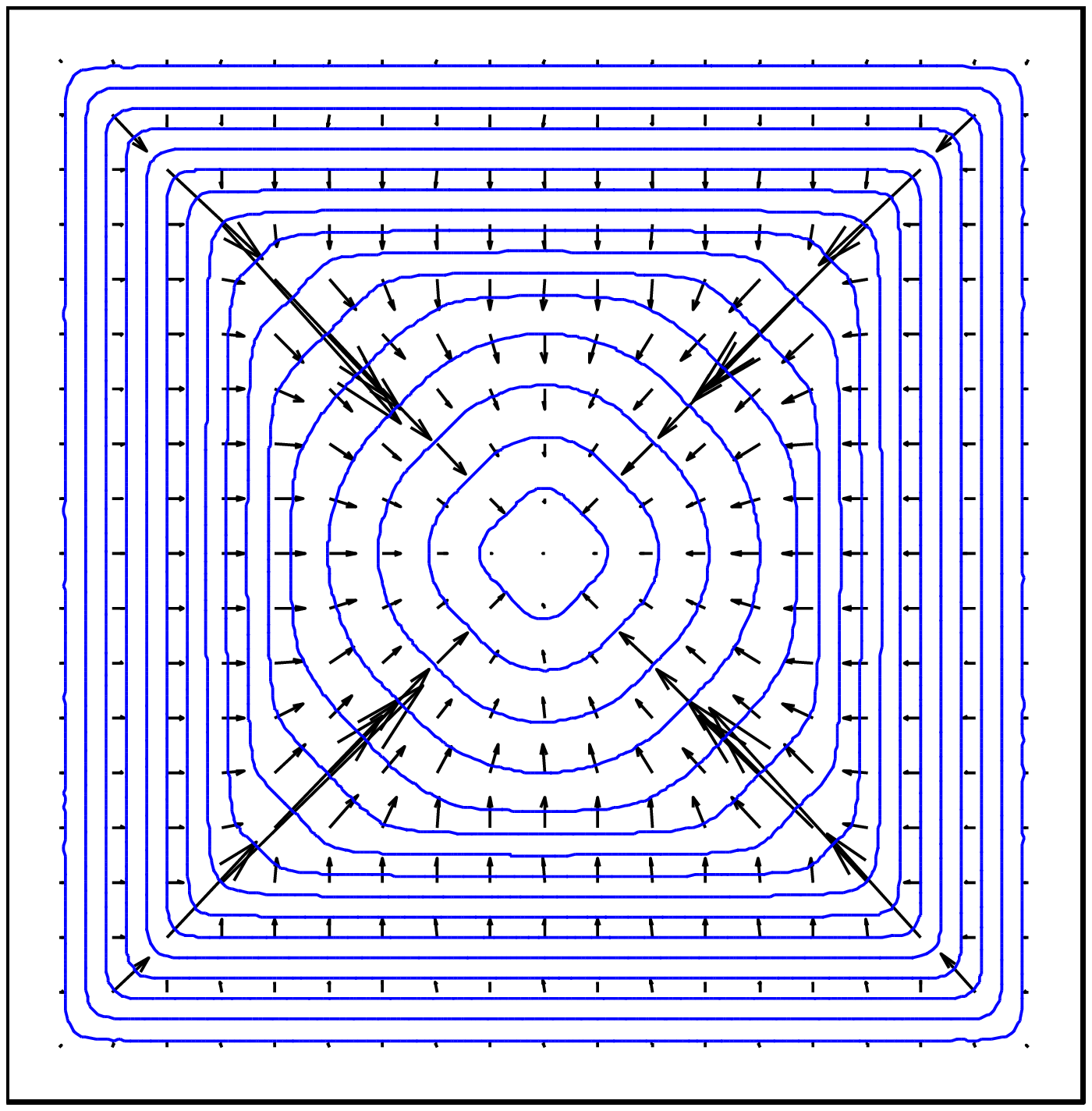}\
}
\caption{Regularized quasi-variational inequality with $\eps=0.02$, 
$(Q^{h,\tau}_B)$ approximation, 
simulation results for  $h=0.02$, $\tau=0.0025$ and $t=0.075$. Left -- 
initial surface $w_0$ (blue lines) and the approximate surface $W^n_{B,r}$ in the 
elements (grey surface). Middle -- the flux modulus $|\vQ^n_{B,r}|$ at the element centers. 
Right --  $\vQ^n_{B,r}$ vector field and levels of $W^n_{B,r}$.}\label{Fig_pyr}
\end{figure}


\end{document}